\documentclass[11pt,letterpaper,english]{article}
\usepackage[latin9]{inputenc}
\usepackage{babel}
\usepackage{mathrsfs}
\usepackage{amsmath}
\usepackage{amssymb}
\usepackage{graphicx}
\usepackage{setspace}
\usepackage[authoryear]{natbib}
\usepackage[unicode=true,pdfusetitle,
 bookmarks=true,bookmarksnumbered=false,bookmarksopen=false,
 breaklinks=false,pdfborder={0 0 0},pdfborderstyle={},backref=false,colorlinks=false]
 {hyperref}
\usepackage[dot]{bibtopic}

\makeatletter

\usepackage{sectsty}
\usepackage{pslatex}
\usepackage{babel}
\usepackage{chngcntr}
\usepackage{apptools}
\usepackage{bm}
\usepackage{amsmath}

\usepackage{endnotes}
\let\footnote=\endnote

\AtAppendix{\counterwithin{lemma}{section}}
\AtAppendix{\counterwithin{proposition}{section}}
\AtAppendix{\counterwithin{theorem}{section}}
\AtAppendix{\counterwithin{corollary}{section}}
\AtAppendix{\counterwithin{assumption}{section}}
\AtAppendix{\counterwithin{definition}{section}}
\AtAppendix{\counterwithin{equation}{section}}

\sectionfont{\sffamily\large}
\subsectionfont{\sffamily\normalsize}
\subsubsectionfont{\sffamily\small}

\newcounter{eqnum}[section] 
\newtheorem{proposition}{Proposition}
\newtheorem{theorem}{Theorem}
\newtheorem{lemma}{Lemma} 

\newtheorem{assumption}{Assumption}

\newcommand{\eproof}{\mbox{}\hfill{\rule{8pt}{8pt}}}

\makeatother

\begin{document}
\title{Investment in the common good: free rider effect and the stability
of mixed strategy equilibria}
\author{Youngsoo Kim\thanks{Culverhouse College of Business, University of Alabama, Tuscaloosa,
Alabama 35487\protect \\
 Email: ykim@cba.ua.edu}, $\text{ }$ H. Dharma Kwon\thanks{Gies College of Business, Department of Business Administration, University
of Illinois at Urbana-Champaign, Champaign, Illinois 61820\protect \\
 Email: dhkwon@illinois.edu}}
\date{January 17, 2022}
\maketitle
\begin{abstract}
In the game of investment in the common good, the free rider problem
can delay the stakeholders' actions in the form of a mixed strategy
equilibrium. However, it has been recently shown that the mixed strategy
equilibria of the stochastic war of attrition are destabilized by
even the slightest degree of asymmetry between the players. Such extreme
instability is contrary to the widely accepted notion that a mixed
strategy equilibrium is the hallmark of the war of attrition. Motivated
by this quandary, we search for a mixed strategy equilibrium in a
stochastic game of investment in the common good. Our results show
that, despite asymmetry, a mixed strategy equilibrium exists if the
model takes into account the repeated investment opportunities. The
mixed strategy equilibrium disappears only if the asymmetry is sufficiently
high. Since the mixed strategy equilibrium is less efficient than
pure strategy equilibria, it behooves policymakers to prevent it by
promoting a sufficiently high degree of asymmetry between the stakeholders
through, for example, asymmetric subsidy. \\
\smallskip{}
\end{abstract}
\vspace{0.5in}

\newpage{}

\begin{btUnit}

\section{Introduction\label{sec:Intro}}

A free rider problem commonly arises in a game of investment in the
common good. For example, if an electronics retailer such as \emph{Best
Buy} offers corporate social responsibility (CSR) programs such as
recycling, tech training, and supplier audit at its own expense, the
other firms in the industry also benefit from them: recycling reduces
material costs across the industry; technology training for the disadvantaged
expands the market coverage for the industry; investments in CSR increase
the industry-wide stock returns \citep{Gunther2015,Serafeim2018}.
A free rider problem also arises in commodity advertising in the orange
juice industry \citep{Lee1988} and the salmon industry \citep{Kinnucan2003}:
one firm's investment in an advertising campaign can benefit its competitors.
The free rider problem can delay the stakeholders' investment actions
and hence diminish the level of the common good. Therefore, it behooves
the policymakers to mitigate the free rider effect. The goal of this
paper is to study a stochastic game of investment in the common good
and examine the impact of asymmetry on the free rider effect.

The game of investment in the common good is a special case of a concession
game. The free rider effect in a concession game often manifests itself
as a mixed strategy equilibrium in which the players delay their actions
by a random time. However, the mixed strategy equilibrium may be destabilized
by even the smallest degree of asymmetry between the players in a
stochastic environment \citep{Georgiadis2019}. This suggests that
even an extremely small degree of asymmetry may mitigate the free
rider effect. This result is surprising because it is contrary to
the commonly accepted notion that a mixed strategy equilibrium is
the hallmark of the games of concession. Furthermore, from the practical
point of view, such extreme sensitivity to asymmetry may not be an
accurate depiction of real life decision making. Motivated by this
quandary, this paper searches for mixed strategy equilibria in the
context of investment in the common good. In particular, we ask the
following questions: Is there a mixed strategy equilibrium in the
model that we study? If so, under what conditions do they exist? Lastly,
what is the managerial implication to the policymakers?

To answer the proposed questions, we analyze a parsimonious model
of a game between two players whose payoffs are under the influence
of a state variable that models a stochastic stock of the common good.
Each player has an infinite number of opportunities to invest in the
common good. The timing and the size of the investment are both under
the player's discretion, but each investment requires a fixed upfront
cost as well as a variable cost. Such a cost structure is typical
in the investment in the common good; for example, a CSR program and
an advertising campaign generally require a fixed upfront cost which
stems from organizing a new team, drawing the plans for execution
of the project, and purchasing equipment.

The game-theoretic model that we examine has a broad class of applications.
In the example of commodity advertising, the state variable is the
stock of goodwill on the commodity in question, and the investment
is on advertising campaign to boost the goodwill \citep{Nerlove1962,Lon2011,Reddy2016}.
In the example of investments in CSR, the players are firms in the
electronics retail industry, and their industry-wide profitability
can be considered as the state variable \citep{Gunther2015,Serafeim2018}.

The paper makes both managerial and methodological contributions.
First, the paper sheds light on the impact of the repeated opportunities
of investment. In contrast to the one-shot stochastic war of attrition,
our model with repeated opportunities of investment possesses a mixed
strategy equilibrium under a moderate degree of asymmetry. We therefore
resolve the discrepancy between the results from the \emph{stochastic}
war of attrition and those from the \emph{canonical} war of attrition.
The result provides useful insights to policymakers: a sufficiently
high degree of asymmetry between the investors can mitigate the free
rider effect by destabilizing the mixed strategy equilibrium.

Second, the paper presents an equilibrium solution to a novel class
of impulse control games possessing a free rider effect. In particular,
it provides a verification theorem for this class of impulse control
games and constructs a mixed strategy subgame perfect equilibrium
for the first time in the literature on impulse control games.

The paper is organized as follows. We provide the context of our paper
in the literature and its contributions to the relevant research domains
in Section \ref{sec:Related-Literature}. As a benchmark, we first
examine the single investment game in Section \ref{sec:One-Time}
and reproduce the result of \citet{Georgiadis2019} that the mixed
strategy Markov perfect equilibrium (MPE) is destabilized by even
the smallest degree of asymmetry. The central results of the paper
are presented in Section \ref{sec:Impulse}: we examine the game allowing
for an infinite number of investments, formulate a verification theorem
for the equilibrium, and obtain mixed strategy equilibria in asymmetric
games. We conclude the paper in Section \ref{sec:Conclusions}. All
mathematical proofs can be found in Electronic Companion (EC).

\section{Related Literature \label{sec:Related-Literature}}

The current paper contributes to the literature on games of concession.
One domain of this literature is on the war of attrition and its mixed
strategy equilibrium. In the seminal paper by \citet{Smith1974},
a game of concession is modeled as the war of attrition which is shown
to possess a mixed strategy equilibrium. \citet{Hendricks1988} obtains
all possible mixed strategy equilibria in the continuous-time deterministic
war of attrition. In the stochastic extensions of the war of attrition,
\citet{Murto2004} characterizes pure strategy Markov perfect equilibria,
and \citet{Steg2015} obtains a mixed strategy equilibrium. The mixed
strategy equilibrium has been also empirically found and examined
\citep{Wang2009,Takahashi2015}. However, \citet{Georgiadis2019}
shows that even the slightest degree of asymmetry between the players
destabilizes a mixed strategy Markov perfect equilibrium in the stochastic
war of attrition. \citet{Kwon2019} also shows a similar absence of
inefficient equilibrium in the game of contribution to the common
good under asymmetry. In stark contrast, our paper shows that the
mixed strategy equilibrium is recovered despite asymmetry if there
are repeated investment opportunities.

The control game model which we study is suitable for modeling commodity
advertising decisions, which is known to suffer from the free rider
problem (\citealt{Lee1988}, \citealt{Kinnucan2003}). Advertising
for a product can be modeled as investment in the stock of goodwill
\citep{Nerlove1962}. Expenditure in advertising under a stochastic
environment has been modeled as a stochastic control problem \citep{Sethi1977}.
Specifically, it has been modeled as a singular control problem \citep{Jack2008},
a singular control combined with a discretionary timing of the product
launch \citep{Lon2011}, or as a singular control game in the context
of a market share competition \citep{Kwon2015}. Just as in our paper,
advertising decision has also been modeled as an impulse control problem
\citep{Reddy2016} because an advertising campaign requires substantial
upfront investment.

The model that we examine is also applicable to the investment decisions
in CSR. The literature views CSR as private provision of public goods
by firms and has attempted to delineate the conditions under which
firms contribute to the production of public goods through CSR despite
free-riding incentives. \citet{Bagnoli2003} shows that firms may
strategically engage in CSR activities to outperform their competition
in the product market. \citet{Kotchen2006} studies the green markets
where products with different levels of greenness are available and
examines the conditions under which the production of green products
increases or decreases. \citet{Morgan2019} finds that the shareholders'
preference for public goods can encourage the firms to be more active
in CSR. We complement this stream of literature by focusing on the
free rider effect from the mixed strategy equilibria and identifying
the conditions under which the provision of public goods such as CSR
is hastened or delayed.

Lastly, the theory of stochastic impulse control has long been applied
to various operations research problems including cash management
\citep{Constantinides1978}, inventory management \citep{Ormeci2008,Cadenillas2010},
interest rate intervention \citep{Mitchell2014}, and production capacity
expansion \citep{Bensoussan2019}. The impulse control framework has
also been applied to game-theoretic models. Amongst others, \citet{Stettner1982},
\citet{Cosso2013}, and \citet{ElAsri2018} examine and solve zero-sum
impulse control games. \citet{Dutta1995} solves a two-sided $(s,S)$
game in the context of a duopolist competition. \citet{Basei2019}
and \citet{Guo2019} examine impulse/singular control games with a
large number of players using the mean-field approach. \citet{Ferrari2019}
models the pollution control as a nonzero-sum impulse control game
between the government and the firm and constructs an equilibrium
under certain conditions. \citet{Zabaljauregui2020} proposes a simple
and efficient algorithm for solving symmetric nonzero-sum impulse
control games. \citet{Campi2020} considers an impulse controller-stopper
nonzero-sum game and characterizes two qualitatively different types
of equilibria. \citet{Aid2020} examines a class of nonzero-sum impulse
control games and obtains a general verification theorem. Our paper
contributes to this domain of research by examining a novel class
of impulse control games possessing a free rider effect which yields
a mixed strategy equilibrium.

\section{Single Investment Game \label{sec:One-Time}}

In this section, we examine a single investment game as a benchmark
model. Specifically, we assume that the common good allows for only
\emph{one} opportunity of investment. The single-investment assumption
is unrealistic for the examples of the CSR program or advertising
campaign because the firms naturally possess multiple investment opportunities.
However, it is a convenient assumption as it considerably simplifies
the analysis. A particularly convenient feature of the single investment
game is that it reduces to a stopping game also known as a Dynkin
game; see Proposition \ref{prop:stopping-game}.  Therefore, most
of the results obtained for stopping games (e.g., \citealt{Georgiadis2019},
\citealt{Touzi2002}) translate to our model. One of the main goals
of this paper is to show that this simplification may not be so innocuous
after all; Section \ref{sec:Impulse} shows that we obtain a markedly
different result if this assumption is relaxed.

The central goal of this section is two-fold. First, we illustrate
mixed strategy equilibria of a symmetric game in Section \ref{subsec:Mixed-Strategy-SPE}.
Second, in Section \ref{subsec:Asymmetric-Game}, we turn our attention
to an asymmetric game and search for a mixed strategy equilibrium
that shares the same simple structures as those illustrated in Section
\ref{subsec:Mixed-Strategy-SPE}.

\subsection{The Model and the Equilibrium Concept \label{subsec:The-Model-and}}

This section introduces a model of two players whose payoffs are under
the influence of the common good. At any point in time, each player
can make costly investment to boost the stock of the common good by
any amount, but the common good allows for only one opportunity of
investment. Thus, if an investment is made by one of the players,
the game is terminated immediately.

\subsubsection{State Variable \label{subsubsec:State-Variable}}

We model the level of the common good (state variable) as a stochastic
process. We let $(\Omega,\mathscr{F},\mathbb{F},\mathbb{P})$ denote
a probability space with the filtration $\mathbb{F}=(\mathscr{F}_{t})_{t\ge0}$,
which satisfies the \emph{usual condition} (p. 172, \citealt{Rogers2000}).
In the absence of the players' control, the \emph{uncontrolled} state
variable $X$ is modeled as a \emph{regular} diffusion process (p.
13, \citealt{Borodin1996}) satisfying the following stochastic differential
equation (SDE):
\begin{equation}
dX_{t}=\mu(X_{t})dt+\sigma(X_{t})dW_{t}\:.\label{eq:SDE}
\end{equation}
Here $W=(W_{t})_{t\ge0}$ is a Wiener process adapted to $\mathbb{F}$.
We assume that $\mu(\cdot)<0$ and $\sigma(\cdot)>0$ to model a stochastically
declining state of the common good in the absence of investments.
The process $X$ takes values within the interval $\mathscr{I}=(a,b)\subseteq\mathbb{R}$,
where $a>-\infty$ and $b\le\infty$. For example, in the case of
a geometric Brownian motion, $a=0$ and $b=\infty$. We assume that
the boundary points $a$ and $b$ are \emph{natural} so that $X$
never reaches $a$ or $b$ in finite time (II.6 of \citealp{Borodin1996}).
To ensure that the strong solution to the SDE (\ref{eq:SDE}) exists,
we make the usual assumption of Lipschitz continuity of $\mu(\cdot)$
and $\sigma(\cdot)$ as well as the inequality $\vert\mu(x)\vert+\vert\sigma(x)\vert\le\delta(1+\vert x\vert)$
for some constant $\delta>0$ (Section 5.2, \citealt{Oksendal2003}).
Lastly, we impose the following transversality assumption (p.330,
\citealt{Alvarez2001}) so that the payoff functions are well-behaved
in the limit $t\rightarrow\infty$:
\begin{equation}
\lim_{t\rightarrow\infty}\mathbb{E}^{x}[\vert X_{t}\vert e^{-rt}]=0\quad\forall x\in\mathscr{I}.\label{eq:transversality}
\end{equation}

We now define the $r$-excessive characteristic differential operator
for the diffusion process $X$ \citep{Oksendal2003}:
\begin{equation}
\mathscr{A}:=\frac{1}{2}\sigma^{2}(x)\frac{d^{2}}{dx^{2}}+\mu(x)\frac{d}{dx}-r\:.\label{eq:char-op}
\end{equation}
By the theory of diffusion (Chapter II of \citealp{Borodin1996}),
there are two solutions to the homogeneous differential equation $\mathscr{A}f(x)=0$.
As a matter of notational convention, we let $\phi(\cdot)>0$ and
$\psi(\cdot)>0$ respectively denote the decreasing and the increasing
solution to the equation $\mathscr{A}f(x)=0$.

\subsubsection{Strategy and Payoff \label{subsubsec:Stragey-and-Payoff}}

As a matter of convention, we let $i$ and $j$ denote the indices
of the two players such that $i\neq j$. Given the single-investment
assumption, each player's strategy can be written as $\nu_{i}=(\tau^{(i)},\zeta_{\tau^{(i)}}^{(i)})$
where $\tau^{(i)}$ is player $i$'s stopping time of investment and
$\zeta_{\tau^{(i)}}^{(i)}\geq0$ is the boost in the state variable
measurable with respect to $\mathscr{F}_{\tau^{(i)}}$. Since the
common good is assumed to allow for only one investment, the game
is terminated at the time $\tau^{(1)}\wedge\tau^{(2)}$. Following
the convention from the literature \citep[e.g.,][]{DLR1995,Grenadier1996,Hoppe2001},
we assume that, if both players attempt to invest simultaneously,
i.e., $\tau^{(i)}=\tau^{(j)}$, then only one player's investment
goes through with a probability of 50\%; this accounts for the factor
of 1/2 in the third term in (\ref{eq:single-game-payoff}) below.

The strategy profile $\nu=(\nu_{1},\nu_{2})$ and the tie-breaking
rule described above determine a stopping time $\tau^{\nu}:=\tau^{(i)}\wedge\tau^{(j)}$
at which an investment is made and the boost $\zeta_{\tau^{\nu}}$
made by one of the players. If $\tau^{(i)}<\tau^{(j)}$, then $\zeta_{\tau^{\nu}}=\zeta_{\tau^{(i)}}^{(i)}$;
if $\tau^{(i)}=\tau^{(j)}$, then $\zeta_{\tau^{\nu}}$ can be either
$\zeta_{\tau^{(i)}}^{(i)}$ or $\zeta_{\tau^{(j)}}^{(j)}$ with 50\%
probability each. Hence, we consider a new binary random variable
$\zeta_{\tau^{\nu}}$, taking values of $\zeta_{\tau^{(i)}}^{(i)}$
or $\zeta_{\tau^{(j)}}^{(j)}$ with 50\% probabilities. Under the
strategy profile $\nu$, the \emph{controlled} state variable is given
by
\[
X_{t}^{\nu}=X_{0}+\int_{0}^{t}\mu(X_{s}^{\nu})ds+\int_{0}^{t}\sigma(X_{s}^{\nu})dW_{s}+\zeta_{\tau^{\nu}}\mathbf{1}_{\{\tau^{\nu}\le t\}}\:.
\]

Next, we formulate the payoff to the players with a common discount
rate $r>0$. We let $\pi(X_{t}^{\nu})$ denote the profit flow per
unit time to each player. We assume that $\pi(\cdot)$ is a continuous
and increasing function. When player $i$ makes an investment at time
$\tau^{(i)}$ to boost $X^{\nu}$ by $\zeta_{\tau^{(i)}}^{(i)}\ge0$,
it costs $c_{i}+k\zeta_{\tau^{(i)}}^{(i)}$ to the player, where $c_{i}>0$
is the upfront cost and $k>0$ is the variable cost of the boost.
Based on this specification, we can write player $i$'s payoff at
time $t\le\tau^{\nu}$ conditional on $\mathscr{F}_{t}$ as follows:
\begin{align}
V_{i,t}^{\nu}= & e^{rt}\mathbb{E}^{x}\biggl[\int_{t}^{\infty}\pi(X_{s}^{\nu})e^{-rs}ds-e^{-r\tau^{(i)}}(k\zeta_{\tau^{(i)}}^{(i)}+c_{i})\mathbf{1}_{\{\tau^{(i)}<\tau^{(j)}\}}-\frac{1}{2}e^{-r\tau^{(i)}}(k\zeta_{\tau^{(i)}}^{(i)}+c_{i})\mathbf{1}_{\{\tau^{(i)}=\tau^{(j)}\}}\vert\mathscr{F}_{t}\biggr]\;.\label{eq:single-game-payoff}
\end{align}
The first term is the cumulative discounted profit, and the second
and third terms are the costs of investment. Note that $X_{t}^{\nu}$
coincides with the uncontrolled process $X_{t}$ for $t<\tau^{\nu}=\tau^{(i)}\wedge\tau^{(j)}$.
Note also that, although $X^{\nu}$ undergoes a boost at time $\tau^{\nu}$
so that $X_{t}^{\nu}\neq X_{t}$ for $t\ge\tau^{\nu}$, the controlled
state variable $X^{\nu}$ is not boosted anymore after time $\tau^{\nu}$
because of the single-investment assumption. This implies that the
single investment game is fundamentally a one-shot stopping game,
which will be established in the following subsection.

\subsubsection{Transformation into a Stopping Game in a Subgame Perfect Equilibrium
\label{subsubsec:Transformation-into-a}}

We now introduce the equilibrium solution concept that will be used
throughout this paper.

\textbf{Definition}. A strategy profile $\nu$ is a \emph{subgame
perfect equilibrium} (SPE) if it is a Nash equilibrium in any subgame.
Formally, it means that $V_{i,\tau}^{(\nu_{i},\nu_{j})}\ge V_{i,\tau}^{(\nu_{i}',\nu_{j})}$
at any stopping time $\tau$ and any admissible strategy $\nu_{i}'$.

\medskip{}

The definition of SPE stipulates that each player should optimize
their payoff at each point in time. Thus, if player $i$ invests at
time $\tau^{(i)}$ in an SPE $\nu$, the magnitude of the impulse
$\zeta_{\tau^{(i)}}^{(i)}$ should be the one that maximizes player
$i$'s payoff at time $\tau^{(i)}$. 

To facilitate the analysis, we make the following assumption for $\pi(\cdot)$
\citep{Alvarez2008}:

\begin{assumption}\label{assum:integrability} (i) $\lim_{x\downarrow a}\pi(x)=\pi_{L}>-\infty$
and (ii) $\mathbb{E}^{x}[\int_{0}^{\infty}\vert\pi(X_{t})\vert e^{-rt}dt]<\infty$
for any $x\in\mathscr{I}$.

\end{assumption} Here $\mathbb{E}^{x}[\cdot]=\mathbb{E}[\cdot\vert X_{0}=x]$
represents the conditional expected value given the initial condition
$X_{0}=x$. The implication of the integrability condition $\mathbb{E}^{x}[\int_{0}^{\infty}\vert\pi(X_{t})\vert e^{-rt}dt]<\infty$
is that the function 
\[
(R_{r}\pi)(x):=\mathbb{E}^{x}[\int_{0}^{\infty}\pi(X_{t})e^{-rt}dt]
\]
 is well-defined \citep{Alvarez2001}. In addition, we make the following
assumption to ensure a unique optimal value of $\zeta_{t}^{(i)}$:

\begin{assumption} \label{assum:single-inv} There is a unique $z^{*}\in\mathscr{I}$
such that $(R_{r}\pi)'(x)-k>0$ if $x<z^{*}$ and $(R_{r}\pi)'(x)-k<0$
if $x>z^{*}$. \end{assumption} Here $z^{*}$ has the meaning of
the optimal end-point of the boost in case the initial value of $x$
is less than $z^{*}$. If player $i$ invests at time $t$ with the
current value of the state variable $X_{t}$, the player's reward
from investment is $(R_{r}\pi)(z)-k(z-X_{t})-c_{i}$ if the player
boosts the state variable up to some value $z$. By Assumption \ref{assum:single-inv},
there is a well-defined unique value of $z$ that maximizes the reward
from investment.

Under the assumptions made thus far, we arrive at the following lemma:

\begin{lemma} \label{lemma:optimal-zeta} Under an SPE, the equilibrium
magnitude of the investment in the single investment game is 
\begin{equation}
\zeta_{t}^{(i)}=\max\{z^{*}-X_{t},0\}\:.\label{eq:SPE-zeta}
\end{equation}
 \end{lemma}

Intuitively, under an SPE, each player's optimal boost must maximize
the reward from investment. Because $\zeta_{t}^{(i)}$ is uniquely
determined by Lemma \ref{lemma:optimal-zeta}, the single investment
game can be transformed into a stopping game as established by Proposition
\ref{prop:stopping-game} below, where the \emph{reward from stopping}
is given by
\begin{equation}
g_{i}(x):=\left\{ \begin{array}{cc}
(R_{r}\pi)(z^{*})-k(z^{*}-x)-c_{i} & \text{if}\:x<z^{*}\\
(R_{r}\pi)(x)-c_{i} & \text{if}\:x\ge z^{*}
\end{array}\right.\;,\label{eq:gi}
\end{equation}
and the reward from the opponent's stopping is given by
\begin{equation}
m_{i}(x):=(R_{r}\pi)(\max\{z^{*},x\})\;.\label{eq:mi}
\end{equation}

\begin{proposition} \label{prop:stopping-game} Under an SPE $\nu$,
the payoff function can be transformed into one of a stopping game:
\begin{align}
V_{i,t}^{\nu}= & e^{rt}\mathbb{E}^{x}\biggl[\int_{t}^{\tau^{(i)}\wedge\tau^{(j)}}\pi(X_{s})e^{-rs}ds+g_{i}(X_{\tau^{(i)}})e^{-r\tau^{(i)}}\mathbf{1}_{\{\tau^{(i)}<\tau^{(j)}\}}+m_{i}(X_{\tau^{(j)}})e^{-r\tau^{(j)}}\mathbf{1}_{\{\tau^{(j)}<\tau^{(i)}\}}\label{eq:V-nu-stopping-time}\\
 & +\frac{1}{2}[g_{i}(X_{\tau^{(i)}})+m_{i}(X_{\tau^{(j)}})]e^{-r\tau^{(i)}}\mathbf{1}_{\{\tau^{(i)}=\tau^{(j)}\}}\vert\mathscr{F}_{t}\biggl]\;.\nonumber 
\end{align}
 \end{proposition}

By virtue of Proposition \ref{prop:stopping-game}, the single investment
game completely transforms into a stopping game with the \emph{uncontrolled}
state variable $X$. For the remainder of this section, we analyze
the game as if it is a stopping game. We therefore succinctly write
$\nu_{i}=\tau^{(i)}$ whenever necessary throughout this section.

\emph{Remark}: For simplicity of the presentation, the probability
space $\Omega$ and the measure $\mathbb{P}$ do not explicitly include
the probabilistic mixture of stopping times, which is necessary to
include mixed strategies in the strategy space. However, it is straightforward
to verify that the expression for the payoff function (\ref{eq:V-nu-stopping-time})
continues to be valid even if the probability space and the measure
are replaced by $\Omega\times L^{(1)}\times L^{(2)}$ and $\hat{\mathbb{P}}$
that incorporate the randomizers of mixed strategies; see Sections
\ref{subsubsec:Formulation_Mixed-strategy} and \ref{subsubsec:Characters_Mixed-strategy-MPE}
for the details.

\subsection{Benchmark: Single-Player Problem\label{subsec:Benchmark:-Single-Player}}

In Section \ref{subsec:Mixed-Strategy-SPE}, it will be shown that
the payoff associated with a mixed strategy MPE is identical to the
optimal payoff function from the single-player problem. Thus, it is
imperative to study the single-player optimal stopping problem first.
Towards this goal, we assume that player $j$ never invests in the
common good ($\tau^{(j)}=\infty$) and solve for player $i$'s best
response and payoff. By virtue of Proposition \ref{prop:stopping-game},
the single-player problem effectively reduces to the optimal stopping
problem of maximizing 
\begin{equation}
\mathbb{E}^{x}[\int_{0}^{\tau}\pi(X_{t})e^{-rt}dt+e^{-r\tau}g_{i}(X_{\tau})]\:,\label{eq:single-p}
\end{equation}
with respect to stopping time $\tau$.

To ensure the existence and uniqueness of the optimal stopping policy,
we assume the following:

\begin{assumption} \label{assump:theta} (i) There exists $\theta_{i}<z^{*}$
such that 
\begin{equation}
\alpha_{i}(x):=[g_{i}(x)-(R_{r}\pi)(x)]/\phi(x)\label{eq:alpha-i}
\end{equation}
 achieves a unique global maximum at $\theta_{i}$. Furthermore, $\alpha_{i}(x)$
increases for $x<\theta_{i}$ and decreases for $x>\theta_{i}$.

(ii) There exists $x_{i}^{c}\in(a,z^{*}]$ such that $\mathscr{A}g_{i}(x)+\pi(x)<0$
if and only if $x<x_{i}^{c}$. \end{assumption}

Assumptions \ref{assump:theta}(i) and (ii) are conventionally made
in the optimal stopping problem literature (see, for example, Theorem
3 of \citealp{Alvarez2001}). The function $\alpha_{i}(\theta_{i})$
has the meaning of the coefficient of $\phi(x)$ in the single player's
payoff function when the threshold of investment is $\theta_{i}$.
Thus, Assumption \ref{assump:theta}(i) ensures that there exists
a unique optimal threshold of investment. We make Assumption \ref{assump:theta}(ii)
to ensure that $\pi(x)$ takes sufficiently low values for $x<x_{i}^{c}$
so that it is always optimal for the players to invest if the state
variable $X$ is sufficiently low. Under the assumptions made so
far, we can obtain the following optimal solution:

\begin{proposition} \label{prop:single-p}(i) The optimal time of
investment for player $i$ is $\tau^{*}=\inf\{t\ge0:X_{t}\le\theta_{i}\}$,
and the optimal boost is $\zeta^{*}=z^{*}-X_{\tau^{*}}$. 

(ii) Player $i$'s optimal payoff given the current value $x$ of
the state variable is 
\begin{equation}
V_{i}^{*}(x)=\begin{cases}
\alpha_{i}(\theta_{i})\phi(x)+(R_{r}\pi)(x) & \text{if}\:x>\theta_{i}\\
g_{i}(x) & \text{otherwise}
\end{cases}\:,\label{eq:V-single}
\end{equation}
where $\alpha_{i}(\cdot)$ and $g_{i}(\cdot)$ are defined by (\ref{eq:alpha-i})
and (\ref{eq:gi}).

(iii) $\theta_{i}$ decreases in $c_{i}$. \end{proposition}

The threshold stopping rule of Proposition \ref{prop:single-p}(i)
is an intuitive result: an investment is made if the value of $X$
is sufficiently low. Proposition \ref{prop:single-p}(ii) follows
from the well-established optimal stopping theory \citep{Alvarez2001}.
Lastly, Proposition \ref{prop:single-p}(iii) holds because a player
with a higher cost of investment $c_{i}$ would have lower incentive
to invest, leading to a lower threshold of investment.

\subsection{Mixed Strategy SPE of Symmetric Games \label{subsec:Mixed-Strategy-SPE}}

Next, we return to the two-player game and obtain mixed strategy SPEs.
We assume that the players are identical, i.e., $c:=c_{1}=c_{2}$,
in which case we have $\theta:=\theta_{1}=\theta_{2}$, and construct
a mixed strategy MPE and asymmetric two-stage SPE (TSSPE).

\subsubsection{Formulation of Mixed Strategies \label{subsubsec:Formulation_Mixed-strategy}}

We first expand the strategy space to encompass probabilistic mixtures
of stopping times. The canonical definition of a mixed strategy requires
an additional random variable for each player which is realized and
revealed to each player in the beginning of the game but unknown to
the opponent. To incorporate the additional random variable into a
player's strategy, we follow the recipe provided by Section 7.1 of
\citet{Touzi2002}. We first re-define $(\Omega,\mathbb{F},\mathscr{F},\mathbb{P})$
as the probability space of the sample paths of $X$ (a Wiener space).
By this definition, $\mathbb{F}$ is the canonical filtration of the
Brownian motion $W$ introduced in (\ref{eq:SDE}), and $\mathbb{P}$
is the Wiener measure. Following \citet{Touzi2002}, we then augment
$(\Omega,\mathbb{F},\mathscr{F},\mathbb{P})$ with the probability
space of player $i$'s randomizer defined on $L^{(i)}=[0,1]$ with
a Borel $\sigma$-algebra $\mathscr{L}^{(i)}$ on $L^{(i)}$ and the
Lebesgue measure $\mathbb{L}^{(i)}$ on $L^{(i)}$. By this recipe,
the expanded probability space is $(\Omega\times L^{(1)}\times L^{(2)},(\mathscr{\hat{F}}_{t})_{t\ge0},\hat{\mathscr{F}},\hat{\mathbb{P}})$,
where the new probability measure is defined as $\hat{\mathbb{P}}:=\mathbb{P}\otimes\mathbb{L}^{(1)}\otimes\mathbb{L}^{(2)}$,
and the new $\sigma$-algebra is given by $\mathscr{\hat{F}}=\mathscr{F}\otimes\mathscr{L}^{(1)}\otimes\mathscr{L}^{(2)}$.
Similarly, the new filtration is defined as the product $\sigma$-algebra
$\mathscr{\hat{F}}_{t}=\mathscr{F}_{t}\otimes\mathscr{L}^{(1)}\otimes\mathscr{L}^{(2)}$.

Next, we formulate the mixed strategy $\nu_{i}$. Following \citet{Touzi2002},
we can construct a mixed strategy by mapping $L^{(i)}=[0,1]$ to the
space of stopping times adapted to $\mathbb{F}=(\mathscr{\mathscr{F}}_{t})_{t\ge0}$.
To construct such a mapping, we introduce a strategy-specific \emph{survival
probability} $M^{(i)}=\{M_{t}^{(i)}:t\ge0,M_{0}^{(i)}=1\}$, which
is a non-negative, non-increasing, and right-continuous process adapted
to $\mathbb{F}$ .  The survival probability $M_{t}^{(i)}$ has
the interpretation of the probability that player $i$ has not yet
invested by time $t$, and hence, player $i$ invests at a rate of
$-dM_{t}^{(i)}/M_{t}^{(i)}$ at any given time $t$. Then we parameterize
player $i$'s investment strategy $\nu_{i}$ as $\hat{\tau}^{(i)}=\inf\{t\ge0:M_{t}^{(i)}\le\hat{l}^{(i)}\}$,
where $\hat{l}^{(i)}$ is player $i$'s randomizer that has a uniform
distribution over $L^{(i)}=[0,1]$.

Note that for any specific realized value of $l^{(i)}\in L^{(i)}$,
the random time $\inf\{t\ge0:M_{t}^{(i)}\le l^{(i)}\}$ is a stopping
time adapted to $\mathbb{F}$ because $M^{(i)}$ is adapted to $\mathbb{F}$.
Hence, because $\hat{l}^{(i)}$ is a random variable of measure $\mathbb{L}^{(i)}$,
the stopping time $\hat{\tau}^{(i)}=\inf\{t\ge0:M_{t}^{(i)}\le\hat{l}^{(i)}\}$
is a stopping time adapted to $(\mathscr{\hat{F}}_{t}^{(i)})_{t\ge0}:=(\mathscr{F}_{t}\otimes\mathscr{L}^{(i)})_{t\ge0}$.
Note also that the probability that player $i$ has not invested by
time $t$ is given by $\mathbb{L}^{(i)}(\{\hat{\tau}^{(i)}>t\})=M_{t}^{(i)}$,
which leads to the required rate of investment $-dM_{t}^{(i)}/M_{t}^{(i)}$.
Finally, we remark that this construction is equivalent to directly
mixing stopping times through a functional analysis method (Section
7 of \citealt{Touzi2002}). 

\subsubsection{Characteristics of Mixed Strategy MPE \label{subsubsec:Characters_Mixed-strategy-MPE}}

We are now in a position to provide the definition of MPE. Below we
let $\mathbb{E}_{\mathbb{L}^{(i)}}$ denote the expectation with respect
to the measure $\mathbb{L}^{(i)}$.

\textbf{Definition}. An SPE $\nu=(\nu_{1},\nu_{2})$ is a \emph{Markov
perfect equilibrium} (MPE) if each player's strategy $\nu_{i}$ of
investment at time $\hat{\tau}^{(i)}$ satisfies the following condition:
\begin{equation}
\mathbb{E}_{\mathbb{L}^{(i)}}[\mathbf{1}_{\{\hat{\tau}^{(i)}\in(\tau,\tau+s]\}}\vert\mathscr{\hat{F}}_{\tau},\hat{\tau}^{(i)}>\tau]=\mathbb{E}_{\mathbb{L}^{(i)}}[\mathbf{1}_{\{\hat{\tau}^{(i)}\in(\tau,\tau+s]\}}\vert X_{\tau},\hat{\tau}^{(i)}>\tau]\:,\label{eq:strong-Markov-strategy}
\end{equation}
 for any stopping time $\tau$ and $s>0$. If $\nu$ is an MPE, we
call $\nu_{1}$ and $\nu_{2}$ \emph{Markov strategies}.

\medskip{}

The Markov property (\ref{eq:strong-Markov-strategy}) implies that
the probability distribution of the investment time $\hat{\tau}^{(i)}$
depends \emph{only on the current value} of the state variable rather
than its full history up to time $\hat{\tau}^{(i)}$ \citep{Maskin2001}.

For the time being, we are solely concerned with MPE, so we focus
on characterizing Markov strategies. A Markov strategy of stopping
is mathematically equivalent to a Markovian \emph{killing} of a Markov
process, which can be equivalently formulated in terms of a \emph{multiplicative
functional} and an \emph{additive functional.} (For the definitions
of multiplicative functionals and additive functionals, see, for example,
Chapters III and IV, \citealt{Blumenthal2007}). Thus, following Chapter
IV of \citet{Blumenthal2007}, we can conveniently express a \emph{Markovian}
survival probability as a\emph{ multiplicative functional} given by
\begin{equation}
M_{t}^{(i)}=\exp(-A_{t}^{(i)})N_{t}^{(i)}\:,\label{eq:Multi_Functional}
\end{equation}
where $A_{t}^{(i)}$ is an \emph{additive functional} of the process
$(X_{s})_{s\in[0,t]}$,  and $N_{t}^{(i)}$ is an auxiliary multiplicative
functional of the form $N_{t}^{(i)}=\mathbf{1}_{\{t<\tau_{E}^{(i)}\}}$.
By the definition of an additive functional, $A_{t}^{(i)}$ is a non-negative,
non-decreasing, and right-continuous process adapted to $\mathbb{F}$
with the initial value $A_{0}^{(i)}=0$. Here $\tau_{E}^{(i)}$ is
a stopping time at which player $i$ invests with probability one.
According to this survival probability, player $i$ will invest at
the probabilistic rate of $dA_{t}^{(i)}$ for $t<\tau_{E}^{(i)}$,
but he will invest at time $\tau_{E}^{(i)}$ with probability one.

By the stipulation that the survival probability (\ref{eq:Multi_Functional})
must be Markov, the time of definite investment $\tau_{E}^{(i)}$
must be a \emph{hitting time} $\tau_{E}^{(i)}=\inf\{t\ge0:X_{t}\in E_{i}\}$
of some closed set $E_{i}$ that has an interpretation of \emph{pure
strategy} investment region; by the diffusive property of $X$, $\tau_{E}^{(i)}$
is equal to the hitting time of the closure of $E_{i}$ even if $E_{i}$
is not closed (3.22 on p. 312, \citealt{Revuz1999}). We also remark
that the definition of additive functionals (II.21, \citealt{Borodin1996})
 dictates that the probability distribution of $A_{s+\tau}^{(i)}-A_{\tau}^{(i)}$
for some $s>0$ conditional on $\mathscr{F}_{\tau}$ depends solely
on the current value $X_{\tau}$ and the time $s$, which naturally
renders the survival probability (\ref{eq:Multi_Functional}) Markovian.
Thus, a Markov strategy can be completely and unambiguously represented
by the set $E_{i}$ and the additive functional $A^{(i)}$.

There are two well-established properties of an additive functional
$A^{(i)}$ of a \emph{regular diffusion} process $X$: (a) $A_{t}^{(i)}$
is continuous almost surely (II.21, \citealt{Borodin1996}), and (b)
$A^{(i)}$ can be expressed as $A_{t}^{(i)}=\int_{0}^{t}\lambda_{i}(X_{s})ds$,
a time-integral of a non-negative Radon-Nikodym derivative $\lambda_{i}(X_{t})$
 (II.23 and II.24, \citealt{Borodin1996}). In addition to these
convenient properties, we impose two regularity conditions to rule
out unrealistic and impractical strategies as well as to ensure analytical
tractability. First, we stipulate that the Radon-Nikodym derivative
$\lambda_{i}(X_{t})$ of $A_{t}^{(i)}$ is without singularities so
that $A_{t}^{(i)}$ is differentiable with respect to time $t$.\footnote{ A singularity in $\lambda_{i}(X_{t})$ implies a singularly high
rate of investment at time $t$. For instance, a \emph{local time}
process does not possess a finite Radon-Nikodym derivative, but it
can be approximated as one with extremely high values of finite Radon-Nikodym
derivatives within a very small region; EC Appendix \ref{App_sec:Singular-Rates-of}
illustrates this point. However, such an excessively large values
of $\lambda_{i}(\cdot)$ is not possible in a mixed strategy MPE because
the equilibrium rate $\lambda_{i}(\cdot)$ must possess a specific
functional form given by (\ref{eq:lambda}). Thus, we exclude the
possibility of non-differentiable $A_{t}^{(i)}$.}  Secondly, we stipulate that the mixed strategy region $\Gamma_{i}:=\{x\in\mathscr{I}:\lambda_{i}(x)>0\}$
is a \emph{union of disjoint intervals}. This condition rules out
any pathological topology of mixed strategy regions such as in a Cantor
set, which is unrealistic and impractical to implement for real-life
decision makers. Under this stipulation, we can always represent $\Gamma_{i}$
as a union of disjoint open intervals because adding a countable set
of points to $\Gamma_{i}$ (or subtracting one from $\Gamma_{i}$)
does not alter the outcome of the game almost surely. Hence, we will
subsequently assume that $\Gamma_{i}$ is an open set.

The next proposition provides a convenient alternative expression
for the equilibrium payoff function. To this end, we first re-express
(\ref{eq:V-nu-stopping-time}) in terms of the expanded probability
space and the strategy profile $\nu=(\hat{\tau}^{(i)},\hat{\tau}^{(j)})$,
where $t=0$ and $X_{0}=x$, as follows:
\begin{align}
V_{i,0}^{\nu}= & \hat{\mathbb{E}}_{j}^{x}\biggr[\int_{0}^{\hat{\tau}^{(i)}\wedge\hat{\tau}^{(j)}}\pi(X_{t})e^{-rt}dt+g_{i}(X_{\hat{\tau}^{(i)}})e^{-r\hat{\tau}^{(i)}}\mathbf{1}_{\{\hat{\tau}^{(i)}<\hat{\tau}^{(j)}\}}+m_{i}(X_{\hat{\tau}^{(j)}})e^{-r\hat{\tau}^{(j)}}\mathbf{1}_{\{\hat{\tau}^{(j)}<\hat{\tau}^{(i)}\}}\label{eq:V-nu-new-probability-space}\\
 & +\frac{1}{2}[g_{i}(X_{\hat{\tau}^{(i)}})+m_{i}(X_{\hat{\tau}^{(j)}})]e^{-r\hat{\tau}^{(i)}}\mathbf{1}_{\{\hat{\tau}^{(i)}=\hat{\tau}^{(j)}\}}\biggl]\;.\nonumber 
\end{align}
Note that this expression is simply an extension of (\ref{eq:V-nu-stopping-time})
except that $\tau^{(i)}$ and $\tau^{(j)}$ are replaced by $\hat{\tau}^{(i)}$
and $\hat{\tau}^{(j)}$, where $\hat{\tau}^{(i)}$ is a $\mathbb{F}$-stopping
time while $\hat{\tau}^{(j)}=\inf\{t\ge0:M_{t}^{(j)}\le\hat{l}^{(j)}\}$
is a $(\mathscr{\hat{F}}_{t}^{(j)})_{t\ge0}$-stopping time. We consider
$\hat{\tau}^{(i)}$ a $\mathbb{F}$-stopping time because $\hat{\tau}^{(i)}$
is a pure strategy stopping time from player $i$'s perspective while
player $j$ employs a mixed strategy from player $i$'s perspective.
As we will show in the remainder of Section \ref{sec:One-Time} below,
it is indeed enough to consider the payoff from a $\mathbb{F}$-stopping
time for player $i$. This is because if a probabilistic mixture of
stopping times is a best response to a given rival's strategy, then
each pure strategy stopping time in the mixture must also be a best
response to the rival's strategy by the definition of the mixed strategy
equilibrium (p.665, \citealt{Hendricks1988}; p.12, \citealt{Steg2015}).
Another difference from (\ref{eq:V-nu-stopping-time}) is that $\mathbb{E}^{x}$
is replaced by $\hat{\mathbb{E}}_{j}^{x}[\cdot]:=\mathbb{E}_{\mathbb{P}\otimes\mathbb{L}^{(j)}}[\cdot\vert X_{0}=x]$
representing the expectation over $\mathbb{P}\otimes\mathbb{L}^{(j)}$
conditional on $X_{0}=x$. 

The following proposition establishes an alternative expression for
(\ref{eq:V-nu-new-probability-space}).

\begin{proposition} \label{prop:alternative-V} Given player $j$'s
Markov strategy $\nu_{j}$ determined by $\lambda_{j}(\cdot)$ and
$E_{j}$,  player $i$'s payoff from investing at a $\mathbb{F}$-stopping
time $\hat{\tau}^{(i)}$ is given by 
\begin{align}
V_{i}^{\nu}(x)= & \mathbb{E}^{x}\biggr[\int_{0}^{\hat{\tau}^{(i)}\wedge\tau_{E}^{(j)}}[\pi(X_{t})+\lambda_{j}(X_{t})m_{i}(X_{t})]e^{-rt-A_{t}^{(j)}}dt+\mathbf{1}_{\{\hat{\tau}^{(i)}<\tau_{E}^{(j)}\}}g_{i}(X_{\hat{\tau}^{(i)}})e^{-r\hat{\tau}^{(i)}-A_{\hat{\tau}^{(i)}}^{(j)}}\label{eq:V-nu-alternative}\\
 & +\mathbf{1}_{\{\hat{\tau}^{(i)}>\tau_{E}^{(j)}\}}m_{i}(X_{\tau_{E}^{(j)}})e^{-r\tau_{E}^{(j)}-A_{\tau_{E}^{(j)}}^{(j)}}+\frac{1}{2}\mathbf{1}_{\{\hat{\tau}^{(i)}=\tau_{E}^{(j)}\}}[g_{i}(X_{\hat{\tau}^{(i)}})+m_{i}(X_{\hat{\tau}^{(i)}})]e^{-r\hat{\tau}^{(i)}-A_{\hat{\tau}^{(i)}}^{(j)}}\biggl]\:.\nonumber 
\end{align}
 \end{proposition} By virtue of Proposition \ref{prop:alternative-V},
the payoff (\ref{eq:V-nu-new-probability-space}) is conveniently
re-expressed as an expectation operator $\mathbb{E}[\cdot]$ over
a new expression for the payoff that incorporates player $j$'s strategy
through $\lambda_{j}(\cdot)$ and $\tau_{E}^{(j)}$.

\subsubsection{Mixed strategy MPE \label{subsubsec:Mixed-MPE}}

We now construct an MPE strategy profile $\nu^{M}=(\nu_{1}^{M},\nu_{2}^{M})$.
As per Section \ref{subsubsec:Characters_Mixed-strategy-MPE}, we
only need to specify $\lambda_{i}(\cdot)$ and $E_{i}$ to determine
$\nu_{i}^{M}$. Specifically, we choose $E_{i}=\emptyset$ for both
$i\in\{1,2\}$ and stipulate the form of $\lambda_{i}(\cdot)$ with
respect to $g(\cdot):=g_{1}(\cdot)=g_{2}(\cdot)$ and $\pi(\cdot)$
as follows:
\begin{equation}
\lambda_{i}(x):=\mathbf{1}_{\{x\in(a,\theta)\}}\frac{-\mathscr{A}g(x)-\pi(x)}{(R_{r}\pi)(z^{*})-g(x)}\:.\label{eq:lambda}
\end{equation}
Recall that $\mathscr{A}$ is the characteristic differential operator
of $X$ defined by (\ref{eq:char-op}). Note that the rate of investment
is non-zero only if $X_{t}\in(a,\theta)$, where $\theta:=\theta_{1}=\theta_{2}$.
Hence, we call $(a,\theta)$ the \emph{common} \emph{mixed strategy
investment region}. The function $\lambda_{i}(x)$ is non-negative
because both the numerator and the denominator are positive for $x<\theta$.
First, $\mathscr{A}g(x)+\pi(x)<0$ holds for $x<\theta$ because of
the inequality $\theta\le x_{i}^{c}$ by the theory of optimal stopping
(Chapter 10 of \citealp{Oksendal2003}) where $x_{i}^{c}$ is given
in Assumption \ref{assump:theta}(ii). Secondly, $(R_{r}\pi)(z^{*})>g(x)$
for all $x\in(a,\theta)$ by the definition of $g(\cdot)$ given in
(\ref{eq:gi}) and $\theta_{i}<z^{*}$ in Assumption \ref{assump:theta}(i).

Below we let $V_{i}^{\nu^{M}}(x)$ denote the payoff function associated
with $\nu^{M}$ when the current state variable is $x$.

\begin{proposition} \label{prop:MixedStrategyMPE} The strategy profile
$\nu^{M}$ defined above is an MPE. Furthermore, $V_{i}^{\nu^{M}}(x)=V_{i}^{*}(x)$
given by (\ref{eq:V-single}). \end{proposition}

Proposition \ref{prop:MixedStrategyMPE} establishes that the strategy
profile $\nu^{M}$ is a mixed strategy MPE. It also establishes that
the payoff under $\nu^{M}$ is identical to the payoff from the single-player
optimal stopping solution. The equality $V_{i}^{\nu^{M}}(x)=V_{i}^{*}(x)$
has an heuristic explanation. According to $\nu^{M}$, none of the
players invests until the hitting time of $(a,\theta)$, where $\theta$
is identical to the threshold of investment in the single-player solution,
so $(\theta,b)$ takes the role of a continuation region for both
players. Furthermore, the equilibrium payoff to player $i$ within
$(a,\theta)$ is equal to the payoff from immediate investment, i.e.,
$g_{i}(x)$; this is because whenever $X$ is in $(a,\theta)$, which
is the common mixed strategy region under $\nu^{M}$, an immediate
investment must be one of player $i$'s best responses as per the
definition of a mixed strategy equilibrium \citep{Hendricks1988,Georgiadis2019}.
Thus, the payoff function has to satisfy the boundary condition $V_{i}^{\nu^{M}}(\theta)=g_{i}(\theta)$.
Since $(\theta,b)$ is a continuation region for player $i$, $V_{i}^{\nu^{M}}(x)$
can be obtained as the solution to a boundary value problem in the
interval $[\theta,b)$ with the boundary condition $V_{i}^{\nu^{M}}(\theta)=g_{i}(\theta)$.
Note that $V_{i}^{*}(x)$ is also the solution to a boundary value
problem in the interval $[\theta,b)$ with the boundary condition
$V_{i}^{*}(\theta)=g_{i}(\theta)$. Hence, we have $V_{i}^{\nu^{M}}(x)=V_{i}^{*}(x)$.

\subsubsection{Two-stage SPE \label{subsubsec:2-stage-SPE}}

Next, we construct a TSSPE $\nu^{S}$. As the name suggests, the equilibrium
path of $\nu^{S}$ consists of two stages. In the first stage, each
player $i$'s strategy is to invest with probability $q_{i}$ upon
reaching $\tau_{\Gamma}:=\inf\{t\ge0:X_{t}\in\Gamma\}$, the hitting
time of the interval $\Gamma:=(a,\theta)$, which we identified as
the common mixed strategy region for the MPE in Section \ref{subsubsec:Mixed-MPE}.
In the second stage, the players employ the mixed strategy MPE profile
presented in Section \ref{subsubsec:Mixed-MPE}. We can thus define
the first stage as the time interval $[0,\tau_{\Gamma}]$ and the
second stage as $(\tau_{\Gamma},\infty)$. We also stipulate that
$q_{1}q_{2}=0$; if $q_{1}$ and $q_{2}$ are both positive, then
there is positive probability of simultaneous investment at $\tau_{\Gamma}$,
and hence, one of the players can improve his payoff by forgoing with
investment at time $\tau_{\Gamma}$. Without loss of generality, we
set $q_{1}=0$ and $q_{2}>0$. Note that player 2's strategy is non-Markov,
and hence, $\nu^{S}$ is \emph{not} an MPE; player 2's plan of action
depends not only on the current value of $X$ but also on whether
$X$ has hit $\Gamma$ in the past history.

Based on the strategy profile described above, we can construct survival
probability $M_{t}^{(i)}$ as follows:
\begin{align}
M_{t}^{(1)} & =\exp[-\mathbf{1}_{\{t>\tau_{\Gamma}\}}\int_{\tau_{\Gamma}}^{t}\lambda_{1}(X_{t})dt]\:,\label{eq:M1}\\
M_{t}^{(2)} & =\begin{cases}
1 & \text{if}\:t<\tau_{\Gamma}\\
(1-q_{2})\exp[-\int_{\tau_{\Gamma}}^{t}\lambda_{2}(X_{t})dt] & \text{if}\:t\ge\tau_{\Gamma}
\end{cases}\;.\label{eq:M2}
\end{align}
The expression (\ref{eq:M1}) reflects player 1's strategy of investment
at the hazard rate of $\lambda_{1}(X_{t})$ in the second stage when
$t>\tau_{\Gamma}$. The expression (\ref{eq:M2}) is also based on
player 2's strategy. Player 2 does nothing until $\tau_{\Gamma}$
during stage 1, so $M_{t}^{(2)}$ is invariant in time for $t<\tau_{\Gamma}$.
At time $\tau_{\Gamma}$, he invests with a probability of $q_{2}$,
which results in the discontinuity of $M_{t}^{(2)}$. In the second
stage, player 2 invests at the hazard rate of $\lambda_{2}(X_{t})$,
and hence, the exponential factor $\exp[-\int_{\tau_{\Gamma}}^{t}\lambda_{2}(X_{t})dt]$.

We let $V_{i}^{\nu^{S}}(x,s)$ denote the payoff to player $i$ associated
with $\nu^{S}$ in stage $s\in\{1,2\}$ when the current state variable
is $x$. Since stage 2 reduces to the mixed strategy MPE, we obtain
$V_{i}^{\nu^{S}}(x,2)=V_{i}^{\nu^{M}}(x)$ where $\nu^{M}$ is the
mixed strategy MPE shown in Section \ref{subsubsec:Mixed-MPE}. The
first-stage payoff functions are given as follows:
\begin{align}
V_{1}^{\nu^{S}}(x,1) & =\mathbb{E}^{x}\left\{ \int_{0}^{\tau_{\Gamma}}\pi(X_{t})e^{-rt}dt+e^{-r\tau_{\Gamma}}[(1-q_{2})V_{1}^{\nu^{M}}(X_{\tau_{\Gamma}})+q_{2}(R_{r}\pi)(z^{*})]\right\} \:,\label{eq:V11}\\
V_{2}^{\nu^{S}}(x,1) & =\mathbb{E}^{x}\left\{ \int_{0}^{\tau_{\Gamma}}\pi(X_{t})e^{-rt}dt+e^{-r\tau_{\Gamma}}[(1-q_{2})V_{2}^{\nu^{M}}(X_{\tau_{\Gamma}})+q_{2}g(X_{\tau_{\Gamma}})]\right\} \label{eq:V21}\\
 & =\mathbb{E}^{x}\left\{ \int_{0}^{\tau_{\Gamma}}\pi(X_{t})e^{-rt}dt+e^{-r\tau_{\Gamma}}g(X_{\tau_{\Gamma}})\right\} \:.\label{eq:V21A}
\end{align}
 The expression for $V_{1}^{\nu^{S}}(x,1)$ in (\ref{eq:V11}) indicates
that player 2 invests at $\tau_{\Gamma}$ with a probability $q_{2}>0$,
resulting in the immediate payoff $(R_{r}\pi)(z^{*})$ to player 1;
with a probability $1-q_{2}$, the game enters stage 2, and player
1 earns the payoff $V_{1}^{\nu^{M}}(X_{\tau_{\Gamma}})$ at time $\tau_{\Gamma}$.
Similarly, player 2's payoff at $\tau_{\Gamma}$ in (\ref{eq:V21})
is $g(X_{\tau_{\Gamma}})$ if he invests at $\tau_{\Gamma}$ and $V_{2}^{\nu^{M}}(X_{\tau_{\Gamma}})$
if he does not. The payoff $V_{2}^{\nu^{S}}(x,1)$ reduces to (\ref{eq:V21A})
because $V_{2}^{\nu^{M}}(X_{\tau_{\Gamma}})=g(X_{\tau_{\Gamma}})$
by (\ref{eq:V-single}) and Proposition \ref{prop:MixedStrategyMPE}.
The following proposition establishes that $\nu^{S}$ is an equilibrium.

\begin{proposition} \label{prop:SPE-mixed} The strategy profile
$\nu^{S}$ defined above is an SPE. \end{proposition}

The mixed strategy TSSPE $\nu^{S}$ is the stochastic analog of the
mixed strategy SPE in Theorem 3 of \citet{Hendricks1988}. In fact,
all SPEs of deterministic wars of attrition take the form of the TSSPE
\citep{Hendricks1988}. This is one of the simplest classes of SPE,
and thus the most likely to be observed in practice. For this reason,
we search for TSSPE in asymmetric games in the remainder of this paper.

\subsection{Absence of Two-stage SPE in Asymmetric Games \label{subsec:Asymmetric-Game}}

We now establish the central result of this section: an asymmetric
game has no TSSPE that can be represented by survival probability
of the form given in (\ref{eq:M1}) and (\ref{eq:M2}). This result
is a slight extension of Theorem 1 of \citet{Georgiadis2019}; the
only difference between the model examined by \citet{Georgiadis2019}
and our current model is in the functional form of $g_{i}(\cdot)$
and $m_{i}(\cdot)$. Nevertheless, for the sake of completeness, we
provide the complete proof in EC Appendix \ref{App_sec:Mathematical-Proofs}.

\begin{theorem} \label{thm:MPE-non-existence} In a single-investment
game with $c_{1}\neq c_{2}$, there exists no MPE or TSSPE. \end{theorem}

Theorem \ref{thm:MPE-non-existence} holds because the two players
must have the common mixed strategy region. If the two players have
differing costs $c_{1}\neq c_{2}$, then the boundary points of their
mixed strategy region cannot coincide due to the asymmetry, and hence,
$\Gamma_{1}=\Gamma_{2}$ cannot be satisfied.

Note, however, that Theorem \ref{thm:MPE-non-existence} does not
claim to have proved that a mixed strategy equilibrium of any kind
is impossible under asymmetry. Instead, the main point of Theorem
\ref{thm:MPE-non-existence} is that an asymmetric single investment
game does not possess a mixed strategy equilibrium that \emph{shares
the same characteristics} with known mixed strategy equilibria obtained
in Sections \ref{subsubsec:Mixed-MPE} and \ref{subsubsec:2-stage-SPE}
as well as in \citet{Hendricks1988} and \citet{Steg2015}.\footnote{In all of the known mixed strategy equilibria, the mixed strategy
regions have simple topologies (unions of disjoint intervals), and
$A_{t}^{(i)}$'s possess a Radon-Nikodym derivative $\lambda_{t}^{(i)}$.} In contrast, the next section demonstrates that there \emph{is} such
an equilibrium in the infinite investment game despite asymmetry,
which is the central goal of this paper.

\section{Impulse Control Game \label{sec:Impulse}}

We now turn to a game with an infinite number of investment opportunities.
The goal of this section is to prove that a moderately asymmetric
game possesses a mixed strategy TSSPE. This result is in stark contrast
to the result from the single investment model (Theorem \ref{thm:MPE-non-existence})
or the stochastic concession games examined by \citet{Georgiadis2019}.
Since we focus on asymmetric games, we assume $c_{1}>c_{2}$ without
loss of generality. We first derive a verification theorem (Theorem
\ref{thm:nu*-equilibrium}) and use it to construct a mixed strategy
SPE (Theorem \ref{thm:Existence_SPE}). We then illustrate an example
of a mixed strategy TSSPE in Section \ref{subsec:Example}. Lastly,
we compare the mixed strategy TSSPE to a pure strategy SPE from the
perspective of a social planner in Section \ref{subsec:Comparison-to-Pure}.

\subsection{The Model \label{subsec:The-Model}}

As a first step, we extend the model introduced in Section \ref{subsec:The-Model-and}
to allow for an infinite number of investment opportunities. We assume
the same properties of the uncontrolled process $X$ as in Section
\ref{subsubsec:State-Variable}, but the payoff and the strategy space
differ.

We let $(\Omega,\mathscr{F},\mathbb{F},\mathbb{P})$ denote the probability
space of the sample paths of a Wiener process $W$ that satisfies
the usual condition. Then we let $\hat{\mathscr{S}}=(\hat{\Omega},\hat{\mathscr{F}},(\hat{\mathscr{F}}_{t})_{t\ge0},\hat{\mathbb{P}})$
denote the probability space $(\Omega,\mathscr{F},\mathbb{F},\mathbb{P})$
augmented by the randomizers of mixed strategies, which we will specify
later in this subsection, and we let $\nu_{i}=(\eta_{n}^{(i)},\xi_{n}^{(i)})_{n\in\mathbb{N}}$
denote player $i$'s strategy that stipulates a series of strictly
increasing $(\hat{\mathscr{F}}_{t})_{t\ge0}$-stopping times of investment
$(\eta_{n}^{(i)})_{n\in\mathbb{N}}$ and the boosts $(\xi_{n}^{(i)})_{n\in\mathbb{N}}$.
Each boost $\xi_{n}^{(i)}$ is measurable with respect to $\hat{\mathscr{F}}_{\eta_{n}^{(i)}}$.
To avoid the pathological scenario of infinitely many investments
within an infinitesimal time interval, we impose the condition $\lim_{n\rightarrow\infty}\eta_{n}^{(i)}=\infty$
\citep{Alvarez2008}.

For convenience, we partition the timeline into \emph{periods} of
investment and specify the player's strategy for each period. We let
$m$ denote the period index which keeps track of the total number
of investments made on the common good. We let $(T_{m})_{m\in\{0,1,2,...\}}$
denote the set of stopping times of investment made by either player
such that $T_{m+1}\ge T_{m}$. As a matter of convention, we define
$T_{0}=0$. By this convention, $T_{m}$ indicates the timing of the
$m$-th investment, and the $m$-th period is defined as $[T_{m},T_{m+1})$.
We also note that $\lim_{m\rightarrow\infty}T_{m}=\infty$ by the
stipulation that $\lim_{n\rightarrow\infty}\eta_{n}^{(i)}=\infty$
for each $i$. We focus on the strategy profiles which have no dependence
on $m$ because this is the simplest class of strategy profiles and
easiest to implement for decision makers in practice.

Next, we construct the probability space. Just as in Section \ref{subsubsec:Formulation_Mixed-strategy},
we utilize the recipe provided by Section 7.1 of \citet{Touzi2002}
for each period and construct an expanded probability space $\hat{\mathscr{S}}$.
Specifically, the new probability space is given by 
\[
\hat{\mathscr{S}}=(\Omega\times\prod_{m\in\mathbb{N}}L_{m}^{(i)}\times\prod_{m\in\mathbb{N}}L_{m}^{(j)},(\mathscr{\hat{F}}_{t})_{t\ge0},\hat{\mathscr{F}},\hat{\mathbb{P}})\:,
\]
where $L_{m}^{(i)}=[0,1]$, $\hat{\mathbb{P}}:=\mathbb{P}\otimes_{m\in\mathbb{N}}\mathbb{L}_{m}^{(i)}\otimes_{m\in\mathbb{N}}\mathbb{L}_{m}^{(j)}$,
where $\mathbb{L}_{m}^{(i)}$ is the Lebesgue measure on $L_{m}^{(i)}$.
We let $\mathscr{L}_{m}^{(i)}$ denote the Borel $\sigma$-algebra
on $L_{m}^{(i)}$, and we define $\hat{\mathscr{F}}=\mathscr{F}\otimes_{m\in\mathbb{N}}\mathscr{L}_{m}^{(i)}\otimes_{m\in\mathbb{N}}\mathscr{L}_{m}^{(j)}$.
Likewise, $\mathscr{\hat{F}}_{t}$ is a product $\sigma$-algebra
given by $\mathscr{F}_{t}\otimes_{m\in\mathbb{N}}\mathscr{L}_{m}^{(i)}\otimes_{m\in\mathbb{N}}\mathscr{L}_{m}^{(j)}$.

Following \citet{Touzi2002}, we let player $i$'s $m$-th period
mixed strategy be characterized by an $m$-th period survival probability
$M_{m}^{(i)}=(M_{m,t}^{(i)})_{t\in[T_{m},T_{m+1})}$ with the initial
condition $M_{m,T_{m}}^{(i)}=1$ and an $m$-th period randomizer
$\hat{l}_{m}^{(i)}\in[0,1]$ with a Lebesgue measure $\mathbb{L}_{m}^{(i)}$
on $[0,1]$. The investment stopping time is then defined as $\tau_{m+1}^{(i)}:=\inf\{t\ge T_{m}:M_{m,t}^{(i)}\le\hat{l}_{m}^{(i)}\}$,
and the period $m$ begins at the stopping time $T_{m}=\tau_{m}^{(i)}\wedge\tau_{m}^{(j)}$.
We also let $\zeta_{m+1}^{(i)}$ denote the $\hat{\mathscr{F}}_{\tau_{m+1}^{(i)}}$-measurable
boost at time $\tau_{m+1}^{(i)}$ executed by player $i$. Just as
in Section \ref{subsubsec:Formulation_Mixed-strategy}, utilizing
the expanded probability space $\hat{\mathscr{S}}$ for each period
is equivalent to directly mixing stopping times  by virtue of the
proof provided in Section 7 of \citet{Touzi2002}. Furthermore we
assume the same tie-breaking rule as in Section \ref{sec:One-Time}:
if $\tau_{m}^{(i)}=\tau_{m}^{(j)}$, then only one of the two players
goes through with the investment with a probability of 50\%. Finally,
under the strategy profile $\nu=(\nu_{1},\nu_{2})$, the controlled
state variable satisfies the following stochastic integral equation:
\[
X_{t}^{\nu}=X_{0}+\int_{0}^{t}\mu(X_{s}^{\nu})ds+\int_{0}^{t}\sigma(X_{s}^{\nu})dW_{s}+\sum_{m=1}^{\infty}\xi_{m}^{\nu}\mathbf{1}_{\{T_{m}\le t\}}\;,
\]
where $\xi_{m}^{\nu}=\xi_{m}^{(i)}$ if $\tau_{m}^{(i)}<\tau_{m}^{(j)}$,
and $\xi_{m}^{\nu}$ can be either $\zeta_{m}^{(i)}$ or $\zeta_{m}^{(j)}$
with 50\% probability each if $\tau_{m}^{(i)}=\tau_{m}^{(j)}$.

We now specify the payoff function. If player $i$ makes an investment
at time $\tau_{m}^{(i)}$ to boost $X^{\nu}$ by $\zeta_{m}^{(i)}\ge0$,
it costs $c_{i}+k\zeta_{m}^{(i)}$. Based on this cost structure,
we can express the payoff to player $i$ at time $t$ as follows:
\begin{align}
V_{i,t}^{\nu}:= & e^{rt}\hat{\mathbb{E}}_{j}\biggl[\int_{t}^{\infty}\pi(X_{s}^{\nu})e^{-rs}ds-\sum_{m=1}^{\infty}e^{-r\tau_{m}^{(i)}}(k\zeta_{m}^{(i)}+c_{i})\mathbf{1}_{\{\tau_{m}^{(i)}<\tau_{m}^{(j)}\}}\nonumber \\
 & -\frac{1}{2}\sum_{m=1}^{\infty}e^{-r\tau_{m}^{(i)}}(k\zeta_{m}^{(i)}+c_{i})\mathbf{1}_{\{\tau_{m}^{(i)}=\tau_{m}^{(j)}\}}\vert\hat{\mathscr{F}}_{t}\biggr]\;.\label{eq:payoff_general}
\end{align}
Here we let $\hat{\mathbb{E}}_{j}$ denote the expectation over the
product measure $\mathbb{P}\otimes_{m\in\mathbb{N}}\mathbb{L}_{m}^{(j)}$.
We remark that the expression implicitly assumes a pure strategy of
player $i$ for simplicity of presentation. However, the expression
(\ref{eq:payoff_general}) will equally well apply even if $\nu$
is a mixed strategy equilibrium because the payoff to player $i$
can be computed using one of the pure-strategy best responses (\citealt{Hendricks1988},
\citealt{Steg2015}).

In general, it is tricky to account for the possibility of simultaneous
moves by two players in impulse control games. This is because infinitely
many impulse control events may take place between two players at
a single point in time; for instance, the state variable could be
trapped in perpetuity within a bounded interval because the two players
push $X$ back and forth forever within an infinitesimal duration
of time. The prior work in the literature circumvents this problem
in various ways. \citet{Guo2019} excludes simultaneous jump controls
by multiple players from the admissible control set. \citet{Dutta1995}
argues that an infinite number of moves at a single instant would
incur an arbitrarily large cost to the players and illustrates a set
of restrictions on the strategy profile to exclude such a possibility.
\citet{Aid2020} prevents the players from accumulating impulse controls
at some finite time by stipulating a condition that is similar to
the assumption $\lim_{n\rightarrow\infty}\eta_{n}^{(i)}=\infty$ in
our paper. We note that we can also impose a set of restrictions on
the equilibrium strategy profile as \citet{Dutta1995} to exclude
this pathology. However, we do not explicitly impose these conditions
because, in a mixed strategy TSSPE we consider in our paper, a player's
investment always moves the state variable to the common continuation
region, and hence, such a complication does not arise.

\subsection{Verification Theorem \label{subsec:Verification-Theorem}}

In this subsection, we construct a specific strategy profile and payoff
functions, and we prove that they constitute a mixed strategy equilibrium
if they satisfy a set of conditions. 

\textbf{Strategy profile}: In the proposed strategy profile $\nu^{*}$,
each period is divided into two stages as in Section \ref{subsubsec:2-stage-SPE}.
Let $\tau_{\theta^{*}}^{m}:=\inf\{t\in[T_{m},T_{m+1}):X_{t}\le\theta^{*}\}$
denote the hitting time of $(a,\theta^{*})$, where $\theta^{*}$
is some threshold that is different from $\theta_{i}$ or $\theta_{j}$
defined in Section \ref{sec:One-Time}. Then we define the first stage
as the time interval $[T_{m},\tau_{\theta^{*}}^{m})$ and the second
stage as $[\tau_{\theta^{*}}^{m},T_{m+1})$. 

In the first stage, if $X_{T_{m}}^{\nu^{*}}\ge\theta^{*}$, then player
2's strategy is to invest at time $\tau_{\theta^{*}}^{m}$ with a
probability of $q\in(0,1)$. On the other hand, player 1's strategy
is not to invest at all in the first stage. If player 2 invests at
$\tau_{\theta^{*}}^{m}$, then period $m$ ends, and the next period
$m+1$ begins. With a probability of $1-q$, player $2$ does not
invest at $\tau_{\theta^{*}}^{m}$, in which case the second stage
of period $m$ begins. In case $X_{T_{m}}^{\nu^{*}}<\theta^{*}$ in
the beginning of the period, player 2's probability of immediate investment
is set to 0. In other words, if $X_{T_{m}}^{\nu^{*}}<\theta^{*}$,
then the period immediately enters its second stage with probability
one.

In the second stage, the investment strategies of $\nu^{*}$ is characterized
by a common mixed strategy region of $\Gamma:=(a,\theta^{*})$ with
the threshold $\theta^{*}$. Each player $i$ invests with a rate
of $\lambda_{t}^{(i)}=\lambda_{i}(X_{t})$ where $\lambda_{i}(\cdot)$
will be defined in (\ref{eq:lambda-impulse}) below. In either stage
of $\nu^{*}$, player $i$'s magnitude of investment is $\max\{z_{i}-x,0\}$,
where $x$ is the value of $X$ at the time of investment. In general,
$z_{1}$ and $z_{2}$ may differ.

See Figure \ref{fig:Sim_MixedStrategy} for a simulated sample path
of $X^{\nu^{*}}$ in periods 0 \textendash{} 2, where player 2 invests
at the end of periods 0 and 2 while player 1 invests at the end of
period 1. In period 0, player 2 invests exactly at $\tau_{\theta^{*}}^{0}$
while in periods 1 and 2, the players invest in the second stage within
the mixed strategy region $\Gamma$.

\begin{figure}
\centering{}\includegraphics[scale=0.3]{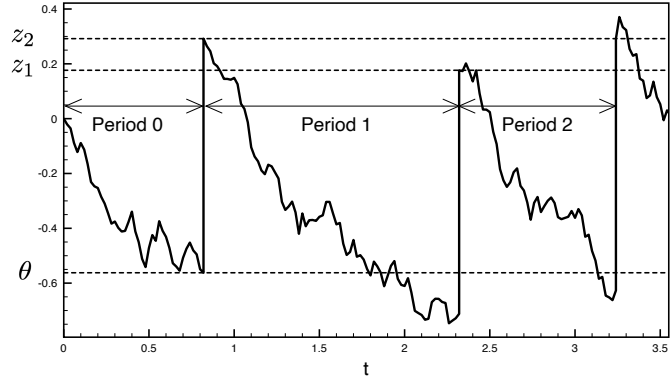}\caption{A simulated sample path of $X^{\nu^{*}}$ under the mixed strategy
profile $\nu^{*}$. \label{fig:Sim_MixedStrategy}}
\end{figure}

\textbf{Payoff functions}: Next, we construct functions $F_{i}(x):=\{F_{i,t}(x)\}_{t\ge0}$
for each $i\in\{1,2\}$, $x\in\mathscr{I}$, and $t\in\mathbb{R}^{+}$.
$F_{i}(x)$ is the candidate for the payoff function for player $i$
associated with $\nu^{*}$. Below we list a number of conditions that
these functions satisfy.

First, we assume that the function $F_{i}(\cdot)$ is expressible
in terms of some functions $U_{1}(\cdot)$ and $V_{i}(\cdot)$ as
follows: 
\begin{align}
F_{1,t}(x) & =\begin{cases}
U_{1}(x) & \text{for}\;t\in[T_{m},\tau_{\theta^{*}}^{m})\\
V_{1}(x) & \text{for}\;t\in[\tau_{\theta^{*}}^{m},T_{m+1})
\end{cases}\;,\label{eq:F1}\\
F_{2,t}(x) & =V_{2}(x)\qquad.\nonumber 
\end{align}

We assume that the game begins at $t=0$ in the first stage of the
0-th period so that $F_{1,0}(x):=U_{1}(x)$ and $F_{2,0}(x):=V_{2}(x)$.
Here we stipulate that $U_{1}(\cdot)$ and $V_{i}(\cdot)$ satisfy
$U_{1}(\cdot)\in C^{2}(\mathscr{I}\backslash\{\theta^{*}\})\cap C^{1}(\mathscr{I}\backslash\{\theta^{*}\})\cap C(\mathscr{I}\backslash\{\theta^{*}\})$,
$V_{i}(\cdot)\in C^{2}(\mathscr{I}\backslash\{\theta^{*}\})\cap C^{1}(\mathscr{I})\cap C(\mathscr{I})$,
and $V_{i}(\cdot)>(R_{r}\pi)(\cdot)$. We also stipulate that $U_{1}(x)>V_{1}(x)$
for $x\ge\theta^{*}$ and $U_{1}(x)=V_{1}(x)$ for $x<\theta^{*}$.
Note that $F_{1,t}(x)$ is discontinuous in time at $\tau_{\theta^{*}}^{m}$
if $x\ge\theta^{*}$; this is to reflect that the first and second
stage payoffs to player 1 differ from each other because player 2
invests at $\tau_{\theta^{*}}^{m}$ with a probability of $q$.

Second, we stipulate that the functions $U_{1}(\cdot)$ and $V_{i}(\cdot)$
satisfy variational inequalities given by
\begin{align}
\mathscr{A}V_{i}(x)+\pi(x) & =0\;\text{for}\;x>\theta^{*}\:,\label{eq:AV_i}\\
\mathscr{A}V_{i}(x)+\pi(x) & <0\;\text{for}\;x<\theta^{*}\;,\label{eq:AV_i-ineq}\\
\mathscr{A}U_{1}(x)+\pi(x) & =0\;\text{for}\;x>\theta^{*}\:.\label{eq:AU_1}
\end{align}
In addition, we impose the quasi-variational inequalities for all
$x\in\mathscr{I}$: $F_{i,t}(x)\ge\sup_{\zeta\ge0}[F_{i,0}(x+\zeta)-k\zeta-c_{i}]$.
Since $F_{i,0}(\cdot)$ is the first-stage functions, this inequality
translates into the following set of inequalities:
\begin{align}
V_{2}(x) & \ge\sup_{\zeta\ge0}[V_{2}(x+\zeta)-k\zeta-c_{2}]\:,\label{eq:quasi-V2}\\
V_{1}(x) & \ge\sup_{\zeta\ge0}[U_{1}(x+\zeta)-k\zeta-c_{1}]\:,\nonumber \\
U_{1}(x) & \ge\sup_{\zeta\ge0}[U_{1}(x+\zeta)-k\zeta-c_{1}]\:.\nonumber 
\end{align}
We also assume that $\vert U_{1}^{\prime\prime}(x)\vert<\infty$ and
$\vert V_{i}^{\prime\prime}(x)\vert<\infty$ are satisfied for $x\not=\theta^{*}$.

Third, we relate $z_{i}$ and the investment rate $\lambda_{i}(\cdot)$
with the functions $F_{i,0}(\cdot)$. We assume that $z_{i}$ is the
unique value that satisfies $z_{i}=\arg\max_{z}[F_{i,0}(z)-kz]$,
which has the meaning of the optimal end point of the boost for player
$i$. Then we can define 
\begin{equation}
g_{i}(x):=\begin{cases}
F_{i,0}(z_{i})-k(z_{i}-x)-c_{i} & \text{for}\;x<z_{i}\\
F_{i,0}(x)-c_{i} & \text{for}\;x\ge z_{i}
\end{cases}\:,\label{eq:g_i}
\end{equation}
which has the interpretation of the reward from investment. We also
stipulate that $V_{i}(x)=g_{i}(x)$ for $x\in\Gamma=(a,\theta^{*})$,
which is reminiscent of the property of a mixed strategy equilibrium.
Finally, we define the investment rate as:
\begin{equation}
\lambda_{i}(x)=\mathbf{1}_{\{x\in(a,\theta^{*})\}}\frac{-\mathscr{A}g_{j}(x)-\pi(x)}{F_{j,0}(z_{i})-g_{j}(x)}\;.\label{eq:lambda-impulse}
\end{equation}

Based on the imposed conditions, we can construct an explicit functional
forms of $U_{1}(\cdot)$ and $V_{i}(\cdot)$. The equations $\mathscr{A}V_{i}(x)+\pi(x)=0$
in (\ref{eq:AV_i}) and $\mathscr{A}U_{1}(x)+\pi(x)=0$ in (\ref{eq:AU_1})
imply the following functional forms: 
\begin{align}
V_{i}(x) & =\begin{cases}
(R_{r}\pi)(x)+w_{i}\phi(x) & \text{for}\quad x\ge\theta^{*}\\
g_{i}(x) & \text{for}\quad x<\theta^{*}
\end{cases}\:,\label{eq:Vi-explicit}\\
U_{1}(x) & =\begin{cases}
(R_{r}\pi)(x)+u_{1}\phi(x) & \text{for}\quad x\ge\theta^{*}\\
g_{1}(x) & \text{for}\quad x<\theta^{*}
\end{cases}\:,\nonumber 
\end{align}
for some positive coefficients $w_{i}$ and $u_{1}$. The coefficients
$w_{i}$ and $u_{1}$ must satisfy the condition $u_{1}>w_{1}$ (from
$U_{1}(x)>V_{1}(x)$ for $x\ge\theta^{*}$). Lastly, we impose two
more conditions that relate $\nu^{*}$ to $U_{1}$ and $V_{i}$. The
first is the boundary condition $U_{1}(\theta^{*})=qU_{1}(z_{2})+(1-q)V_{1}(\theta^{*})$
which we stipulate because player $2$ invests at the hitting time
of $\theta^{*}$ with a probability of $q$. The second condition
is $U_{1}(z_{2})\ge V_{1}(\theta^{*})=g_{1}(\theta^{*})$, which has
the interpretation that player 1 earns higher payoff when player 2
invests than when player 1 invests at $\theta^{*}$.

\textbf{Survival probability}: Based on the strategy profile $\nu^{*}$,
we can construct per-period survival probability $M_{m,t}^{(i)}$
as follows:
\begin{align}
M_{m,t}^{(1)} & =\exp[-\mathbf{1}_{\{t>\tau_{\theta^{*}}^{m}\}}\int_{\tau_{\theta^{*}}^{m}}^{t}\lambda_{1}(X_{t})dt]\:,\label{eq:Multiplicative-1}\\
M_{m,t}^{(2)} & =\begin{cases}
1 & \text{if}\:t<\tau_{\theta^{*}}^{m}\\
(1-\mathbf{1}_{\{X_{\tau_{\theta^{*}}^{m}}=\theta^{*}\}}q)\exp[-\int_{\tau_{\theta^{*}}^{m}}^{t}\lambda_{2}(X_{t})dt] & \text{if}\:t\ge\tau_{\theta^{*}}^{m}
\end{cases}\;.\label{eq:Multiplicative-2}
\end{align}
The expressions $M_{m,t}^{(i)}$ are analogous to (\ref{eq:M1}) and
(\ref{eq:M2}) from the single investment game. The expression (\ref{eq:Multiplicative-1})
clearly follows the strategy $\nu_{1}^{*}$ of investment at the hazard
rate of $\lambda_{1}(X_{t})$. The expression (\ref{eq:Multiplicative-2})
requires some explanation. Player 2 does not invest until $\tau_{\theta^{*}}^{m}$,
so $M_{m,t}^{(2)}$ does not vary in time for $t<\tau_{\theta^{*}}^{m}$.
At the hitting time of $\theta^{*}$, if $X_{\tau_{\theta^{*}}^{m}}=\theta^{*}$,
player 2 invests with a probability of $q$, so $M_{m,t}^{(2)}$ is
reduced by $q$ at time $\tau_{\theta^{*}}^{m}$. However, if $X_{\tau_{\theta^{*}}^{m}}<\theta^{*}$,
then the game immediately enters the second stage, so there is no
immediate discontinuous reduction in $M_{m,t}^{(2)}$. After $\tau_{\theta^{*}}^{m}$,
player 2's rate of investment is $\lambda_{2}(X_{t})$, which explains
the factor $\exp[-\int_{\tau_{\theta^{*}}^{m}}^{t}\lambda_{2}(X_{t})dt]$.

\medskip{}

In summary, we have specified a strategy profile $\nu^{*}$ and hypothesized
a pair of functions $F_{1}(\cdot)$ and $F_{2}(\cdot)$ that satisfy
a set of conditions. We are now ready to present Theorem \ref{thm:nu*-equilibrium}
which establishes that $\nu^{*}$ is an equilibrium given the existence
of $F_{1}(\cdot)$ and $F_{2}(\cdot)$.

\begin{theorem} \label{thm:nu*-equilibrium} Under Assumptions \ref{assum:integrability}
and \ref{assum:single-inv}, if $F_{1}(\cdot)$ and $F_{2}(\cdot)$
that satisfy the conditions above exist, then $\nu^{*}$ is an SPE.
\end{theorem}

Theorem \ref{thm:nu*-equilibrium} establishes a sufficient condition
for a mixed strategy SPE. To our knowledge, this is the first verification
theorem in the literature for a mixed strategy equilibrium of an impulse
control game.

\subsection{Existence of Mixed Strategy TSSPE \label{subsec:Existence-of-Mixed}}

In this section, we utilize Theorem \ref{thm:nu*-equilibrium} and
show that a mixed strategy TSSPE exists in our model under certain
conditions. Closely following the conventions used in \citet{Alvarez2008},
we define the following auxiliary functions associated with the uncontrolled
process $X$ and the flow profit $\pi(\cdot)$:
\begin{align*}
S'(x):=\exp\left[-\int_{x_{0}}^{x}\frac{2\mu(y)}{\sigma^{2}(y)}dy\right]\:,\: & m'(x):=\frac{2}{\sigma^{2}(x)S'(x)}\:,\\
\rho(x):=\pi(x)+k[\mu(x)-rx]\:,\: & L(x):=-\frac{\rho(x)\phi'(x)}{S'(x)}-r\int_{x}^{b}\phi(y)\rho(y)m'(y)dy\;,
\end{align*}
where $x_{0}\in\mathscr{I}$ is an arbitrarily chosen reference point,
$S(\cdot)$ is the \emph{scale function}, and $m(\cdot)$ is the \emph{speed
measure} (II.4 in \citealp{Borodin1996}).

We now make the following assumptions to ensure that the desired equilibrium
exists.

\begin{assumption} \label{assum:rho_L} (i) $\rho(\cdot)\in C^{1}(\mathscr{I})$,
and there exists $x^{*}\in\mathscr{I}$ such that $\rho'(x)>0$ for
$x<x^{*}$ and $\rho'(x)<0$ for $x>x^{*}$. (ii) $\lim_{x\rightarrow b}\frac{\rho(x)\phi'(x)}{S'(x)}=0$.
(iii) $L(x)<0$ for some $x\in\mathscr{I}$.

\end{assumption}

Assumption \ref{assum:rho_L} can be found to be satisfied in many
examples of stochastic impulse control problems. Assumption \ref{assum:rho_L}(i)
is made by \citet{Alvarez2008}, and Assumption \ref{assum:rho_L}(ii)
is only slightly stronger than the condition $\lim_{x\rightarrow b}\frac{\phi'(x)}{S'(x)}=0$
which must be satisfied by a natural boundary $b$. Assumption \ref{assum:rho_L}(iii)
is needed to ensure that the function $L(\cdot)$ has a unique zero
in $\mathscr{I}$, which is used to derive an optimal impulse control
in \citet{Alvarez2008}.

In addition, we define two more auxiliary functions:
\begin{align}
I(x) & :=\frac{(R_{r}\pi)'(x)-k}{\phi'(x)}\label{eq:I-exp}\\
J(x) & :=(R_{r}\pi)(x)-kx-I(x)\phi(x)\:.\label{eq:J-exp}
\end{align}
Based on Assumptions \ref{assum:integrability}, \ref{assum:single-inv},
and \ref{assum:rho_L}, we obtain the following useful properties
of $I(\cdot)$ and $J(\cdot)$:

\begin{lemma} \label{lemm:I0-J0} (i) There exists $\hat{x}<\min\{x^{*},z^{*}\}$
such that $I'(x)<0$ and $J'(x)>0$ for $x<\hat{x}$, and $I'(x)>0$
and $J'(x)<0$ for $x>\hat{x}$.

(ii) $I(x)<0$ for $x<z^{*}$ and $I(x)>0$ for $x>z^{*}$. \end{lemma}

Lastly, we make the following additional assumption:

\begin{assumption} \label{assump:single-player-impulse} There exist
$\theta^{*}$ and $z_{2}$ that satisfy
\begin{align}
I(z_{2}) & =I(\theta^{*})\label{eq:I-z2-th}\\
J(z_{2})-J(\theta^{*}) & =c_{2}\:.\label{eq:J-c2}
\end{align}
\end{assumption}

In the context of our model, Assumption \ref{assump:single-player-impulse}
can be shown to ensure the optimal solution to the single-player impulse
control problem (Section 4.2 of \citealt{Alvarez2008}; Lemma \ref{lemm:1-player-impulse}
of our paper). Thus, this is not a particularly stringent assumption
because we are essentially stipulating that the optimal solution to
the single-player impulse control problem exists.

Finally, under Assumptions \ref{assum:integrability}, \ref{assum:single-inv},
\ref{assum:rho_L}, and \ref{assump:single-player-impulse}, we establish
the following theorem:

\begin{theorem}\label{thm:Existence_SPE} A mixed strategy TSSPE
$\nu^{*}$ characterized by $\theta^{*}$, $z_{1}$, $z_{2}$, and
$q\in(0,1)$ exists if $c_{1}\in(c_{2},\hat{c})$ for some $\hat{c}>c_{2}$.
Furthermore, the parameters satisfy $\theta^{*}<\hat{x}<z_{1}<z_{2}$.
\end{theorem}

Theorem \ref{thm:Existence_SPE} stands in stark contrast to the result
from the single investment game. In the model that allows for an infinite
number of investment opportunities, a moderate amount of asymmetry
does not destabilize the mixed strategy SPE. In our impulse control
game, the mixed strategy equilibrium can be sustained despite asymmetry
because player 2's active investment strategy with $q>0$ in the future
periods increases player 1's payoff in the current period. To elaborate,
recall that in the single investment game, a mixed strategy MPE is
destabilized by asymmetry because the mixed strategy regions $\Gamma_{1}$
and $\Gamma_{2}$ cannot coincide due to asymmetry. Namely, because
$c_{1}>c_{2}$, player 1 naturally has lower threshold of investment.
In the repeated game of investment, however, player 2 can allocate
a sufficiently high probability $q$ of investment at the hitting
time of $\theta^{*}$, which then drives player 1's investment threshold
to $\theta^{*}$. For the same reason, $q$ increases in $c_{1}$
because, as $c_{1}$ increases, $q$ needs to increase to compensate
for player 1's disincentive to invest.

Theorem \ref{thm:Existence_SPE} also establishes that player 1's
magnitude of the investment is comparatively less than that of player
2, i.e., $z_{1}<z_{2}$. Because of the probability $q>0$ that player
2 would invest at the hitting time of $\theta^{*}$, player 1's payoff
function is particularly high when $X^{\nu^{*}}$ is close to the
threshold $\theta^{*}$, but the effect of player 2's investment at
$\theta^{*}$ tapers off as $X^{\nu^{*}}$ moves far above $\theta^{*}$.
Thus, player 1's incentive to boost $X^{\nu^{*}}$ to a high value
is not as strong as that of player 2, and hence, $z_{1}$ is less
than $z_{2}$.

\subsection{Example \label{subsec:Example}}

We now illustrate an example of the mixed strategy TSSPE. As the uncontrolled
state variable process, we consider a geometric Brownian motion with
SDE $dX_{t}=\mu X_{t}dt+\sigma X_{t}dW_{t}$ for some constants $\mu<0$
and $\sigma>0$. In this case, $\mathscr{I}=(0,\infty)$, and the
solutions to the homogeneous differential equation $\mathscr{A}f=0$
are given by $\phi(x)=x^{\gamma_{-}}$ and $\psi(x)=x^{\gamma_{+}}$,
where the power indices are given by 
\[
\gamma_{\pm}=\frac{1}{2}-\frac{\mu}{\sigma^{2}}\pm\sqrt{(\frac{1}{2}-\frac{\mu}{\sigma^{2}})^{2}+\frac{2r}{\sigma^{2}}}\;.
\]
We assume a flow profit function of $\pi(x)=x^{\alpha}$ for some
$\alpha\in(0,1)$. Then it follows that $(R_{r}\pi)(x)=x^{\alpha}/(r-\delta(\alpha))$
where $\delta(\alpha)=\alpha\mu+\sigma^{2}\alpha(\alpha-1)/2<0$ \citep{Alvarez2008}.

It is straightforward to verify that this example satisfies Assumptions
\ref{assum:integrability}, \ref{assum:single-inv}, and \ref{assum:rho_L}.
In addition, the following proposition ensures that Assumption \ref{assump:single-player-impulse}
holds if $c_{2}$ is not too large.

\begin{proposition}\label{prop:NumericalEx} The example above satisfies
Assumption \ref{assump:single-player-impulse} if $c_{2}\in(0,J(z^{*}))$.\end{proposition}

By Proposition \ref{prop:NumericalEx} and Theorem \ref{thm:Existence_SPE},
this example possesses a mixed strategy TSSPE if $c_{2}\in(0,J(z^{*}))$
and the difference $c_{1}-c_{2}$ is not too large. If we set $r=1,\alpha=0.5,\mu=-0.5,\sigma=0.25,k=1,c_{2}=0.015$,
then we obtain $\theta^{*}=0.0439$ and $z_{2}=0.1480$. Under this
parameter set, a mixed strategy TSSPE $\nu^{*}$ exists as long as
$c_{1}\le0.01729$.\footnote{If $c_{1}>0.01729$, condition (\ref{eq:OptimalU}) of EC Appendix
fails to hold, so the mixed strategy TSSPE does not exist. Note that
condition (\ref{eq:OptimalU}) of EC Appendix assures that there is
a unique optimal value of the end point $z_{1}$ of the boost; if
this condition is violated, then player 1 may improve his payoff by
choosing a different point $z_{1}'$ closer to $\theta^{*}$ as the
end point, which thus destabilizes the equilibrium.}

\subsection{Comparison to Pure Strategy Equilibrium \label{subsec:Comparison-to-Pure}}

In this subsection, we demonstrate that there exists a pure strategy
equilibrium which is more efficient than the mixed strategy SPE $\nu^{*}$.
Letting $V_{i}^{\nu}(x)$ denote the payoff function associated with
a strategy profile $\nu$ conditional on $X_{0}=x$, we say that an
equilibrium $\nu$ is \emph{more efficient} than another equilibrium
$\nu'$ if $V_{1}^{\nu}(x)+V_{2}^{\nu}(x)>V_{1}^{\nu'}(x)+V_{2}^{\nu'}(x)$
for all $x\in\mathscr{I}$.

As a preliminary step, we establish the solution to the single-player
model, which reduces to the standard impulse control problem \citep{Alvarez2008}.
For convenience, we let $c$ denote the upfront cost of investment
by the single player. The following lemma directly follows from \citet{Alvarez2008}:

\begin{lemma} \label{lemm:1-player-impulse} Under Assumptions \ref{assum:integrability},
\ref{assum:single-inv}, and \ref{assum:rho_L}, suppose that there
exist $\theta^{*}$ and $z$ such that $\theta^{*}<z$, $I(\theta^{*})=I(z)$,
and $J(z)-J(\theta^{*})=c$. In each period, the optimal policy $\nu$
of the single player is to invest at $\tau_{\theta^{*}}=\inf\{t\ge0:X_{t}^{\nu}\le\theta^{*}\}$
to boost $X^{\nu}$ up to $z$. Furthermore, the payoff function is
given by 
\begin{equation}
V_{s}(x)=\begin{cases}
-I(\theta^{*})\phi(x)+(R_{r}\pi)(x) & x>\theta^{*}\\
V_{s}(z)-k(z-x)-c & x\le\theta^{*}
\end{cases}\ .\label{eq:V-nu}
\end{equation}
\end{lemma}

Next, we establish a pure strategy equilibrium which is more efficient
than a mixed strategy equilibrium. We again assume that $c_{1}\ge c_{2}$
and let $z_{2}$ and $\theta^{*}$ denote the solution to $I(\theta^{*})=I(z_{2})$
and $J(z_{2})-J(\theta^{*})=c_{2}$. Recall that $\hat{x}$ is the
global minimizer of $I(\cdot)$.

\begin{proposition}\label{prop:pure-mixed} Suppose that 
\begin{equation}
\beta:=\frac{(R_{r}\pi)(z_{2})-(R_{r}\pi)(\theta^{*})}{\phi(\theta^{*})-\phi(z_{2})}>-I(\hat{x})\label{eq:pure-suff}
\end{equation}
holds. Then there exists a pure strategy equilibrium $\nu^{P}$ in
which only player 2 invests. Furthermore, $\nu^{P}$ is more efficient
than the mixed strategy SPE $\nu^{*}$ obtained in Theorem \ref{thm:Existence_SPE}.
\end{proposition}

The condition (\ref{eq:pure-suff}) ensures that the variational inequality
holds for the pure strategy equilibrium in which player 1 never invests
and player 2 is the only player who invests. In the numerical example
illustrated in Section \ref{subsec:Existence-of-Mixed}, we confirm
that (\ref{eq:pure-suff}) holds: we obtain $\beta=7.60\times10^{-4}>-I(\hat{x})=2.95\times10^{-4}$.
It is not easy to verify whether (\ref{eq:pure-suff}) generally holds,
but it is beyond the scope of this paper to conduct an exhaustive
search for pure strategy equilibria.

Proposition \ref{prop:pure-mixed} establishes that the pure strategy
equilibrium is relatively more efficient. Intuitively, the investments
are never delayed in a pure strategy equilibrium, so it is natural
that this equilibrium is more efficient than the mixed strategy equilibrium.
Thus, it behooves the policymaker or the social planner to attempt
to thwart a mixed strategy equilibrium in favor of a pure strategy
equilibrium.

To further elaborate on this point, we compare the payoffs associated
with the mixed strategy equilibrium, the pure strategy equilibrium,
and the social planner's solution. For this purpose, we need to formulate
and solve the social planner's problem whose objective is to maximize
the sum of the two players' payoff functions. We first observe that
the social planner should always have player 2 invest because $c_{2}<c_{1}$.
Thus, the upfront cost of investment is effectively $c_{2}$. Secondly,
the total profit flow per unit time is $2\pi(\cdot)$ because it is
the sum of the two player's profit flows. Hence, given a social planner's
optimal policy $\nu=(\tau_{n},\zeta_{n})_{n\in\mathbb{N}}$ of impulse
control on $X$, the total payoff function is defined as follows:
\[
V^{S}(x):=\sup_{\nu}\mathbb{E}^{x}[\int_{0}^{\infty}2\pi(X_{t}^{\nu})e^{-rt}dt-\sum_{n=1}^{\infty}e^{-r\tau_{n}}(k\zeta_{n}+c_{2})]\:.
\]
 Now suppose that the conditions for Lemma \ref{lemm:1-player-impulse}
hold with the modified profit flow function $2\pi(\cdot)$ in place
of $\pi(\cdot)$ and $c_{2}$ in place of $c$. We let $z_{S}$ and
$\theta_{S}$ denote the solutions to $I(\theta_{S})=I(z_{S})$, and
$J(z_{S})-J(\theta_{S})=c$ where $I(\cdot)$ and $J(\cdot)$ are
appropriately modified to accommodate the profit flow $2\pi(\cdot)$.
Then the social planner's optimal policy is to boost $X$ up to $z_{S}$
whenever $X$ falls below $\theta_{S}$.

We also define $V^{P}(x):=V_{1}^{\nu^{P}}(x)+V_{2}^{\nu^{P}}(x)$
and $V^{M}(x):=V_{1}(x)+V_{2}(x)$ respectively as the total payoffs
(to the social planner) associated with the pure strategy equilibrium
and the mixed strategy equilibrium. For simplicity, we define $V^{M}$
as the total payoff of the second stage, but the qualitative features
of the first-stage payoff $U_{1}(x)+V_{2}(x)$ are similar.

To compare the three kinds of total payoffs, we revisit the numerical
example given in Section \ref{subsec:Existence-of-Mixed} with $c_{1}=0.0165$
and illustrate the payoffs in Figure \ref{fig:MPS_value}. Not surprisingly,
the social planner's optimal total payoff is higher than the payoffs
associated with either equilibrium. Most importantly, $V^{M}$ is
strictly less than $V^{P}$, which confirms Proposition \ref{prop:pure-mixed}.

\begin{figure}
\centering{}\includegraphics[scale=0.5]{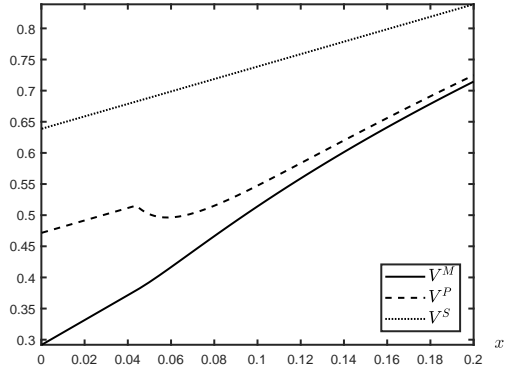}\caption{The total payoffs associated with the mixed strategy equilibrium,
the pure strategy equilibrium, and the social planner's solution.
\label{fig:MPS_value}}
\end{figure}

\section{Conclusions \label{sec:Conclusions}}

The recent theoretical result establishes that asymmetry destabilizes
mixed strategy MPE in stochastic wars of attrition. This result is
vexing because the extant literature views the mixed strategy equilibrium
as the hallmark of the war of attrition. From the practical point
of view, the mixed strategy equilibrium is regarded as the realistic
outcome of a concession game. Thus, the non-existence result begs
the question of whether the war of attrition model is an adequate
framework for games with a free rider problem. In this paper, we establish
that repeated opportunities of concession stabilizes the mixed strategy
equilibrium despite asymmetry. Therefore, we offer one possible resolution
of the apparent discrepancy between the non-existence result from
the stochastic war of attrition and the conventional view.

The result of this paper has practical implications to stakeholders
of common goods. First, the free rider effect from the mixed strategy
equilibrium is much more robust to asymmetry than suggested by the
recent theoretical results on the stochastic war of attrition. Second,
in light of our result that the mixed strategy equilibrium disappears
under a sufficiently high degree of asymmetry, it may behoove policymakers
to promote a sufficiently high degree of asymmetry to mitigate the
inefficiency arising from the mixed strategy equilibrium. For instance,
in the context of CSR, the government may induce sufficiently high
asymmetry between the firms through asymmetric incentives such as
selective subsidies.\footnote{For the body of work that examines the efficacy of selective subsidy
to facilitate investments in public goods, see \citet{Ayres1996},
\citet{Rothkopf2003}, and \citet{Kotowski2018}. }

This paper is potentially useful to future research on free rider
problems. First, it will be fruitful to investigate whether the mixed
strategy equilibrium can be stabilized by any other realistic features
of a stochastic concession game. Second, the verification theorem
constructed in this paper may be extended to other stochastic control
games with free rider problems which might possess mixed strategy
equilibria.

\theendnotes
\noindent \begin{flushleft}
\textsf{\textbf{\large{}References}}{\large\par}
\par\end{flushleft}

{\small{}\begin{btSect}[INFORMS2011]{ImpulseGame}
\btPrintCited
\end{btSect}
}{\small\par}

\end{btUnit}

\newpage{}

\appendix
\begin{btUnit}
\setcounter{page}{1}

\section*{{\LARGE{}Electronic Companion}}

\section{Case of Singular Rates of Investment \label{App_sec:Singular-Rates-of}}

In this Appendix, we consider the possibility that $A_{t}^{(i)}$
is not differentiable with respect to time because its Radon-Nikodym
derivative $\lambda_{i}(\cdot)$ has singularities, and we argue that
this case can be approximated by finite but arbitrarily high rates
of investment. One such example is a local time process denoted by
$L_{t}(x)$ at some point $x\in\mathscr{I}$ because it is well-known
that $L_{t}(x)$ is continuous but does not possess a time-derivative
at every point in time.

By (6.2) on p. 202 of \citet{Karatzas1998}, 
\[
L_{t}(x)=\lim_{\lambda\rightarrow\infty}\frac{\lambda}{4}\mathbb{L}(\{0\le s\le t:\vert W_{s}-x\vert\le\lambda^{-1}\})\:,
\]
where $\mathbb{L}$ is the Lebesgue measure over time. The right hand
side is the total duration of time that a Brownian motion $W_{s}$
stays within the interval $[x-\lambda^{-1},x+\lambda^{-1}]$ as $\lambda\rightarrow\infty$.

Assuming that $X=W$, if this process were to be added to $A^{(i)}$,
it implies that player $i$ invests at a hazard rate of $\lambda$
within the interval $[x-\lambda^{-1},x+\lambda^{-1}]$ for a very
large value of $\lambda$. Therefore, $L_{t}(\cdot)$ is a limiting
case of arbitrarily large values of the rate of investment within
an arbitrarily small interval.

We remark that this argument generally applies to any additive functionals
of regular diffusion processes. Specifically, an arbitrary additive
functional $A_{t}$ of a regular diffusion process has the following
representation:
\[
A_{t}=\int_{0}^{t}\lambda(X_{s})ds\:,
\]
where $\lambda(\cdot)$ is the \emph{killing measure} density that
characterizes $A_{t}$ (II.23 and II.24, \citealt{Borodin1996}).
The killing measure density $\lambda(\cdot)$ may or may not have
singularities, but if it is finite everywhere, then $A_{t}$ is time-differentiable.
Even if $\lambda(\cdot)$ has a singularity (such as in a Dirac measure)
at some point $x_{0}$, one can approximate it as a finite function
taking very high values near $x_{0}$.

\section{Mathematical Proofs\label{App_sec:Mathematical-Proofs}}

\textbf{Proof of Lemma \ref{lemma:optimal-zeta}}: By the definition
of SPE, the equilibrium value of $\zeta_{t}^{(i)}$ must be the one
that maximizes the payoff to player $i$ at the time of investment
$t$. We thus examine player $i$'s payoff at its time of investment.
Towards this goal, we first formulate the time of investment in case
$\tau^{(i)}=\tau^{(j)}$ by utilizing the following mathematical machinery.
We introduce a Bernoulli random variable $Z_{i}\sim\text{Ber}(\frac{1}{2})$
that is independent of the uncontrolled state variable $X$ such that
$Z_{1}+Z_{2}=1$. We stipulate that the value of $Z_{i}$ is revealed
at the stopping time $\tau^{(i)}$ only in the event $\tau^{(i)}=\tau^{(j)}$,
in which case player $i$ becomes an investor if and only if $Z_{i}=1$,
which occurs with probability 50\%. Then we define an augmented probability
space $(\Omega\times\{0,1\},\tilde{\mathscr{F}},\tilde{\mathbb{F}},\tilde{\mathbb{P}})$
where $\tilde{\mathbb{F}}$ is the filtration generated by $X$ and
$Z_{i}$. Furthermore, the new probability measure is given by $\tilde{\mathbb{P}}=\mathbb{P}\otimes\mathbb{Z}$,
where $\mathbb{Z}$ is the probability measure of $Z_{i}$. Under
the new probability space, we extend the stopping time $\tau^{(i)}$,
which is adapted to $\mathbb{F}$, to a stopping time $\tilde{\tau}^{(i)}$
adapted to $\tilde{\mathbb{F}}$ as follows: 
\begin{align*}
\tilde{\tau}^{(i)} & =\tau^{(i)}\;\text{if }\:\tau^{(i)}\neq\tau^{(j)}\\
 & =\tau^{(i)}\;\text{if }\:\tau^{(i)}=\tau^{(j)}\;\text{and}\;Z_{i}=1\\
 & =\infty\;\text{if}\:\;\tau^{(i)}=\tau^{(j)}\;\text{and}\;Z_{i}=0\;.
\end{align*}
By this formulation, player $i$ invests at the stopping time $\tilde{\tau}^{(i)}$
only if $\tilde{\tau}^{(i)}<\tilde{\tau}^{(j)}$; the event $\tilde{\tau}^{(i)}=\tilde{\tau}^{(j)}$
never takes place by construction.

Next, we examine player $i$'s payoff at its investment time $\tilde{\tau}^{(i)}$
conditional on that player $i$ is an investor at time $\tilde{\tau}^{(i)}<\infty$:
\begin{align}
 & \mathbb{E}^{x}[\int_{\tilde{\tau}^{(i)}}^{\infty}\pi(X_{s}^{\nu})e^{-rs}ds-(k\zeta_{\tilde{\tau}^{(i)}}^{(i)}+c_{i})e^{-r\tilde{\tau}^{(i)}}|\mathscr{F}_{\tilde{\tau}^{(i)}},\tilde{\tau}^{(i)}<\tilde{\tau}^{(j)}]\nonumber \\
 & =e^{-r\tilde{\tau}^{(i)}}\{\mathbb{E}^{X_{\tilde{\tau}^{(i)}}+\zeta_{\tilde{\tau}^{(i)}}^{(i)}}[\int_{0}^{\infty}\pi(X_{s})e^{-rs}ds]-(k\zeta_{\tilde{\tau}^{(i)}}^{(i)}+c_{i})\}=e^{-r\tilde{\tau}^{(i)}}\{(R_{r}\pi)(z)-k(z-X_{\tilde{\tau}^{(i)}})-c_{i}\}\:,\label{eq:pf-payoff}
\end{align}
where $z\equiv\zeta_{\tilde{\tau}^{(i)}}^{(i)}+X_{\tilde{\tau}^{(i)}}$.
Here we use the fact that $X_{s}^{\nu}$ satisfies the SDE (\ref{eq:SDE})
for all $s\geq\tilde{\tau}^{(i)}$ with the condition $X_{\tilde{\tau}^{(i)}}^{\nu}=X_{\tilde{\tau}^{(i)}}+\zeta_{\tilde{\tau}^{(i)}}^{(i)}$
in the single investment game. Then by Assumption \ref{assum:single-inv},
because $\zeta_{\tilde{\tau}^{(i)}}^{(i)}\ge0$ must be satisfied,
the payoff (\ref{eq:pf-payoff}) is uniquely maximized by $z=z^{*}$
if $X_{\tilde{\tau}^{(i)}}<z^{*}$, and it is maximized by $z=X_{\tilde{\tau}^{(i)}}$
if $X_{\tilde{\tau}^{(i)}}\geq z^{*}$. It thus follows that $\zeta_{t}^{(i)}=z^{*}-X_{t}$
if $X_{t}<z^{*}$ and $\zeta_{t}^{(i)}=0$ if $X_{t}\ge z^{*}$ for
any time $t$. \eproof

\textbf{Proof of Proposition \ref{prop:stopping-game}}: We inspect
the payoff function given by (\ref{eq:single-game-payoff}) and replace
$\zeta_{t}^{(i)}$ and $\zeta_{t}^{(j)}$ with the SPE investment
strategies given by (\ref{eq:SPE-zeta}). In the right-hand side of
(\ref{eq:single-game-payoff}), we can decompose $\int_{t}^{\infty}\pi(X_{s}^{\nu})e^{-rs}ds$
as $\int_{t}^{\tau^{(i)}\wedge\tau^{(j)}}\pi(X_{s})e^{-rs}ds+\int_{\tau^{(i)}\wedge\tau^{(j)}}^{\infty}\pi(X_{s}^{\nu})e^{-rs}ds$
where we use the fact that the controlled process $X_{s}^{\nu}$ coincides
with the uncontrolled process $X$ for $s<\tau^{\nu}=\tau^{(i)}\wedge\tau^{(j)}$.
Furthermore, $X_{s}^{\nu}$ satisfies the SDE (\ref{eq:SDE}) for
$s\ge\tau^{\nu}$ with the condition that $X_{\tau^{\nu}}^{\nu}=\max\{z^{*},X_{\tau^{\nu}}\}$.
Hence, we can utilize the strong Markov property of the uncontrolled
state variable $X$ and replace $\mathbb{E}^{x}[\int_{\tau^{\nu}}^{\infty}\pi(X_{s}^{\nu})e^{-rs}ds\vert\mathscr{F}_{\tau^{\nu}}]$
with $e^{-r\tau^{\nu}}(R_{r}\pi)(X_{\tau^{\nu}}^{\nu})$.  Then
we can use the definitions of $g_{i}(\cdot)$ and $m_{i}(\cdot)$
given in (\ref{eq:gi}) and (\ref{eq:mi}) to rewrite (\ref{eq:single-game-payoff})
as (\ref{eq:V-nu-stopping-time}). \eproof

\textbf{Proof of Proposition \ref{prop:single-p}}. (i) By Proposition
\ref{prop:stopping-game}, it suffices to consider an optimal stopping
problem of the form (\ref{eq:single-p}). We first note that it is
never optimal to invest when $X_{t}\ge z^{*}$. If a player invests
when $X_{t}\ge z^{*}$, the optimal boost is zero, and hence, the
player would only waste the upfront cost $c_{i}$. Thus, we limit
our attention to the policy of investment when $X_{t}<z^{*}$.

Note that we can re-write the payoff function (\ref{eq:single-p})
as follows:
\[
\mathbb{E}^{x}\{(R_{r}\pi)(x)+e^{-r\tau}[g_{i}(X_{\tau})-(R_{r}\pi)(X_{\tau})\}\:.
\]
Thus, we have transformed the problem into an optimal stopping problem
with a bequest function $g_{i}(x)-(R_{r}\pi)(x)$. By definition of
$g_{i}(\cdot)$ in (\ref{eq:gi}), $g_{i}(x)-(R_{r}\pi)(x)$ decreases
in $x$ for $x<z^{*}$ and stays constant for $x\ge z^{*}$. Under
Assumption \ref{assump:theta}, we can invoke Theorem 3(B) of \citet{Alvarez2001}
to conclude that $\tau^{*}=\inf\{t\ge0:X_{t}\le\theta_{i}\}$.

(ii) This statement directly follows from Theorem 3(B) of \citet{Alvarez2001}.

(iii) Lastly, we prove that $\theta_{i}$ decreases in $c_{i}$. Define
$f(x,c_{i})=[g_{i}(x)-(R_{r}\pi)(x)]/\phi(x)$. For $x<z^{*}$, $f(x,c_{i})$
has well-defined derivatives with respect to $x$. We find that 
\[
\frac{\partial^{2}f(x,c_{i})}{\partial x\partial c_{i}}=\frac{\phi'(x)}{\phi^{2}(x)}<0\;.
\]
It implies that $\partial_{x}f(x,c_{i})$ decreases in $c_{i}$. We
also know that $\partial_{\theta_{i}}f(\theta_{i},c_{i})=0$ by the
definition of $\theta_{i}$ and that $\partial_{\theta_{i}}^{2}f(\theta_{i},c_{i})<0$
since $\theta_{i}$ is the maximizer of $f(\theta_{i},c_{i})$. Therefore,
we have 
\[
\frac{d\partial_{\theta_{i}}f(\theta_{i},c_{i})}{dc_{i}}=\frac{\partial^{2}f(\theta_{i},c_{i})}{\partial\theta_{i}^{2}}\frac{d\theta_{i}}{dc_{i}}+\frac{\partial^{2}f(\theta_{i},c_{i})}{\partial\theta_{i}\partial c_{i}}=0\:,
\]
which implies that $d\theta_{i}/dc_{i}<0$. \eproof

\textbf{Proof of Proposition \ref{prop:alternative-V}}: Recall that
the survival probability is given by $M_{t}^{(j)}=\exp(-A_{t}^{(j)})N_{t}^{(j)}$
where $N_{t}^{(j)}=\mathbf{1}_{\{t<\tau_{E}^{(j)}\}}$. For analytical
simplicity, we introduce an associated stopping time $\hat{v}:=\inf\{t\ge0:\exp(-A_{t}^{(j)})\le\hat{l}^{(j)}\}$
and re-formulate $\hat{\tau}^{(j)}$ as $\hat{\tau}^{(j)}=\hat{v}\wedge\tau_{E}^{(j)}$;
it is straightforward to verify that this expression is equivalent
to $\inf\{t\ge0:M_{t}^{(j)}\le\hat{l}^{(j)}\}$.

To derive (\ref{eq:V-nu-alternative}), we apply Tonelli's theorem
(2.37, \citealt{Folland1999}) to the right-hand side of (\ref{eq:V-nu-new-probability-space})
and take the integral over the measure $\mathbb{L}^{(j)}$, leaving
behind the dependence on the sample path of $X$. Tonelli's theorem
is applicable because all the terms in (\ref{eq:V-nu-new-probability-space})
are bounded from below on account of the fact that $\pi(x)\ge\pi_{L}$,
$m_{i}(x)\ge(R_{r}\pi)(z^{*})$, and $g_{i}(x)\ge g_{i}(a)>-\infty$.
Thus, we fix an arbitrary sample path $X(\omega)$ for some $\omega\in\Omega$
and proceed to integrate out $\hat{l}^{(j)}$. Since $\hat{l}^{(j)}$
has a uniform distribution over $[0,1]$, the cumulative distribution
function for the stopping time $\hat{v}$ is of the form $F(v)=1-\exp(-A_{v}^{(j)})$
and $F'(v)=\lambda_{j}(X_{v})\exp(-A_{v}^{(j)})$.  Note that $F(0)=0$
and $F(\infty)=1$ because $\lim_{t\rightarrow\infty}\exp(-A_{t}^{(j)})=0$.
(For the sample paths in which $\lim_{t\rightarrow\infty}\exp(-A_{t}^{(j)})>0$,
see the remark at the end of this proof). Using the notation $\hat{u}:=\hat{\tau}^{(i)}\wedge\tau_{E}^{(j)}$,
we obtain 
\begin{align*}
\mathbb{E}_{\mathbb{L}^{(j)}}[\int_{0}^{\hat{\tau}^{(i)}\wedge\hat{\tau}^{(j)}}\pi(X_{t})e^{-rt}dt] & =\mathbb{E}_{\mathbb{L}^{(j)}}[\int_{0}^{\hat{v}\wedge\hat{u}}\pi(X_{t})e^{-rt}dt]\\
 & =\int_{0}^{\infty}(\int_{0}^{v\wedge\hat{u}}\pi(X_{t})e^{-rt}dt)F'(v)dv\\
 & =\int_{0}^{\hat{u}}[1-F(v)]\pi(X_{v})e^{-rv}dv=\int_{0}^{\hat{\tau}^{(i)}\wedge\tau_{E}^{(j)}}\pi(X_{t})e^{-rt-A_{t}^{(j)}}dt\:,
\end{align*}
where the third equality is obtained via partial integration. Next,
we obtain 
\begin{align*}
\mathbb{E}_{\mathbb{L}^{(j)}}[m_{i}(X_{\hat{\tau}^{(j)}})e^{-r\hat{\tau}^{(j)}}\mathbf{1}_{\{\hat{\tau}^{(j)}<\hat{\tau}^{(i)}\}}] & =\int_{0}^{\hat{u}}m_{i}(X_{v})e^{-rv}F'(v)dv+\mathbb{L}^{(j)}(\hat{v}>\tau_{E}^{(j)})m_{i}(X_{\tau_{E}^{(j)}})e^{-r\tau_{E}^{(j)}}\mathbf{1}_{\{\tau_{E}^{(j)}<\hat{\tau}^{(i)}\}}\\
 & =\int_{0}^{\hat{\tau}^{(i)}\wedge\tau_{E}^{(j)}}\lambda_{j}(X_{t})m_{i}(X_{t})e^{-rt-A_{t}^{(j)}}dt+\mathbf{1}_{\{\hat{\tau}^{(i)}>\tau_{E}^{(j)}\}}m_{i}(X_{\tau_{E}^{(j)}})e^{-r\tau_{E}^{(j)}-A_{\tau_{E}^{(j)}}^{(j)}}\:.
\end{align*}
Next, we note that $\mathbf{1}_{\{\hat{\tau}^{(i)}<\hat{\tau}^{(j)}\}}=\mathbf{1}_{\{\tau^{(i)}<\hat{v}\}}\mathbf{1}_{\{\tau^{(i)}<\tau_{E}^{(j)}\}}$,
so we obtain
\begin{align*}
\mathbb{E}_{\mathbb{L}^{(j)}}[g_{i}(X_{\hat{\tau}^{(i)}})e^{-r\hat{\tau}^{(i)}}\mathbf{1}_{\{\hat{\tau}^{(i)}<\hat{\tau}^{(j)}\}}] & =\mathbb{E}_{\mathbb{L}^{(j)}}[\mathbf{1}_{\{\tau^{(i)}<\hat{v}\}}]g_{i}(X_{\hat{\tau}^{(i)}})e^{-r\hat{\tau}^{(i)}}\mathbf{1}_{\{\tau^{(i)}<\tau_{E}^{(j)}\}}\\
 & =\mathbf{1}_{\{\hat{\tau}^{(i)}<\tau_{E}^{(j)}\}}g_{i}(X_{\hat{\tau}^{(i)}})e^{-r\hat{\tau}^{(i)}-A_{\hat{\tau}^{(i)}}^{(j)}}\:.
\end{align*}
Lastly, 
\begin{align*}
\mathbb{E}_{\mathbb{L}^{(j)}}[\mathbf{1}_{\{\hat{\tau}^{(i)}=\hat{\tau}^{(j)}\}}] & =\mathbb{E}_{\mathbb{L}^{(j)}}[\mathbf{1}_{\{\hat{\tau}^{(i)}=\tau_{E}^{(j)}\}}\mathbf{1}_{\{\hat{v}>\tau_{E}^{(j)}\}}+\mathbf{1}_{\{\hat{\tau}^{(i)}=\hat{v}\}}\mathbf{1}_{\{\hat{v}\le\tau_{E}^{(j)}\}}]\\
 & =\mathbf{1}_{\{\hat{\tau}^{(i)}=\tau_{E}^{(j)}\}}e^{-A_{\hat{\tau}^{(i)}}^{(j)}}+\int_{0}^{\tau_{E}^{(j)}}\mathbf{1}_{\{\hat{\tau}^{(i)}=v\}}F'(v)dv=\mathbf{1}_{\{\hat{\tau}^{(i)}=\tau_{E}^{(j)}\}}e^{-A_{\hat{\tau}^{(i)}}^{(j)}}\:.
\end{align*}
 This completes the derivation of the payoff expression (\ref{eq:V-nu-alternative}).

\emph{Remark}: Depending on the sample path $\omega\in\Omega$ of
$X$, $\exp(-A_{t}^{(j)})$ may not always converge to 0 in the limit
$t\rightarrow\infty$. However, we can always construct a process
$A_{t,n}^{(j)}$ such that $A_{t,n}^{(j)}=A_{t}^{(j)}$ for $t<n$
and $\lim_{t\rightarrow\infty}A_{t,n}^{(j)}=\infty$. If we suppose
that player $j$ adopts an investment strategy following the new additive
functional $A_{t,n}^{(j)}$, it means that player $j$ ultimately
makes investment in the time interval $(n,\infty)$ if he has not
made any in the interval $[0,n]$. With the new strategy of player
$j$, the new payoff function is identical to the old expression (\ref{eq:V-nu-new-probability-space})
in the limit $n\rightarrow\infty$ because $\lim_{n\rightarrow\infty}\mathbb{E}^{x}[(R_{r}\pi)(X_{n})e^{-rn}]=0$
and $\lim_{n\rightarrow\infty}\mathbb{E}^{x}[m_{i}(X_{n})e^{-rn}]=0$;
these transversality conditions originate from (\ref{eq:transversality})
because both $\vert(R_{r}\pi)(x)\vert$ and $\vert m_{i}(x)\vert=\vert(R_{r}\pi)(z^{*}\lor x)\vert$
are bounded by $A+B\vert x\vert$ for some $A>0$ and $B>0$ due to
Assumption \ref{assum:single-inv}. Thus, for the purpose of evaluating
the payoff function, we can safely assume that $\lim_{t\rightarrow\infty}\exp(-A_{t}^{(j)})=0$
for each sample path of $X$. \eproof

\textbf{Proof of Proposition \ref{prop:MixedStrategyMPE}}: Since
$\nu^{M}$ is already a Markov strategy profile, it suffices to prove
that $\nu_{i}^{M}$ and $\nu_{j}^{M}$ are best responses to each
other. We achieve this task in two steps: We first show that $V_{i}^{*}(x)$
dominates any $V_{i}^{(\nu_{i},\nu_{j}^{M})}(x)$ for any $\nu_{i}$,
and then we show that $V_{i}^{\nu^{M}}(x)=V_{i}^{*}(x)$.

(i) To prove that $V_{i}^{*}(x)$ dominates $V_{i}^{(\nu_{i},\nu_{j}^{M})}(x)$
for any $\nu_{i}$, it suffices to prove that $V_{i}^{*}(\cdot)$
satisfies the conditions of the verification theorem of optimal stopping
problem for player $i$ given player $j$'s strategy $\nu_{j}^{M}$.

By Proposition \ref{prop:alternative-V}, player $i$'s payoff can
be expressed as follows:
\begin{equation}
V_{i}^{(\nu_{i},\nu_{j}^{M})}(x)=\mathbb{E}^{x}\{\int_{0}^{\hat{\tau}^{(i)}}[\pi(X_{t})+\lambda_{j}(X_{t})(R_{r}\pi)(z^{*})]e^{-rt-A_{t}^{(j)}}dt+g(X_{\hat{\tau}^{(i)}})e^{-r\hat{\tau}^{(i)}-A_{\hat{\tau}^{(i)}}^{(j)}}\}\:.\label{eq:V-nu-nuM}
\end{equation}
Recall that we set $E_{1}=E_{2}=\emptyset$. Thus, as shown in equation
III.18.13 of \citet{Rogers2000}, we can define a modified characteristic
differential operator as 
\[
\hat{\mathscr{A}}_{i}:=\frac{1}{2}\sigma^{2}(x)\frac{d^{2}}{dx^{2}}+\mu(x)\frac{d}{dx}-(r+\lambda_{j}(x))=\mathscr{A}-\lambda_{j}(x)\:.
\]
Notice that player $i$'s payoff (\ref{eq:V-nu-nuM}) is equivalent
to that of an optimal stopping problem with a modified operator $\hat{\mathscr{A}}_{i}$
and a modified flow profit of $\pi(x)+\lambda_{j}(x)(R_{r}\pi)(z^{*})$.

We now verify that $V_{i}^{*}(x)$ satisfies the conditions of the
verification theorem (Theorem 10.4.1 of \citealp{Oksendal2003}) and
therefore, $V_{i}^{*}(x)\ge V_{i}^{(\nu_{i},\nu_{j}^{M})}(x)$ for
any arbitrary $\nu_{i}$. From the expression (\ref{eq:V-single}),
$V_{i}^{*}(\cdot)$ is twice continuously differentiable except at
$x=\theta(=\theta_{i})$ and continuously differentiable everywhere
\citep{Alvarez2008}. Since $V_{i}^{*}(\cdot)$ is the optimal value
function of a stopping problem, $V_{i}^{*}(x)\ge g_{i}(x)$ for all
$x\in\mathscr{I}$. Finally, note that $V_{i}^{*}(x)$ satisfies the
following differential equation:
\[
\hat{\mathscr{A}}_{i}V_{i}^{*}(x)+\pi(x)+\lambda_{j}(x)(R_{r}\pi)(z^{*})=\begin{cases}
\mathscr{A}V_{i}^{*}(x)+\pi(x) & \text{for}\:x>\theta\:,\\
\mathscr{A}V_{i}^{*}(x)+\pi(x)+\lambda_{j}(x)[(R_{r}\pi)(z^{*})-V_{i}^{*}(x)] & \text{for}\:x<\theta\:.
\end{cases}
\]
The first line of the right hand side results because $\lambda_{j}(x)=0$
for $x>\theta$. By the expression (\ref{eq:V-single}), note that
$\mathscr{A}V_{i}^{*}(x)+\pi(x)=0$ for $x>\theta$. The second line
also vanishes by the functional form of $\lambda_{j}(\cdot)$ and
because $V_{i}^{*}(x)=g(x)$ for $x\le\theta$. It follows that $\hat{\mathscr{A}}_{i}V_{i}^{*}(x)+\pi(x)+\lambda_{j}(x)(R_{r}\pi)(z^{*})=0$
for all $x$ except at $\theta$. Therefore, the conditions (i)\textendash (vi)
of Theorem 10.4.1 from \citet{Oksendal2003} are all satisfied.

(ii) It remains to prove that $V_{i}^{\nu^{M}}(x)=V_{i}^{*}(x)$.
Given that player $j$ employs $\nu_{j}^{M}$, suppose that player
$i$ adopts a strategy $\nu_{i}$ of investing at some arbitrary stopping
time $\hat{\tau}^{(i)}$ such that $X_{\hat{\tau}^{(i)}}\le\theta$.
(Note that $\hat{\tau}^{(i)}$ is not necessarily a hitting time).
Our goal is to prove that $V_{i}^{*}(x)$ has exactly the same expression
for $V_{i}^{(\nu_{i},\nu_{j}^{M})}(x)$ given by (\ref{eq:V-nu-nuM})
for any such stopping time $\hat{\tau}^{(i)}$.

We note that $\hat{\mathscr{A}}_{i}V_{i}^{*}(x)+\pi(x)+\lambda_{j}(x)(R_{r}\pi)(z^{*})=0$
for all $x$ except at $\theta$ and that $V_{i}^{*}\in C^{2}(\mathscr{I}\backslash\{\theta\})\cap C(\mathscr{I})$.
Closely following the proof of Theorem 10.4.1 of \citet{Oksendal2003},
we introduce a sequence of stopping times $\{u_{n}\}$ where $u_{n}=\min(n,\inf_{t>0}\{X_{t}\not\in G_{n}\})$
where $\{G_{n}\}$ is an increasing sequence of compact subsets of
$\mathscr{I}$ such that $\lim_{n\rightarrow\infty}G_{n}=\mathscr{I}$.
Theorem 10.4.1 of \citet{Oksendal2003} establishes that 
\[
V_{i}^{*}(x)=\mathbb{E}^{x}\{\int_{0}^{u_{n}}[\pi(X_{t})+\lambda_{j}(X_{t})(R_{r}\pi)(z^{*})]e^{-rt-A_{t}^{(j)}}dt+V_{i}^{*}(X_{u_{n}})e^{-ru_{n}-A_{u_{n}}^{(j)}}\}
\]
for any $n$. From the functional form of $V_{i}^{*}(\cdot)$ given
in (\ref{eq:V-single}) and by Assumption \ref{assum:single-inv},
we note that $\vert V_{i}^{*}(x)\vert<A\vert x\vert+B$ for some $A>0$
and $B>0$.  Hence, by the condition (\ref{eq:transversality})
that $\lim_{t\rightarrow\infty}\mathbb{E}^{x}[\vert X_{t}\vert e^{-rt}]=0$,
we have $\lim_{n\rightarrow\infty}\mathbb{E}^{x}[e^{-ru_{n}-A_{u_{n}}^{(j)}}V_{i}^{*}(X_{u_{n}})]=0$.
Furthermore, by Assumption \ref{assum:integrability}, $\pi(\cdot)$
is absolutely integrable, so we can take the limit $n\rightarrow\infty$
and use the dominated convergence theorem to obtain
\[
V_{i}^{*}(x)=\mathbb{E}^{x}\{\int_{0}^{\infty}[\pi(X_{t})+\lambda_{j}(X_{t})(R_{r}\pi)(z^{*})]e^{-rt-A_{t}^{(j)}}dt\}\:.
\]
It follows that the following Dynkin's formula holds
\[
V_{i}^{*}(x)=\mathbb{E}^{x}\{\int_{0}^{\tau}[\pi(X_{t})+\lambda_{j}(X_{t})(R_{r}\pi)(z^{*})]e^{-rt-A_{t}^{(j)}}dt+V_{i}^{*}(X_{\tau})e^{-r\tau-A_{\tau}^{(j)}}\}\:.
\]
for any arbitrary stopping time $\tau>0$.

We now set $\tau=\hat{\tau}^{(i)}$, in which case $V_{i}^{*}(X_{\tau})=g(X_{\tau})$
by the functional form of $V_{i}^{*}(\cdot)$ because $X_{\hat{\tau}^{(i)}}\le\theta$.
Then we obtain 
\[
V_{i}^{*}(x)=\mathbb{E}^{x}\{\int_{0}^{\hat{\tau}^{(i)}}[\pi(X_{t})+\lambda_{j}(X_{t})(R_{r}\pi)(z^{*})]e^{-rt-A_{t}^{(j)}}dt+g(X_{\hat{\tau}^{(i)}})e^{-r\hat{\tau}^{(i)}-A_{\hat{\tau}^{(i)}}^{(j)}}\}\:,
\]
which coincides with $V_{i}^{(\nu_{i},\nu_{j}^{M})}(x)$. Hence, $V_{i}^{(\nu_{i},\nu_{j}^{M})}(x)=V_{i}^{*}(x)$
as long as $X_{\hat{\tau}^{(i)}}\le\theta$. Finally, even if $\hat{\tau}^{(i)}$
is a probabilistic mixture of $\mathbb{F}$-stopping times, the same
argument holds as long as $X_{\hat{\tau}^{(i)}}\le\theta$. \eproof

\textbf{Proof of Proposition \ref{prop:SPE-mixed}}: In the second
stage subgame, the strategy profile reduces to the MPE of Proposition
\ref{prop:MixedStrategyMPE}. Hence, it is enough to show that $\nu^{S}$
is an SPE in the first stage.

Suppose that $X_{0}=x>\theta$ so that $\tau_{\Gamma}>0$. Our goal
is to prove that $\nu_{1}$ and $\nu_{2}$ are best responses to each
other. We first inspect player 2's best response to $\nu_{1}$. Let
$f_{i}(x;\tau)$ denote player $i$'s first-stage payoff conditional
on the current state variable $x$ associated with the investment
time $\tau\le\tau_{\Gamma}$ given that player $j$ employs strategy
$\nu_{j}$. As a matter of convention, if player $i$'s strategy is
to never invest in the first stage, we say $\tau>\tau_{\Gamma}$.

(i) We first examine player 2's best response to $\nu_{1}^{S}$. For
$\tau\le\tau_{\Gamma}$, player 2's payoff is given by
\begin{equation}
f_{2}(x;\tau)=\mathbb{E}^{x}[\int_{0}^{\tau}\pi(X_{t})e^{-rt}dt+e^{-r\tau}g(X_{\tau})]\label{eq:U2}
\end{equation}
 Note that the reward from investment exactly at time $\tau_{\Gamma}$
(at the commencement of stage 2) is still $g(X_{\tau_{\Gamma}})$
because $V_{i}^{\nu^{M}}(x)=V_{i}^{*}(x)$ by Proposition \ref{prop:MixedStrategyMPE},
where $V_{i}^{*}(x)$ is given by (\ref{eq:V-single}).  Upon comparing
(\ref{eq:U2}) to (\ref{eq:single-p}), we note that $f_{2}(x;\tau)$
coincides with the payoff from a single-player investment problem.
By Proposition \ref{prop:single-p}, therefore, we conclude that $f_{2}(x;\tau)$
is optimized when $\tau=\tau_{\Gamma}$, i.e., $\sup_{\tau\leq\tau_{\Gamma}}f_{2}(x;\tau)=f_{2}(x;\tau_{\Gamma})=V_{2}^{*}(x)$.

Next, suppose that player 2 never invests in the first stage of the
strategy profile. Note that the mixed strategy MPE commences at time
$\tau_{\Gamma}$ and that $X_{\tau_{\Gamma}}\le\theta$. Hence, the
payoff to player 2 at time $\tau_{\Gamma}$ is given by $V_{2}^{*}(X_{\tau_{\Gamma}})$
by Proposition \ref{prop:MixedStrategyMPE}. In turn, because $X_{\tau_{\Gamma}}\le\theta$,
we have $V_{2}^{*}(X_{\tau_{\Gamma}})=g(X_{\tau_{\Gamma}})$ by the
form of $V_{i}^{*}(X_{\tau_{\Gamma}})$. Thus, player 2's payoff is
given by
\begin{align*}
\mathbb{E}^{x}[\int_{0}^{\tau_{\Gamma}}\pi(X_{t})e^{-rt}dt+e^{-r\tau_{\Gamma}}V_{2}^{*}(X_{\tau_{\Gamma}})] & =\mathbb{E}^{x}[\int_{0}^{\tau_{\Gamma}}\pi(X_{t})e^{-rt}dt+e^{-r\tau_{\Gamma}}g(X_{\tau_{\Gamma}})]\\
 & =f_{2}(x;\tau_{\Gamma})=V_{2}^{*}(x)\;.
\end{align*}
Recall that $V_{2}^{*}(x)$ is player 2's payoff function if he invests
at $\tau_{\Gamma}$, which dominates the payoff from investing at
any time $\tau<\tau_{\Gamma}$ within stage 1. This implies that investment
at time $\tau_{\Gamma}$ and no investment in stage 1 are both best
responses of player 2 to $\nu_{1}^{S}$. Therefore, any probabilistic
mix of investment at $\tau_{\Gamma}$ and entering stage 2 MPE without
investing at $\tau_{\Gamma}$ is a best response for player 2. It
follows that $\nu_{2}^{S}$ is a best response to $\nu_{1}^{S}$.

(ii) Next, we examine player 1's best response to $\nu_{2}^{S}$.
We first consider the payoff $f_{1}(x;\tau)$ associated with an investment
time $\tau<\tau_{\Gamma}$. Note that $f_{1}(x;\tau)$ coincides with
the expression of $f_{2}(x;\tau)$ when $\tau<\tau_{\Gamma}$. On
the other hand, player 1's payoff associated with $\nu_{1}^{S}$ where
player 1 invests only in stage 2 is given by (\ref{eq:V11}), which
is better than $f_{2}(x;\tau)$ for any $\tau<\tau_{\Gamma}$ because
of Proposition \ref{prop:single-p} for a single-player problem. We
conclude that investment at $\tau<\tau_{\Gamma}$ is not a best response
for player 1 to $\nu_{2}^{S}$.

Lastly, we show that investment at $\tau_{\Gamma}$ is also not a
best response for player 1. The payoff from investment at $\tau_{\Gamma}$
is given by
\[
\mathbb{E}^{x}\left\{ \int_{0}^{\tau_{\Gamma}}\pi(X_{t})e^{-rt}dt+e^{-r\tau_{\Gamma}}(1-q_{2})g(X_{\tau_{\Gamma}})+\frac{1}{2}e^{-r\tau_{\Gamma}}q_{2}[(R_{r}\pi)(z^{*})+g(X_{\tau_{\Gamma}})]\right\} \:.
\]
Because $g(X_{\tau_{\Gamma}})=V_{1}^{\nu^{M}}(X_{\tau_{\Gamma}})<(R_{r}\pi)(z^{*})$
and $\mathbb{P}(e^{-r\tau_{\Gamma}}>0)>0$, we conclude that this
payoff is worse than (\ref{eq:V11}). Therefore, player 1's best response
is to never invest in the first stage. \eproof

\textbf{Proof of Theorem \ref{thm:MPE-non-existence}}: Because TSSPE
reduces to a mixed strategy MPE in the second stage subgame, the goal
of this proof is to prove the non-existence of an MPE. For this purpose,
we utilize Proposition \ref{prop:alternative-V} and the expression
(\ref{eq:V-nu-alternative}) for the payoff function. In order to
obtain Theorem \ref{thm:MPE-non-existence}, we first need to prove
five lemmas. As a first lemma, we establish that $E_{i}$ does not
overlap with $\Gamma_{j}$ or with $E_{j}$.

\begin{lemma} \label{lemm:No_overlap} If $\nu$ is an MPE, then
$E_{i}\cap\Gamma_{j}=\emptyset$ and $E_{i}\cap E_{j}=\emptyset$.
\end{lemma}

\textbf{Proof of Lemma \ref{lemm:No_overlap}}. (i) First, we prove
$E_{i}\cap E_{j}=\emptyset$. Suppose, on the contrary, that there
exists a point $x\in E_{1}\cap E_{2}$. Then the payoff to player
$i$ at $x$ is given by $\frac{1}{2}[(R_{r}\pi)(\max\{z^{*},x\})+g_{i}(x)]$.
However, if player $i$ deviates from the strategy and takes up an
alternative strategy of not investing at $X_{t}=x$, then he earns
$(R_{r}\pi)(\max\{z^{*},x\})$ instead, which is a higher payoff than
$\frac{1}{2}[(R_{r}\pi)(\max\{z^{*},x\})+g_{i}(x)]$ because $(R_{r}\pi)(\max\{z^{*},x\})>g_{i}(x)$.
This contradicts the assumption that $\nu$ is an equilibrium. We
conclude that $E_{1}\cap E_{2}=\emptyset$.

(ii) Next, we prove $E_{i}\cap\Gamma_{j}=\emptyset$. For the sake
of contradiction, we suppose that there exists $x\in E_{i}\cap\Gamma_{j}$.
Under this assumption, one of the best responses of player $j$ is
to invest immediately at $x$ because it belongs to the mixed strategy
region $\Gamma_{j}$ for player $j$. Then, given the current value
$x$ of the state variable, the reward from immediate investment is
\[
V_{j}^{\nu}(x)=\frac{1}{2}[g_{j}(x)+(R_{r}\pi)(\max\{z^{*},x\})]\:.
\]
However, if player $j$ does \emph{not} invest, player $i$ is the
only player who invests, so player $j$'s reward will be $(R_{r}\pi)(\max\{z^{*},x\})>V_{j}^{\nu}(x)$
because $(R_{r}\pi)(\max\{z^{*},x\})>g_{j}(x)$. This contradicts
the assumption that immediate investment is one of the best responses
of player $j$. \hfill{}$\square$

Lemma \ref{lemm:No_overlap} holds because a player has no incentive
to invest when his opponent invests with probability one. Next, we
establish a lemma regarding the mixed strategy regions and the equilibrium
investment rate.

\begin{lemma} \label{lemm:CommonGamma} If $\nu$ is an MPE, then
$\Gamma_{1}=\Gamma_{2}\subset(a,\min\{x_{1}^{c},x_{2}^{c}\})$, and
\begin{equation}
\lambda_{t}^{(i)}=\lambda_{i}(X_{t}):=\mathbf{1}_{\{X_{t}\in\Gamma_{i}\}}\frac{-\mathscr{A}g_{j}(X_{t})-\pi(X_{t})}{(R_{r}\pi)(z^{*})-g_{j}(X_{t})}\:.\label{eq:lambda_asymm}
\end{equation}
\end{lemma}

\textbf{Proof of Lemma \ref{lemm:CommonGamma}}. (i) We first show
that $\Gamma_{i}\cap(x_{i}^{c},b)=\emptyset$. For the purpose of
leading to a contradiction, suppose that $\Gamma_{i}\cap(x_{i}^{c},b)\neq\emptyset$.
Then because $\Gamma_{i}$ is an open set, for any $x_{0}\in\Gamma_{i}\cap(x_{i}^{c},b)$,
there is an open neighborhood $N\subset\Gamma_{i}\cap(x_{i}^{c},b)$
of $x_{0}$. We let $\tau_{N}:=\{t\ge0:X_{t}\not\in N\}$ be the escape
time from $N$.

Now, we claim that $V_{i}^{\nu}(x_{0})=g_{i}(x_{0})$. To see this,
note that whenever $X_{t}\in\Gamma_{i}$, which is player $i$'s mixed
strategy region, an immediate investment is one of the best responses
for player $i$. Recall also that $\Gamma_{i}\cap E_{j}=\emptyset$
by Lemma \ref{lemm:No_overlap}. Therefore, player $i$'s equilibrium
payoff at $x_{0}$ must be equal to the reward from immediate investment,
which is $g_{i}(x_{0})$.

Furthermore, player $i$ should earn the same payoff, $g_{i}(x_{0})$,
by investing at time $\tau_{N}$. This is because $X_{\tau_{N}}\in\Gamma_{i}$,
which implies that investing at $\tau_{N}$ must be one of the best
responses for player $i$ and should yield the same equilibrium payoff
to player $i$ as an immediate investment so that player $i$ mixes
those pure strategies in equilibrium. Therefore, we can obtain
\begin{equation}
V_{i}^{\nu}(x_{0})=g_{i}(x_{0})=\mathbb{E}^{x}\biggl[\int_{0}^{\tau_{N}}[\pi(X_{t})+\lambda_{j}(X_{t})m_{i}(X_{t})]e^{-rt-A_{t}^{(j)}}dt+g_{i}(X_{\tau_{N}})e^{-r\tau_{N}-A_{\tau_{N}}^{(j)}}\biggr]\:,\label{eq:V-nu-HJB}
\end{equation}
where the expected value in equation \eqref{eq:V-nu-HJB} is player
$i$'s payoff from investing at time $\tau_{N}$ because $\Gamma_{i}\cap E_{j}=\emptyset$
by Lemma \ref{lemm:No_overlap}.

Finally, because $N$ is an arbitrary neighborhood of $x_{0}$, we
can obtain the following HJB equation by letting $N\rightarrow\{x_{0}\}$
along with equation (7.5.1) of \citet{Oksendal2003}:
\begin{equation}
(\mathscr{A}-\lambda_{j}(x_{0}))g_{i}(x_{0})+\pi(x_{0})+\lambda_{j}(x_{0})m_{i}(x_{0})=0\:.\label{eq:A-lambda-HJB}
\end{equation}
However, because $\lambda_{j}(\cdot)\ge0$, $m_{i}(\cdot)>g_{i}(\cdot)$,
and $\mathscr{A}g_{i}(x_{0})+\pi(x_{0})>0$ for $x_{0}>x_{i}^{c}$
by Assumption \ref{assump:theta}(ii), equation \eqref{eq:A-lambda-HJB}
can never be satisfied for any $x_{0}\in\Gamma_{i}\cap(x_{i}^{c},b)$,
which leads to a contradiction. Therefore, we can conclude that $\Gamma_{i}\cap(x_{i}^{c},b)=\emptyset$.

(ii) We now prove $\Gamma_{1}=\Gamma_{2}$ and derive the functional
form of $\lambda_{i}(\cdot)$. First, choose any $x_{0}\in\Gamma_{i}$.
Then because $\Gamma_{i}\cap(x_{i}^{c},b)=\emptyset$ by part (i)
above and $\Gamma_{i}$ is an open set, it must be the case that $x_{0}<x_{i}^{c}$.
Next, similarly to part (i) above, we can pick an open neighborhood
$N\subset\Gamma_{i}\cap(a,x_{i}^{c})$ of $x_{0}$ and use an identical
argument as in part (i) to arrive at the same expression (\ref{eq:V-nu-HJB}).
Hence, if we take the limit $N\rightarrow\{x_{0}\}$ use equation
(7.5.1) of \citet{Oksendal2003}, we again obtain the expression (\ref{eq:A-lambda-HJB}),
leading to 
\[
\lambda_{j}(x_{0})=\frac{-\mathscr{A}g_{i}(x_{0})-\pi(x_{0})}{m_{i}(x_{0})-g_{i}(x_{0})}\ ,
\]
which is positive because $\mathscr{A}g_{i}(x_{0})+\pi(x_{0})<0$
for $x_{0}<x_{i}^{c}$ by Assumption \ref{assump:theta}(ii) and $m_{i}(\cdot)>g_{i}(\cdot)$.
Then because $x_{0}$ is an arbitrary point of $\Gamma_{i}$, it follows
that $\lambda_{j}(x)>0$ whenever $x\in\Gamma_{i}$, which implies
that $\Gamma_{i}\subseteq\Gamma_{j}$ because $\lambda_{j}(x)>0$
if and only if $x\in\Gamma_{j}$ by the definition of $\Gamma_{j}$.
Moreover, we can use a symmetric argument to obtain $\Gamma_{j}\subseteq\Gamma_{i}$.
Therefore, we can conclude that $\Gamma:=\Gamma_{i}=\Gamma_{j}$.
Note that this also implies that $\Gamma_{1}=\Gamma_{2}\subset(a,\min\{x_{1}^{c},x_{2}^{c}\})$
because $\Gamma_{i}\cap(x_{i}^{c},b)=\emptyset$ by part (i) above
and $\Gamma_{i}$ is an open set.

Finally, by using (a) our derivation of $\lambda_{j}(\cdot)$ above,
(b) $\Gamma=\Gamma_{i}=\Gamma_{j}\subset(a,\min\{x_{i}^{c},x_{j}^{c}\})$,
(c) $m_{i}(x)=(R_{r}\pi)(\max\{z^{*},x\})$ by definition, and (d)
$x_{i}^{c}\leq z^{*}$ by Assumption \ref{assump:theta}(ii), we can
obtain 
\[
\lambda_{i}(x)=\frac{-\mathscr{A}g_{j}(x)-\pi(x)}{(R_{r}\pi)(z^{*})-g_{j}(x)}\mathbf{1}_{\{x\in\Gamma\}}\:,
\]
which completes the derivation of the functional form of $\lambda_{i}(\cdot)$.\hfill{}$\square$

We now establish that the mixed strategy region $\Gamma:=\Gamma_{1}=\Gamma_{2}$
and the investment region $E_{i}$ must not share boundary points.
This statement is different from $\Gamma\cap E_{i}=\emptyset$, which
is already established by Lemma \ref{lemm:No_overlap}.

\begin{lemma} \label{lemm:boundary-G-E} $\partial E_{i}\cap\partial\Gamma=\emptyset$
in an MPE. \end{lemma}

\textbf{Proof of Lemma \ref{lemm:boundary-G-E}}: Suppose that $y\in\partial E_{i}\cap\partial\Gamma$.
Without loss of generality, suppose that $(c,y)\subset\Gamma$ for
some $c<y$. Then $V_{j}^{\nu}(x)=g_{j}(x)$ for all $x\in(c,y)$.
Since $y$ is at the boundary of $E_{i}$, we also have $V_{j}^{\nu}(y)=m_{j}(y)$
because player $i$ would immediately invest when $X_{t}=y$ due to
the diffusive nature of $X$ and by the definition of $E_{i}$ being
a closed set.

Suppose that player $j$ adopts an alternative strategy $\nu_{j}'$
which is identical to $\nu_{j}$ except that $(d,y)$ is his continuation
region for some $d>c$. We let $\nu'=(\nu_{i},\nu_{j}')$ denote this
new strategy profile. Because $\lambda_{i}(\cdot)$ is a continuous
function within $\Gamma$, player $j$'s payoff function under $\nu'$
within the interval $(d,y)$ can be obtained as the solution to the
boundary value problem (Chapter 9, \citealt{Oksendal2003}) which
satisfies the differential equation
\[
\mathscr{A}V_{j}^{\nu'}(x)+\pi(x)+\lambda_{i}(x)[(m_{j}(x)-V_{j}^{\nu'}(x)]=0\;,
\]
and the boundary conditions $V_{j}^{\nu'}(d)=g_{j}(d)$ and $V^{\nu'}(y)=m_{j}(y)$.
Since the solution should be continuous, there is a subinterval of
$(d,y)$ within which $V_{j}^{\nu'}(x)>g_{j}(x)$ holds because $V^{\nu'}(y)=m_{j}(y)>g_{j}(y)$.
This implies that player $j$'s payoff is higher than the equilibrium
payoff for some subinterval of $(d,y)$, and it contradicts the assumption
that $\nu$ is an equilibrium. We can construct an identical argument
for the case $c>y$. We conclude that $\partial E_{i}\cap\partial\Gamma=\emptyset$
must hold. \hfill{}$\square$

The next lemma establishes that no player invests when $X$ is above
$x_{i}^{c}$.

\begin{lemma} \label{lemm:Not_above_xc} $E_{i}\cap(x_{i}^{c},b)=\emptyset$
in an MPE. \end{lemma}

\textbf{Proof of Lemma \ref{lemm:Not_above_xc}}: Suppose, on the
contrary, there exists $y\in E_{i}\cap(x_{i}^{c},b)$. Because player
$i$ immediately invests at $y$, it follows that $V_{i}^{\nu}(y)=g_{i}(y)$.
We consider a bounded stopping time $\tau$ such that $\tau<\inf\{t\ge0:X_{t}\le x_{i}^{c}\}\wedge\tau_{E}^{(j)}$.
In other words, $X$ does not hit either $x_{i}^{c}$ or $E_{j}$
before $\tau$. Because $E_{i}\cap E_{j}=\emptyset$ and $y>x_{i}^{c}$,
we have $\tau>0$ almost surely. We let $\nu_{i}'$ denote player
$i$'s alternative strategy of investment at time $\tau$ instead.
Then the payoff associated with the new strategy profile $\nu'=(\nu_{i}',\nu_{j})$
is given by
\begin{align*}
V_{i}^{\nu'}(y) & =\mathbb{E}^{y}[\int_{0}^{\tau}\pi(X_{t})e^{-rt}dt+e^{-r\tau}g_{i}(X_{\tau})]\\
 & =\mathbb{E}^{y}[\int_{0}^{\tau}(\pi(X_{t})+\mathscr{A}g_{i}(X_{t}))e^{-rt}dt]+g_{i}(y)>g_{i}(y)=V_{i}^{\nu}(y)
\end{align*}
Note that the process $X$ is an uncontrolled state variable for $t<\tau$
because $\Gamma$ does not intersect with $(x_{i}^{c},b)$ by virtue
of Lemma \ref{lemm:CommonGamma}. The inequality is due to the fact
that $\pi(X_{t})+\mathscr{A}g_{i}(X_{t})>0$ for $t\in(0,\tau)$ and
that $\tau>0$. This contradicts the assumption that $\nu$ is an
equilibrium. We conclude that $E_{i}\cap(x_{i}^{c},b)=\emptyset$.
\hfill{}$\square$

In the final lemma before completing the proof, we establish that
$E_{1}$ and $E_{2}$ are absent in a mixed strategy MPE.

\begin{lemma} \label{lemm:EmptyE} $E_{1}=E_{2}=\emptyset$ in a
mixed strategy MPE. \end{lemma}

\textbf{Proof of Lemma \ref{lemm:EmptyE}}: For the purpose of leading
to a contradiction, suppose that $E_{1}\cup E_{2}\neq\emptyset$.
By virtue of Lemma \ref{lemm:No_overlap}, there exists an interval
$(c,d)$ which is adjacent to $\Gamma$ and $E_{i}$ for one of $i\in\{1,2\}$
but which does not overlap with $E_{j}$. Without loss of generality,
suppose that $c$ belongs to the boundary of $\Gamma$ and $d$ belongs
to the boundary of $E_{i}$. It follows that $V_{i}^{\nu}(c)=g_{i}(c)$
and $V_{i}^{\nu}(d)=g_{i}(d)$.

By virtue of Lemma \ref{lemm:CommonGamma}, $\mathscr{A}g_{i}(x)+\pi(x)<0$
for all $x\in(c,d)$. For any $y\in(c,d)$ as the initial point, we
consider $\tau_{c}:=\inf\{t\ge0:X_{t}\not\in(c,d)\}$, the time of
escape from $(c,d)$. Because $V_{i}^{\nu}(c)=g_{i}(c)$ and $V_{i}^{\nu}(d)=g_{i}(d)$,
$V_{i}^{\nu}(y)$ coincides with the reward from investment at the
time $\tau_{c}$ as follows: 
\begin{align*}
V_{i}^{\nu}(y) & =\mathbb{E}^{y}\{\int_{0}^{\tau_{c}}\pi(X_{t})e^{-rt}dt+e^{-r\tau_{c}}g_{i}(X_{\tau_{c}})\}\\
 & =\mathbb{E}^{y}\{\int_{0}^{\tau_{c}}[\pi(X_{t})+\mathscr{A}g(X_{t})]e^{-rt}dt\}+g_{i}(y)<g_{i}(y)\:.
\end{align*}
The inequality holds because $\pi(X_{t})+\mathscr{A}g(X_{t})<0$.
This implies that player $i$ is better off investing immediately
at $y$, which contradicts the assumption that $\nu$ is an equilibrium.
The contradiction results from the assumption that $E_{1}\cup E_{2}\neq\emptyset$,
and therefore, we conclude $E_{1}=E_{2}=\emptyset$. \hfill{}$\square$

We finally complete the proof of Theorem \ref{thm:MPE-non-existence}
by showing that $\sup\Gamma=\theta_{i}$; since $\sup\Gamma$ cannot
take two different values, this is possible only if $\theta_{1}=\theta_{2}$,
thus proving the theorem. To lead to a contradiction, suppose that
$\sup\Gamma:=\gamma\neq\theta_{i}$ for one player $i$. Since $\Gamma$
is union of disjoint intervals, there is some $\epsilon<\gamma$ such
that $[\epsilon,\gamma]\subset\Gamma$. Furthermore, $(\gamma,b)$
is a continuation region for both players, i.e., no player invests
in the interval $(\gamma,b)$. Since the interval $(\epsilon,\gamma)$
is a mixed strategy region, by the definition of a mixed strategy,
one possible payoff-equivalent strategy of player $i$ is to never
invest within the interval because it is one of the best responses.
This implies that we can consider $(\epsilon,b)$ a subset of a continuation
region for the purpose of computing the payoff function. Then $V_{i}^{\nu}(\cdot)$
must be a solution to a boundary value problem (Chapter 9, \citealt{Oksendal2003})
within the interval $(\epsilon,b)$ satisfying the equation (\ref{eq:A-lambda-HJB})
and the boundary conditions $V_{i}^{\nu}(\epsilon)=g_{i}(\epsilon)$
and $\lim_{x\rightarrow b}[V_{i}^{\nu}(x)-(R_{r}\pi)(x)]=0$, and
hence, it must be continuously differentiable within $(\epsilon,b)$.
In particular, the first derivative of $V_{i}^{\nu}(\cdot)$ must
be continuous at $x=\gamma$.

We formally prove this statement below. We let $\tau_{\epsilon}:=\inf\{t\ge0:X_{t}\not\in(\epsilon,b)\}$
denote the escape time from the interval $(\epsilon,b)$. By virtue
of \citet{Lon2011}, there exists a well-defined solution to the second-order
differential equation 
\[
\mathscr{A}f(x)+\pi(x)+\lambda_{j}(x)[(m_{i}(x)-f(x)]=0
\]
within a bounded interval $x\in(\epsilon,\gamma)$. Suppose that $f(x)$
is a continuously differentiable and continuous function in the interval
$[\epsilon,b)$ that is constructed to satisfy the equation above
for $x\in(\epsilon,\gamma)$ as well as $\mathscr{A}f(x)+\pi(x)=0$
for $x\in(\gamma,b)$. Furthermore, we suppose that $f(\epsilon)=g_{i}(\epsilon)$
and $\lim_{x\rightarrow b}[f(x)-(R_{r}\pi)(x)]=0$. Since $[\epsilon,\gamma]$
is compact and $f(x)=\beta\phi(x)+(R_{r}\pi)(x)$ for $x\in(\gamma,b)$
where $\beta$ is some coefficient, we find that $f(x)-(R_{r}\pi)(x)$
is a bounded function of $x$ in the interval $[\epsilon,b)$.

Following the same line of argument as in the proof part (ii) of Proposition
\ref{prop:MixedStrategyMPE}, we find that 
\[
f(x)=\mathbb{E}^{x}\{\int_{0}^{u_{n}\wedge\tau_{\epsilon}}[\pi(X_{t})+\lambda_{j}(X_{t})(R_{r}\pi)(z^{*})]e^{-rt-A_{t}^{(j)}}dt+f(X_{u_{n}\wedge\tau_{\epsilon}})e^{-ru_{n}\wedge\tau_{\epsilon}-A_{u_{n}\wedge\tau_{\epsilon}}^{(j)}}\}\:,
\]
where we employ the same notation $u_{n}$ as in the proof of Proposition
\ref{prop:MixedStrategyMPE}. Note that $f(x)$ is bounded within
the interval $[\epsilon,\gamma]$ and that $f(x)-(R_{r}\pi)(x)\rightarrow0$
in the limit $x\rightarrow\infty$. Hence, in the limit $n\rightarrow\infty$,
$f(X_{u_{n}\wedge\tau_{\epsilon}})\rightarrow f(\epsilon)$ if $X_{\tau_{\epsilon}}=\epsilon$
and $f(X_{u_{n}\wedge\tau_{\epsilon}})-(R_{r}\pi)(X_{u_{n}\wedge\tau_{\epsilon}})\rightarrow0$
if $X_{\tau_{\epsilon}}=\infty$ in case $\tau_{\epsilon}=\infty$.
(By the definition of $\tau_{\epsilon}$, these are the only two possibilities).
Since $\lim_{\tau\rightarrow\infty}\mathbb{E}^{x}[e^{-r\tau}\vert(R_{r}\pi)(X_{\tau})\vert]=0$
from the transversality assumption that we made in \eqref{eq:transversality},
we can employ the bounded convergence theorem to obtain 
\[
\lim_{n\rightarrow\infty}\mathbb{E}^{x}[f(X_{u_{n}\wedge\tau_{\epsilon}})e^{-ru_{n}\wedge\tau_{\epsilon}-A_{u_{n}\wedge\tau_{\epsilon}}^{(j)}}]=f(\epsilon)e^{-r\tau_{\epsilon}-A_{\tau_{\epsilon}}^{(j)}}\:.
\]
Using the fact that $f(\epsilon)=g_{i}(\epsilon)$ we finally obtain
the following:
\[
f(x)=\mathbb{E}^{x}\{\int_{0}^{\tau_{\epsilon}}[\pi(X_{t})+\lambda_{j}(X_{t})(R_{r}\pi)(z^{*})]e^{-rt-A_{t}^{(j)}}dt+g_{i}(X_{\tau_{\epsilon}})e^{-r\tau_{\epsilon}-A_{\tau_{\epsilon}}^{(j)}}\}\:.
\]
The right-hand-side of the equation coincides with the expression
for $V_{i}^{\nu}(x)$. Since $f(\cdot)$ was constructed to be continuously
differentiable, so is $V_{i}^{\nu}(x)$.

We now compute the discontinuity in the first derivative $(V_{i}^{\nu})'(\gamma^{+})-(V_{i}^{\nu})'(\gamma^{-})$
at $\gamma$ and see if it is non-zero. The payoff function $V_{i}^{\nu}(\cdot)$
in the interval $[\gamma,b)$ can be obtained from the solution to
the boundary value problem ( Chapter 9, \citealt{Oksendal2003}):
\[
V_{i}^{\nu}(x)=\alpha_{i}(\gamma)\phi(x)+(R_{r}\pi)(x)\:,
\]
while $V_{i}^{\nu}(x)=g_{i}(x)$ for $x\in(\epsilon,\gamma)$. From
the definition of $\alpha_{i}(\cdot)$ in (\ref{eq:alpha-i}), it
is straightforward to verify the following:
\[
(V_{i}^{\nu})'(\gamma^{+})-(V_{i}^{\nu})'(\gamma^{-})=\alpha_{i}(\gamma)\phi'(\gamma)+(R_{r}\pi)'(\gamma)-g_{i}'(\gamma)=-\phi(\gamma)\alpha_{i}'(\gamma)\;.
\]
By Assumption \ref{assump:theta}(i), $\alpha_{i}'(\gamma)=0$ only
at $\gamma=\theta_{i}$, and therefore, the first derivative of $V_{i}^{\nu}(\cdot)$
is discontinuous if $\gamma\neq\theta_{i}$. It follows that $V_{i}^{\nu}(\cdot)$
cannot be continuously differentiable if $\gamma\neq\theta_{i}$.
By Lemma \ref{lemm:CommonGamma}, this implies that $\gamma=\theta_{i}=\theta_{j}$
must hold, which is not possible if $\theta_{i}\neq\theta_{j}$. We
conclude that there is no mixed strategy MPE if $c_{1}\neq c_{2}$.
\eproof 

\textbf{Proof of Theorem \ref{thm:nu*-equilibrium}}: To complete
the proof of the theorem, we establish two lemmas. First, we show
that the strategy profile $\nu^{*}$ is an SPE if it induces the payoff
functions $F_{i}(\cdot)$ that satisfy the conditions stipulated in
Section \ref{subsec:Verification-Theorem}. Second, we prove that
the functions $F_{i}(\cdot)$, if they exist, are indeed the payoff
functions associated with the strategy profile $\nu^{*}$.

\begin{lemma} \label{lemm:verification} If the strategy profile
$\nu^{*}=(\nu_{1}^{*},\nu_{2}^{*})$ induces the payoff functions
$F_{i}(\cdot)$ satisfying the conditions given in Section \ref{subsec:Verification-Theorem},
then it is an SPE. \end{lemma}

\textbf{Proof of Lemma \ref{lemm:verification}}: Our goal is to show
that $F_{i}(\cdot)$ dominates any payoff to player $i$ given player
$j$'s strategy $\nu_{j}^{*}$.

(i) We first prove the optimality of player 2's payoff function $F_{2}(\cdot)$
given player 1's strategy of $\nu_{1}^{*}$. From the form of $\lambda_{i}(\cdot)$
in (\ref{eq:lambda-impulse}) and the fact that $V_{2}(x)=g_{2}(x)$
for $x\in\Gamma=(a,\theta^{*})$ as well as (\ref{eq:AV_i}), it is
straightforward to verify that 
\begin{equation}
\mathscr{A}V_{2}(x)+\pi(x)+\lambda_{1}(x)[V_{2}(z_{1})-V_{2}(x)]=0\:,\label{eq:V2-HJB}
\end{equation}
for $x\in\mathscr{I}\backslash\{\theta^{*}\}$.

Next, we will show that $V_{2}(\cdot)=R_{2}(\cdot)$, where $R_{2}(\cdot)$
is defined as follows:

\begin{equation}
R_{2}(x):=\mathbb{E}^{x}\{\int_{0}^{\infty}e^{-rt-A_{t}^{(1)}}[\pi(X_{t})+\lambda_{1}(X_{t})V_{2}(z_{1})]dt\}\:,\label{eq:R2}
\end{equation}
where $A_{t}^{(1)}:=\int_{0}^{t}\lambda_{1}(X_{s})ds$ and $X$ is
the \emph{uncontrolled} state variable. Note that $\lambda_{1}(\cdot)$
is a continuous function well-defined in the interval $[a,\theta^{*}]$
and zero elsewhere, so $\vert\lambda_{i}(x)\vert$ is bounded. Furthermore,
Assumption \ref{assum:integrability} ensures that the integral above
is well-defined.

To show $V_{2}(\cdot)=R_{2}(\cdot)$, we closely follow the proof
of Theorem 10.4.1 of \citet{Oksendal2003} and introduce a sequence
of stopping times $\{u_{n}\}$ defined as $u_{n}=\min(n,\inf_{t>0}\{X_{t}\not\in G_{n}\})$,
where $(G_{n})_{n\in\mathbb{N}}$ is an increasing sequence of compact
subsets of $\mathscr{I}$ such that $\lim_{n\rightarrow\infty}G_{n}=\mathscr{I}$.
Then because $V_{2}(\cdot)$ satisfies (\ref{eq:V2-HJB}),  the arguments
in the proof of Theorem 10.4.1 of \citet{Oksendal2003} establish
that 
\begin{equation}
V_{2}(x)=\mathbb{E}^{x}\{\int_{0}^{u_{n}}e^{-rt-A_{t}^{(1)}}[\pi(X_{t})+\lambda_{1}(X_{t})V_{2}(z_{1})]dt+e^{-ru_{n}-A_{u_{n}}^{(1)}}V_{2}(X_{u_{n}})\}\:,\label{eq:V2-un}
\end{equation}
for any $n$. From the functional form of $V_{2}(x)$ in (\ref{eq:Vi-explicit})
and the fact that $(R_{r}\pi)'(x)<k$ for sufficiently large $x$
by Assumption \ref{assum:single-inv}, we note that $\vert V_{2}(x)\vert<A\vert x\vert+B$
for some $A>0$ and $B>0$. By the transversality condition (\ref{eq:transversality})
and the dominated convergence theorem, we have $\lim_{n\rightarrow\infty}\mathbb{E}^{x}[e^{-ru_{n}-A_{1,u_{n}}}V_{2}(X_{u_{n}})]=0$.
Thus, if we take the limit $n\rightarrow\infty$ of (\ref{eq:V2-un}),
using the dominated convergence theorem once again, the right hand
side converges to the expression for $R_{2}(\cdot)$, so we have $V_{2}(x)=R_{2}(x)$.

Next, we let $\nu:=(\nu_{1}^{*},\nu_{2})$ denote the strategy profile
of player 1 employing $\nu_{1}^{*}$ and player 2 employing some arbitrary
strategy $\nu_{2}$. Our goal is to show that the payoff to player
2 is dominated by $V_{2}(\cdot)$. We let $\tau_{m+1}$ denote player
2's stopping time of investment after $T_{m}$ under the strategy
$\nu_{2}$, and we let $\upsilon_{m+1}$ denote player 1's time of
investment after $T_{m}$ under the strategy $\nu_{1}^{*}$. For notational
convenience, we let $\zeta^{(2)}(\tau_{m})$ denote player 2's boost
in $X^{\nu}$ at the time $\tau_{m}$. Under this notation, we note
$T_{m+1}=\tau_{m+1}\wedge\upsilon_{m+1}$.

Because $V_{2}(x)=R_{2}(x)$, we can apply the Dynkin's theorem to
$V_{2}(\cdot)$ for any arbitrary stopping time \citep{Alvarez2008}
and obtain the following:
\begin{equation}
e^{-rT_{m}}V_{2}(X_{T_{m}})-\mathbb{E}^{X_{T_{m}}}[e^{-r\tau_{m+1}^{-}-A_{m,\tau_{m+1}^{-}}^{(1)}}V_{2}(X_{\tau_{m+1}^{-}})]=\mathbb{E}^{X_{T_{m}}}\{\int_{T_{m}}^{\tau_{m+1}^{-}}e^{-rt-A_{m,t}^{(1)}}[\pi(X_{t})+\lambda_{1}(X_{t})V_{2}(z_{1})]dt\}\:,\label{eq:V2-Dynkin}
\end{equation}
where we used the notation $A_{m,t}^{(1)}:=\int_{T_{m}}^{t}\lambda_{1}(X_{s})ds$.
In the expression above, $X$ is an \emph{uncontrolled} state variable
which evolves within period $m$ with the initial value $X_{T_{m}}$.
We define $\hat{T}_{m+1}=\tau_{m+1}^{-}\wedge\upsilon_{m+1}$. Since
$\upsilon_{m+1}$ is the random time with an arrival rate of $\lambda_{1}(X_{t}^{\nu})$,
we can re-express it in terms of the \emph{controlled} state variable
$X^{\nu}$ and $\upsilon_{m+1}$ as follows:

\begin{align}
V_{2}(X_{T_{m}}^{\nu})e^{-rT_{m}}= & \hat{\mathbb{E}}_{1}^{X_{T_{m}}^{\nu}}[\int_{T_{m}}^{\tau_{m+1}^{-}\wedge\upsilon_{m+1}}e^{-rt}\pi(X_{t}^{\nu})dt+\mathbf{1}_{\{\upsilon_{m+1}<\tau_{m+1}\}}e^{-r\upsilon_{m+1}}V_{2}(z_{1})\nonumber \\
 & +\mathbf{1}_{\{\tau_{m+1}\le\upsilon_{m+1}\}}e^{-r\tau_{m+1}^{-}}V_{2}(X_{\tau_{m+1}^{-}}^{\nu})]\nonumber \\
= & \hat{\mathbb{E}}_{1}^{X_{T_{m}}^{\nu}}[\int_{T_{m}}^{\hat{T}_{m+1}}e^{-rt}\pi(X_{t}^{\nu})dt+e^{-r\hat{T}_{m+1}}V_{2}(X_{\hat{T}_{m+1}}^{\nu})]\:,\label{eq:V2-Tm}
\end{align}
where $X_{\hat{T}_{m+1}}^{\nu}=z_{1}$ if $\hat{T}_{m+1}=\upsilon_{m+1}$.
Recall that $\hat{\mathbb{E}}_{1}$ is the expectation over the measure
$\mathbb{P}\otimes_{m\in\mathbb{N}}\mathbb{L}_{m}^{(1)}$ when player
1's survival probability $M_{m,t}^{(1)}$ is given by (\ref{eq:Multiplicative-1}).
Note that the expression (\ref{eq:V2-Tm}) reduces to (\ref{eq:V2-Dynkin})
by the same argument that derived Proposition \ref{prop:alternative-V}.

From (\ref{eq:V2-Tm}), we derive the following expression:
\begin{align*}
\hat{\mathbb{E}}_{1}^{x}[V_{2}(X_{T_{m}}^{\nu})e^{-rT_{m}}-e^{-r\hat{T}_{m+1}}V_{2}(X_{\hat{T}_{m+1}}^{\nu})] & =\hat{\mathbb{E}}_{1}^{x}[\hat{\mathbb{E}}_{1}^{X_{T_{m}}^{\nu}}[\int_{T_{m}}^{\hat{T}_{m+1}}V_{2}(X_{T_{m}}^{\nu})e^{-rT_{m}}-e^{-r\hat{T}_{m+1}}V_{2}(X_{\hat{T}_{m+1}}^{\nu})]]\\
 & =\hat{\mathbb{E}}_{1}^{x}[\int_{T_{m}}^{\hat{T}_{m+1}}e^{-rt}\pi(X_{t}^{\nu})dt]\:.
\end{align*}
Then we sum the expression above for $m$ from 0 up to some finite
integer $M$: 
\[
V_{2}(x)+\hat{\mathbb{E}}_{1}^{x}\{\sum_{m=1}^{M}e^{-rT_{m}}[V_{2}(X_{T_{m}}^{\nu})-V_{2}(X_{\hat{T}_{m}}^{\nu})]-e^{-rT_{M+1}}V_{2}(X_{\hat{T}_{M+1}}^{\nu})\}=\hat{\mathbb{E}}_{1}^{x}[\int_{0}^{T_{M+1}}e^{-rt}\pi(X_{t}^{\nu})dt]\:,
\]
where we exploit the fact that $e^{-rT_{m}}=e^{-r\hat{T}_{m}}$ and
that $\int_{0}^{T_{M+1}}(...)dt=\int_{0}^{\hat{T}_{M+1}}(...)dt$.
We note that $V_{2}(X_{T_{m}}^{\nu})=V_{2}(X_{\hat{T}_{m}}^{\nu})$
whenever $\hat{T}_{m}=\upsilon_{m}$, i.e., whenever $m$-th period
begins with player 1's investment. Thus, we can re-express $V_{2}(x)$
as follows:
\begin{align}
V_{2}(x) & =\hat{\mathbb{E}}_{1}^{x}\{\sum_{m=1}^{M}e^{-rT_{m}}\mathbf{1}_{\{\tau_{m}<\upsilon_{m}\}}[V_{2}(X_{\tau_{m}^{-}}^{\nu})-V_{2}(X_{\tau_{m}}^{\nu})]+\int_{0}^{T_{M+1}}e^{-rt}\pi(X_{t}^{\nu})dt\}+\hat{\mathbb{E}}^{x}[e^{-rT_{M+1}}V_{2}(X_{\hat{T}_{M+1}}^{\nu})]\label{eq:M-equality}\\
 & \ge\hat{\mathbb{E}}_{1}^{x}\{\sum_{m=1}^{M}e^{-rT_{m}}\mathbf{1}_{\{\tau_{m}<\upsilon_{m}\}}[-k\zeta^{(2)}(\tau_{m})-c_{2}]+\int_{0}^{T_{M+1}}e^{-rt}\pi(X_{t}^{\nu})dt\}+\hat{\mathbb{E}}^{x}[e^{-rT_{M+1}}\pi_{L}/r]\:.\label{eq:M-limit}
\end{align}
The inequality holds because of the quasi-variational inequality $V_{2}(X_{\tau_{m}^{-}}^{\nu})\ge\sup_{\zeta\ge0}[V_{2}(X_{\tau_{m}^{-}}^{\nu}+\zeta)-k\zeta-c_{2}]$
from (\ref{eq:quasi-V2}) and Assumption \ref{assum:integrability}(i)
which leads to $V_{2}(x)\ge(R_{r}\pi)(x)\ge\pi_{L}/r$.

Now we take the limit $M\rightarrow\infty$ on (\ref{eq:M-limit}).
Noting that $-k\zeta^{(2)}(\tau_{m})-c_{2}<0$, we employ the monotone
convergence theorem to find that
\[
\lim_{M\rightarrow\infty}\hat{\mathbb{E}}_{1}^{x}\{\sum_{m=1}^{M}e^{-rT_{m}}\mathbf{1}_{\{\tau_{m}<\upsilon_{m}\}}[-k\zeta^{(2)}(\tau_{m})-c_{2}]\}=\hat{\mathbb{E}}_{1}^{x}\{\sum_{m=1}^{\infty}e^{-rT_{m}}\mathbf{1}_{\{\tau_{m}<\upsilon_{m}\}}[-k\zeta^{(2)}(\tau_{m})-c_{2}]\}\:.
\]
 Similarly, because $\pi(x)>\pi_{L}$ and $\lim_{m\rightarrow\infty}T_{m}=\infty$,
we can employ the monotone convergence theorem and find $\hat{\mathbb{E}}_{1}^{x}[\int_{0}^{T_{M+1}}e^{-rt}\pi(X_{t}^{\nu})dt]$
converges to $\hat{\mathbb{E}}_{1}^{x}[\int_{0}^{\infty}e^{-rt}\pi(X_{t}^{\nu})dt]$
in the same limit. Lastly, because $\lim_{M\rightarrow\infty}\hat{\mathbb{E}}_{1}^{x}[e^{-rT_{M+1}}\pi_{L}/r]=0$,
we find the following inequality:
\[
V_{2}(x)\ge\hat{\mathbb{E}}_{1}^{x}\{\sum_{m=1}^{\infty}e^{-rT_{m}}\mathbf{1}_{\{\tau_{m}<\upsilon_{m}\}}[-k\zeta^{(2)}(\tau_{m})-c_{2}]+\int_{0}^{\infty}e^{-rt}\pi(X_{t}^{\nu})dt\}\:.
\]
 Note that the right-hand-side of the inequality is the payoff to
player $2$ under the strategy profile $\nu:=(\nu_{1}^{*},\nu_{2})$.
Since $\nu_{2}$ is an arbitrary strategy of player 2, it follows
that $V_{2}(\cdot)$ dominates all possible payoffs associated with
any investment strategies of player 2 conditional on player 1's strategy
$\nu_{1}^{*}$. Since we assume that $V_{2}(\cdot)$ is player 2's
payoff from $\nu^{*}$, we conclude that $\nu_{2}^{*}$ is a best
response to $\nu_{1}^{*}$.

(ii) Next, we prove that $\nu_{1}^{*}$ is player 1's best response
to $\nu_{2}^{*}$. Using the same notational convention as in (i),
we let $\nu:=(\nu_{1},\nu_{2}^{*})$ denote the strategy profile of
player 1 employing an arbitrary investment strategy $\nu_{1}$ and
player 2 employing $\nu_{2}^{*}$. We let $\tau_{m+1}$ denote player
1's stopping time of investment after $T_{m}$ under the strategy
$\nu_{1}$, and we let $\upsilon_{m+1}$ denote player 2's time of
investment after $T_{m}$ under the strategy $\nu_{2}^{*}$. We also
let $\zeta^{(1)}(\tau_{m})$ denote player 1's boost in $X^{\nu}$
at the time $\tau_{m}$.

We now define the following: 
\[
R_{1}(x):=\mathbb{E}^{x}\{\int_{0}^{\infty}e^{-rt-A_{t}^{(2)}}[\pi(X_{t})+\lambda_{2}(X_{t})U_{1}(z_{2})]dt\}\:.
\]
Just as in (i), it is straightforward to verify that $V_{1}(\cdot)$
satisfies 
\begin{equation}
\mathscr{A}V_{1}(x)+\pi(x)+\lambda_{2}(x)[U_{1}(z_{2})-V_{1}(x)]=0\:.\label{eq:V1-HJB}
\end{equation}
Furthermore, we can utilize the same line of arguments as in (i) to
derive $V_{1}(x)=R_{1}(x)$. We now define $\hat{T}_{m+1}=\tau_{m+1}^{-}\wedge\upsilon_{m+1}$
just as in (i), where $\tau_{m+1}$ is player 1's timing of investment.
From equation (\ref{eq:AU_1}), which reduces to $\mathscr{A}U_{1}(x)+\pi(x)=0$
for $x>\theta^{*}$, we obtain the following if $X_{T_{m}}>\theta^{*}$:
\[
U_{1}(X_{T_{m}})e^{-rT_{m}}=\mathbb{E}^{X_{T_{m}}}[e^{-r\tau_{m+1}^{-}\wedge\tau_{\theta^{*}}^{m}}U_{1}(X_{\tau_{m+1}^{-}\wedge\tau_{\theta^{*}}^{m}})+\int_{T_{m}}^{\tau_{m+1}^{-}\wedge\tau_{\theta^{*}}^{m}}e^{-rt}\pi(X_{t})dt]\:.
\]
Since $U_{1}(\theta^{*})=qU_{1}(z_{2})+(1-q)V_{1}(\theta^{*})$, we
have 
\begin{align}
e^{-r\tau_{m+1}^{-}\wedge\tau_{\theta^{*}}^{m}}U_{1}(X_{\tau_{m+1}^{-}\wedge\tau_{\theta^{*}}^{m}})= & \mathbf{1}_{\{\tau_{m+1}<\tau_{\theta^{*}}^{m}\}}e^{-r\tau_{m+1}}U_{1}(X_{\tau_{m+1}^{-}})+\mathbf{1}_{\{\tau_{m+1}\ge\tau_{\theta^{*}}^{m}\}}e^{-r\tau_{\theta^{*}}^{m}}[qU_{1}(z_{2})+(1-q)V_{1}(\theta^{*})]\label{eq:U1-Tm}\\
= & \mathbf{1}_{\{\tau_{m+1}<\tau_{\theta^{*}}^{m}\}}e^{-r\tau_{m+1}}U_{1}(X_{\tau_{m+1}^{-}})+\mathbf{1}_{\{\tau_{m+1}>\tau_{\theta^{*}}^{m}\}}qe^{-r\tau_{\theta^{*}}^{m}}U_{1}(z_{2})\nonumber \\
 & +\mathbf{1}_{\{\tau_{m+1}=\tau_{\theta^{*}}^{m}\}}e^{-r\tau_{\theta^{*}}^{m}}[qU_{1}(z_{2})+(1-q)V_{1}(\theta^{*})]\nonumber \\
 & +\mathbf{1}_{\{\tau_{m+1}>\tau_{\theta^{*}}^{m}\}}(1-q)\mathbb{E}^{X_{\tau_{\theta^{*}}^{m}}}\left\{ e^{-r\tau_{m+1}^{-}-A_{m,\tau_{m+1}^{-}}^{(2)}}V_{1}(X_{\tau_{m+1}^{-}})\right.\nonumber \\
 & \left.+\int_{\tau_{\theta^{*}}^{m}}^{\tau_{m+1}^{-}}e^{-rt-A_{m,t}^{(2)}}[\pi(X_{t})+\lambda_{2}(X_{t})U_{1}(z_{2})]dt\right\} \:\nonumber 
\end{align}
where we use the notation $A_{m,t}^{(2)}:=\int_{T_{m}}^{t}\lambda_{2}(X_{t})dt$.
In the last equality, we applied the Dynkin's theorem to $V_{1}(\cdot)$.
As in (i), $X$ is an \emph{uncontrolled} state variable which evolves
within period $m$.

From the inequality $U_{1}(z_{2})\ge V_{1}(\theta^{*})$, we have
the following:
\begin{align*}
\mathbf{1}_{\{\tau_{m+1}=\tau_{\theta^{*}}^{m}\}}e^{-r\tau_{\theta^{*}}^{m}}[qU_{1}(z_{2})+(1-q)V_{1}(\theta^{*})] & \ge\mathbf{1}_{\{\tau_{m+1}=\tau_{\theta^{*}}^{m}\}}e^{-r\tau_{\theta^{*}}^{m}}[\frac{1}{2}qU_{1}(z_{2})+(1-\frac{1}{2}q)V_{1}(\theta^{*})]\\
 & =\mathbf{1}_{\{\tau_{m+1}=\tau_{\theta^{*}}^{m}\}}e^{-r\tau_{\theta^{*}}^{m}}[\frac{1}{2}qU_{1}(z_{2})+(1-\frac{1}{2}q)(U_{1}(z_{1})-k(z_{1}-\theta^{*})-c_{1})]
\end{align*}
which emulates what takes place in case $\tau_{m+1}=\tau_{\theta^{*}}^{m}$:
there is a probability of $q/2$ that player 2 invests which yields
$U_{1}(z_{2})$, and a probability of $1-q/2$ that player 1 invests,
resulting in $U_{1}(z_{1})-k(z_{1}-\theta^{*})-c_{1}$.

By the same argument as in (i), switching to the controlled state
variable $X^{\nu}$, $V_{1}(\cdot)$ satisfies the following:
\begin{align*}
V_{1}(X_{T_{m}}^{\nu})e^{-rT_{m}}= & \hat{\mathbb{E}}_{2}^{X_{T_{m}}^{\nu}}[\int_{T_{m}}^{\tau_{m+1}^{-}\wedge\upsilon_{m+1}}e^{-rt}\pi(X_{t}^{\nu})dt+\mathbf{1}_{\{\upsilon_{m+1}<\tau_{m+1}\}}e^{-r\upsilon_{m+1}}U_{1}(z_{2})\\
 & +\mathbf{1}_{\{\tau_{m+1}\le\upsilon_{m+1}\}}e^{-r\tau_{m+1}^{-}}V_{1}(X_{\tau_{m+1}^{-}}^{\nu})]\\
= & \hat{\mathbb{E}}_{2}^{X_{T_{m}}^{\nu}}[\int_{T_{m}}^{\hat{T}_{m+1}}e^{-rt}\pi(X_{t}^{\nu})dt+e^{-r\hat{T}_{m+1}}F_{1,\hat{T}_{m+1}}(X_{\hat{T}_{m+1}}^{\nu})]\:,
\end{align*}
because we can identify $F_{1,\hat{T}_{m+1}}(X_{\hat{T}_{m+1}}^{\nu})=U_{1}(z_{2})$
if $\upsilon_{m+1}<\tau_{m+1}$ and $F_{1,\hat{T}_{m+1}}(X_{\hat{T}_{m+1}}^{\nu})=V_{1}(X_{\hat{T}_{m+1}}^{\nu})=V_{1}(X_{\tau_{m+1}^{-}}^{\nu})$
if $\tau_{m+1}\le\upsilon_{m+1}$. Here, we have used the notation
$\hat{\mathbb{E}}_{2}^{X_{T_{m}}^{\nu}}$ which represents the expectation
over the measure $\mathbb{P}\otimes_{m\in\mathbb{N}}\mathbb{L}_{m}^{(2)}$
conditional on the value of $X_{T_{m}}^{\nu}$ with the survival probability
given by $M_{m,t}^{(2)}$ in (\ref{eq:Multiplicative-2}). Recall
that player 2's strategy $\nu_{2}^{*}$ is to invest at $\tau_{\theta^{*}}^{m}$
with a probability of $q$, and to invest at the rate of $\lambda_{2}(X_{t})$
for $t\in(\tau_{\theta^{*}}^{m},T_{m+1})$. Hence, from (\ref{eq:U1-Tm}),
we obtain a similar expression for $U_{1}(X_{T_{m}}^{\nu})$ for $X_{T_{m}}^{\nu}\ge\theta^{*}$:
\begin{align}
U_{1}(X_{T_{m}}^{\nu})e^{-rT_{m}}\ge & \hat{\mathbb{E}}_{2}^{X_{T_{m}}^{\nu}}\{\int_{T_{m}}^{\hat{T}_{m+1}}e^{-rt}\pi(X_{t}^{\nu})dt+e^{-r\hat{T}_{m+1}}F_{1,\hat{T}_{m+1}}(X_{\hat{T}_{m+1}}^{\nu})\mathbf{1}_{\{\tau_{m+1}\neq\tau_{\theta^{*}}^{m}\}}\nonumber \\
 & +\mathbf{1}_{\{\tau_{m+1}=\tau_{\theta^{*}}^{m}\}}e^{-r\tau_{\theta^{*}}^{m}}[(1-\chi_{m})U_{1}(z_{2})+\chi_{m}(U_{1}(z_{1})-k(z_{1}-\theta^{*})-c_{1})]\}\:.\label{eq:U1-Tm-1}
\end{align}
Here it is understood that $F_{1,\hat{T}_{m+1}}(\cdot)=U_{1}(\cdot)$
for $\hat{T}_{m+1}<\tau_{\theta^{*}}^{m}$ and $F_{1,\hat{T}_{m+1}}(\cdot)=V_{1}(\cdot)$
for $\hat{T}_{m+1}\ge\tau_{\theta^{*}}^{m}$. Furthermore, $F_{1,\hat{T}_{m+1}}(X_{\hat{T}_{m+1}}^{\nu})=U_{1}(z_{2})$
if $\upsilon_{m+1}<\tau_{m+1}$.

The last term following $\mathbf{1}_{\{\tau_{m+1}=\tau_{\theta^{*}}^{m}\}}$
requires some explanation. We define $\chi_{m}$ a Bernoulli random
variable measurable with respect to $\hat{\mathscr{\mathscr{F}}}$,
but independent of the Wiener process $W$, such that $\hat{\mathbb{P}}(\chi_{m}=0)=q/2$
and $\hat{\mathbb{P}}(\chi_{m}=1)=1-q/2$. We interpret $\chi_{m}=0$
as the event in which player 2 invests in case $\tau_{m+1}=\tau_{\theta^{*}}^{m}$.
Similarly, $\chi_{m}=1$ is the event in which player 1 invests in
case $\tau_{m+1}=\tau_{\theta^{*}}^{m}$. We note that $F_{1,T_{m+1}}(X_{m+1}^{\nu})=U_{1}(z_{2})$
in case player 2 invests at $\tau_{\theta^{*}}^{m}$, and $F_{1,T_{m+1}}(X_{m+1}^{\nu})=U_{1}(z_{1})$
in case player 1 invests at $\tau_{\theta^{*}}^{m}$. Thus, we can
re-express the last line of (\ref{eq:U1-Tm-1}) as
\[
\mathbf{1}_{\{\tau_{m+1}=\tau_{\theta^{*}}^{m}\}}e^{-r\tau_{\theta^{*}}^{m}}[F_{1,T_{m+1}}(X_{T_{m+1}}^{\nu})-\chi_{m}(k(z_{1}-\theta^{*})+c_{1})]\:.
\]

Now we can combine the two expressions to obtain the following:

\begin{align}
F_{1,T_{m}}(X_{T_{m}}^{\nu})e^{-rT_{m}}\ge & \hat{\mathbb{E}}_{2}^{X_{T_{m}}^{\nu}}[\int_{T_{m}}^{\hat{T}_{m+1}}e^{-rt}\pi(X_{t}^{\nu})dt+e^{-r\hat{T}_{m+1}}F_{1,\hat{T}_{m+1}}(X_{\hat{T}_{m+1}}^{\nu})\mathbf{1}_{\{\tau_{m+1}\neq\tau_{\theta^{*}}^{m}>T_{m}\:\text{or}\:\tau_{\theta^{*}}^{m}=T_{m}\}}\label{eq:F1-equality}\\
 & +\mathbf{1}_{\{\tau_{m+1}=\tau_{\theta^{*}}^{m}>T_{m}\}}e^{-r\tau_{\theta^{*}}^{m}}[F_{1,T_{m+1}}(X_{T_{m+1}}^{\nu})-\chi_{m}(k(z_{1}-\theta^{*})+c_{1})]\:.\nonumber 
\end{align}
which is analogous to (\ref{eq:V2-Tm}) except for the case of $\tau_{m+1}=\tau_{\theta^{*}}^{m}$.
We now define the following random variable:
\begin{align*}
\hat{F}_{\hat{T}_{m+1}}:= & F_{1,\hat{T}_{m+1}}(X_{\hat{T}_{m+1}}^{\nu})\mathbf{1}_{\{\tau_{m+1}\neq\tau_{\theta^{*}}^{m}>T_{m}\:\text{or}\:\tau_{\theta^{*}}^{m}=T_{m}\}}\\
 & +\mathbf{1}_{\{\tau_{m+1}=\tau_{\theta^{*}}^{m}>T_{m}\}}[F_{1,T_{m+1}}(X_{T_{m+1}}^{\nu})-\chi_{m}(k(z_{1}-\theta^{*})+c_{1})]\:,
\end{align*}
so that we can re-express (\ref{eq:F1-equality}) as 
\[
F_{1,T_{m}}(X_{T_{m}}^{\nu})e^{-rT_{m}}\ge\hat{\mathbb{E}}_{2}^{X_{T_{m}}^{\nu}}[\int_{T_{m}}^{\hat{T}_{m+1}}e^{-rt}\pi(X_{t}^{\nu})dt+e^{-r\hat{T}_{m+1}}\hat{F}_{\hat{T}_{m+1}}]\:.
\]
 Then we can employ the same line of logic as in (i) and sum over
the terms $\hat{\mathbb{E}}_{2}^{X_{T_{m}}^{\nu}}[e^{-rT_{m}}[\hat{F}_{\hat{T}_{m}}-F_{1,T_{m}}(X_{T_{m}}^{\nu})]]$
to arrive at the following expression:

\[
F_{1,t}(x)e^{-rt}\ge\hat{\mathbb{E}}_{2}^{x}\{\sum_{m=1}^{\infty}e^{-rT_{m}}(\mathbf{1}_{\{\tau_{m}<\upsilon_{m}\}}+\mathbf{1}_{\{\tau_{m}=\tau_{\theta^{*}}^{m-1}>T_{m-1}\}}\chi_{m})[-k\zeta^{(2)}(\tau_{m})-c_{2}]+\int_{t}^{\infty}e^{-rt}\pi(X_{t}^{\nu})dt\}\:.
\]
Note that the right-hand side is the payoff to player 1 associated
with the strategy profile $(\nu_{1},\nu_{2}^{*})$. We conclude that
$F_{1,t}(x)$ dominates the payoff from any arbitrary strategy of
player 1 against $\nu_{2}^{*}$. Thus, $\nu_{1}^{*}$ is the best
response to $\nu_{2}^{*}$. \hfill{}$\square$

Under the same set of assumptions, the next lemma establishes that
the payoff functions associated with $\nu^{*}$ are $F_{i}(\cdot)$
if it satisfies all the stipulated conditions.

\begin{lemma} \label{lemm:equil-payoff} The functions $F_{1}(\cdot)$
and $F_{2}(\cdot)$, which satisfy the conditions given in Section
\ref{subsec:Verification-Theorem}, are the payoff functions associated
with the strategy profile $\nu^{*}$. \end{lemma}

\textbf{Proof of Lemma \ref{lemm:equil-payoff}}. We prove the statement
of the lemma for $F_{2}(\cdot)$ in part (i) and for $F_{1}(\cdot)$
in part (ii).

(i) We show that $F_{2,t}(x)=V_{2}(x)$ is player 2's payoff associated
with $\nu^{*}$. Adopting the notational convention used in the proof
of Lemma \ref{lemm:verification}, we let $\upsilon_{m}$ and $\tau_{m}$
respectively denote the stopping time of investment by players 1 and
2 after time $T_{m-1}$ under the strategy profile $\nu^{*}$. We
also use the definition $\hat{T}_{m}=\tau_{m}^{-}\wedge\upsilon_{m}$.
Under strategy $\nu_{2}^{*}$, if player 2 invests at time $\tau_{m}$,
then $X_{\tau_{m}^{-}}^{\nu^{*}}\in(a,\theta^{*}]$ and $X_{\tau_{m}}^{\nu^{*}}=z_{2}$
must be satisfied, which implies that $V_{2}(X_{\tau_{m}^{-}}^{\nu^{*}})-V_{2}(X_{\tau_{m}}^{\nu^{*}})=-k(z_{2}-X_{\tau_{m}^{-}}^{\nu^{*}})-c_{2}$
from the functional form of $V_{2}(\cdot)$ because $V_{2}(x)=g_{2}(x)=V_{2}(z_{2})-k(z_{2}-x)-c_{2}$
for $x\in(a,\theta^{*}]$. From (\ref{eq:M-equality}), we then obtain
the following relation:
\[
V_{2}(x)=\hat{\mathbb{E}}^{x}\{\sum_{m=1}^{M}e^{-rT_{m}}\mathbf{1}_{\{\tau_{m}<\upsilon_{m}\}}[-k(z_{2}-X_{\tau_{m}^{-}}^{\nu^{*}})-c_{2}]+\int_{0}^{T_{M+1}}e^{-rt}\pi(X_{t}^{\nu^{*}})dt\}+\hat{\mathbb{E}}^{x}[e^{-rT_{M+1}}V_{2}(X_{\hat{T}_{M+1}}^{\nu^{*}})]\:.
\]
Here $\hat{\mathbb{E}}:=\mathbb{E}_{\mathbb{P}\otimes\mathbb{L}^{(1)}\otimes\mathbb{L}^{(2)}}$
is an expectation over $\mathbb{P}$ and the two players' mixed strategies.
Note that 
\[
\lim_{M\rightarrow\infty}\hat{\mathbb{E}}^{x}[e^{-rT_{M+1}}V_{2}(X_{\hat{T}_{M+1}}^{\nu^{*}})]=0
\]
because the times between each period $(T_{m+1}-T_{m})$ are mutually
independent random variables so that $\hat{\mathbb{E}}[e^{-r(T_{M+1}-T_{M})}]=...=\hat{\mathbb{E}}[e^{-r(T_{2}-T_{1})}]<1$.
From Assumption \ref{assum:integrability} (i) and (ii) and using
the monotone convergence theorem, we arrive at the following expression
in the limit $M\rightarrow\infty$:
\[
V_{2}(x)=\hat{\mathbb{E}}^{x}\{\sum_{m=1}^{\infty}e^{-rT_{m}}\mathbf{1}_{\{\tau_{m}<\upsilon_{m}\}}[-k(z_{2}-X_{\tau_{m}^{-}}^{\nu^{*}})-c_{2}]+\int_{0}^{\infty}e^{-rt}\pi(X_{t}^{\nu^{*}})dt\}\:.
\]
It follows that $F_{2,t}(x)=V_{2}(x)$ is the payoff induced by $\nu_{2}^{*}$
given that player 1 employs $\nu_{1}^{*}$.

(ii) Next, we prove that $F_{1}(\cdot)$ is player 1's payoff induced
by $\nu^{*}$. In this proof, we do not need to be concerned about
the event $\tau_{m}=\tau_{\theta^{*}}^{m-1}$ because player 1's strategy
$\nu_{1}^{*}$ does not allow it. Note that $F_{1,t}(\cdot)$ satisfies
(\ref{eq:F1-equality}) which is analogous to (\ref{eq:V2-Tm}). By
the same line of arguments as in (i), it follows that $F_{1}(\cdot)$
satisfies 
\[
e^{-rt}F_{1,t}(x)=\hat{\mathbb{E}}^{x}\{\sum_{m=1}^{\infty}e^{-rT_{m}}\mathbf{1}_{\{\tau_{m}<\upsilon_{m}\}}[-k(z_{2}-X_{\tau_{m}^{-}}^{\nu^{*}})-c_{2}]+\int_{t}^{\infty}e^{-rs}\pi(X_{t}^{\nu^{*}})ds\}\;,
\]
 which implies that $F_{1,t}(\cdot)$ is player 1's payoff associated
with $\nu^{*}$. \hfill{}$\square$

Finally, by virtue of Lemmas \ref{lemm:verification} and \ref{lemm:equil-payoff},
we immediately obtain Theorem \ref{thm:nu*-equilibrium}. \eproof

\textbf{Proof of Lemma \ref{lemm:I0-J0}}. (i) Because of the relationship
\[
J'(x)=(R_{r}\pi)'(x)-k-\phi'(x)I(x)-I'(x)\phi(x)=-I'(x)\phi(x),
\]
 it suffices to prove the statement for $I'(\cdot)$.

From the property of diffusion processes, there exists a constant
parameter $B$ (p. 706, \citealp{Alvarez2008}) given by
\begin{equation}
B=\frac{\psi'(x)\phi(x)-\phi'(x)\psi(x)}{S'(x)}\;,\label{eq:B-const}
\end{equation}
 and $(R_{r}f)(x):=\mathbb{E}^{x}[\int_{0}^{\infty}e^{-rt}f(X_{t})dt]$
can be expressed as follows \citep{Alvarez2008}: 
\begin{equation}
(R_{r}f)(x)=B^{-1}\phi(x)\int_{a}^{x}\psi(y)f(y)m'(y)dy+B^{-1}\psi(x)\int_{x}^{b}\phi(y)f(y)m'(y)dy\;,\label{eq:Rf}
\end{equation}
for a given integrable function $f(\cdot)$. From (\ref{eq:B-const})
and (\ref{eq:Rf}), the following expression can be derived:
\begin{align}
I'(x) & =\frac{2S'(x)}{\sigma^{2}(x)[\phi'(x)]^{2}}L(x)\;.\label{eq:I-prime}
\end{align}
From the definition of $L(\cdot)$, it is also straightforward to
obtain
\[
L'(x)=-\rho'(x)\frac{\phi'(x)}{S'(x)}\:.
\]
 By Assumption \ref{assum:rho_L}(i), we note that $L'(x)>0$ for
$x\in(a,x^{*})$ and $L'(x)<0$ for $x\in(x^{*},b)$.

Next, we examine the limiting forms of $L(\cdot)$. By Assumption
\ref{assum:rho_L}(ii), we find that $\lim_{x\rightarrow b}L(x)=0$.
Furthermore, we obtained above that $L'(x)>0$ for $x\in(a,x^{*})$
and $L'(x)<0$ for $x\in(x^{*},b)$. From Assumption \ref{assum:rho_L}(iii),
we deduce that there exists a unique point $\hat{x}\in(a,x^{*})$
such that $L(x)<0$ for all $x\in(a,\hat{x})$ and $L(x)>0$ for all
$x\in(\hat{x},b)$. From (\ref{eq:I-prime}), it follows that $I'(x)<0$
for all $x\in(a,\hat{x})$ and $I'(x)>0$ for all $x\in(\hat{x},b)$.

Lastly, we note that $I(\cdot)$ achieves its global minimum value
at $\hat{x}$. From the form of $I(\cdot)$ in (\ref{eq:I-exp}),
its minimum value must be negative, which means that $\hat{x}<z^{*}$
by Assumption \ref{assum:single-inv}. Thus, $\hat{x}<\min\{z^{*},x^{*}\}$.

(ii) This statement immediately follows from Assumption \ref{assum:single-inv},
the expression (\ref{eq:I-exp}), and the fact that $\phi'(x)<0$.
\eproof

\textbf{Proof of Theorem \ref{thm:Existence_SPE}}: Before we prove
this theorem, we state and prove the following lemma which lays out
the sufficient conditions for the verification theorem (Theorem \ref{thm:nu*-equilibrium}).
In practice, we only need to consult the following lemma to see if
a mixed strategy equilibrium exists. For the purpose of this lemma,
we continue to suppose Assumptions \ref{assum:integrability}, \ref{assum:single-inv},
\ref{assum:rho_L}, and \ref{assump:single-player-impulse}.

\begin{lemma} \label{lemm:suff-conds} Consider the strategy profile
$\nu^{*}$ given in Section \ref{subsec:Verification-Theorem} which
is characterized by $q,\theta^{*},z_{1}$, and $z_{2}$. Suppose that
$q\in(0,1)$, $\theta^{*}$, $z_{1}$, and $z_{2}$ satisfy $\theta^{*}<\hat{x}<z_{1}<z_{2}$
and the following conditions:

\begin{align}
J(z_{1})-J(\theta^{*}) & =c_{1}\label{eq:J-c1}\\
\frac{\phi(\theta^{*})[I(\theta^{*})-I(z_{1})]}{-I(z_{1})\phi(z_{2})+(R_{r}\pi)(z_{2})+I(\theta^{*})\phi(\theta^{*})-(R_{r}\pi)(\theta^{*})} & =q\;.\label{eq:q2}\\
-I(z_{1})\phi(\theta^{*})+(R_{r}\pi)(\theta^{*})-k\theta^{*} & <J(z_{1})\label{eq:OptimalU}
\end{align}
Then $\nu^{*}$ is a mixed strategy TSSPE. \end{lemma}

\textbf{Proof of Lemma \ref{lemm:suff-conds}}: Our goal is to prove
that $F_{i,t}(\cdot)$ satisfying all the conditions specified in
Section \ref{subsec:Verification-Theorem} exists if equations (\ref{eq:I-z2-th})
\textendash{} (\ref{eq:q2}) are satisfied. Then the statement of
the proposition follows from Theorem \ref{thm:nu*-equilibrium}.

(i) We first prove that the function $F_{2,t}(x)=V_{2}(x)$ characterized
in Section \ref{subsec:Verification-Theorem} exists if (\ref{eq:I-z2-th})
and (\ref{eq:J-c2}) are satisfied. Our plan is to construct a function
which satisfies all the conditions of $V_{2}(x)$.

We define a coefficient $\alpha=-I(\theta^{*})=-I(z_{2})$, which
is positive from Lemma \ref{lemm:I0-J0}(ii), and we construct a function
\[
\hat{V}_{2}(x):=\begin{cases}
\hat{g}_{2}(x) & \text{for}\:x\in(a,\theta^{*}]\\
\alpha\phi(x)+(R_{r}\pi)(x) & \text{for}\:x\in(\theta^{*},b)
\end{cases}\:,
\]
where $\hat{g}_{2}(x):=\alpha\phi(z_{2})+(R_{r}\pi)(z_{2})-k(z_{2}-x)-c_{2}=\hat{V}_{2}(z_{2})-k(z_{2}-x)-c_{2}$.

We now show that $\hat{V}_{2}(\cdot)\in C^{2}(\mathscr{I}\backslash\{\theta^{*}\})\cap C^{1}(\mathscr{I})\cap C(\mathscr{I})$.
By construction, $\hat{V}_{2}(\cdot)$ is already twice continuously
differentiable everywhere except at $\theta^{*}$, so it remains to
prove that it is continuously differentiable everywhere. From $\alpha=-I(\theta^{*})$,
we have $\lim_{x\downarrow\theta^{*}}\hat{V}_{2}(x)=-I(\theta^{*})\phi(\theta^{*})+(R_{r}\pi)(\theta^{*})$.
Similarly, from $\alpha=-I(z_{2})$, we have $\hat{V}_{2}(\theta^{*})=\hat{g}_{2}(\theta^{*})=-I(z_{2})\phi(z_{2})+(R_{r}\pi)(z_{2})-k(z_{2}-\theta^{*})-c_{2}$.
From the definition of $J(\cdot)$ in (\ref{eq:J-exp}) and the assumed
relation (\ref{eq:J-c2}), we arrive at $\lim_{x\downarrow\theta^{*}}\hat{V}_{2}(x)=\hat{V}_{2}(\theta^{*})$,
and hence, $\hat{V}_{2}(\cdot)$ is continuous everywhere in $\mathscr{I}$.
Next, we note that $\hat{V}_{2}'(\theta^{*-})=k$ and, by the definition
of $I(\cdot)$ from (\ref{eq:I-exp}), 
\[
\hat{V}_{2}'(\theta^{*+})=\alpha\phi'(\theta^{*})+(R_{r}\pi)'(\theta^{*})=-I(\theta^{*})\phi'(\theta^{*})+(R_{r}\pi)'(\theta^{*})=k\:.
\]
Thus, the first derivative of $\hat{V}_{2}(\cdot)$ is continuous
at $\theta^{*}$ and thus continuous everywhere in $\mathscr{I}$.

We also prove that $\hat{V}_{2}(\cdot)$ satisfies (\ref{eq:AV_i})
and (\ref{eq:AV_i-ineq}). By construction, $\mathscr{A}\hat{V}_{2}(x)+\pi(x)=0$
holds for $x>\theta^{*}$. It remains to show that $\mathscr{A}\hat{V}_{2}(x)+\pi(x)<0$
for $x<\theta^{*}$. From the definition of $I(\cdot)$ in (\ref{eq:I-exp})
and the fact that $\alpha=-I(z_{2})$, the first derivative of $\hat{V}_{2}(\cdot)$
is given by
\begin{equation}
\hat{V}_{2}'(x)=\begin{cases}
k & \text{for}\;x\in(a,\theta^{*})\\
k+\phi'(x)[I(x)-I(z_{2})] & \text{for}\;x\in[\theta^{*},b)
\end{cases}\;.\label{eq:V'}
\end{equation}
Because $\theta^{*}<\hat{x}<z_{2}$ by assumption, we have $I(x)<I(z_{2})$
for $x\in(\theta^{*},z_{2})$ and $I(x)>I(z_{2})$ for $x\in(z_{2},b)$.
Thus, $\hat{V}_{2}'(x)>k$ for $x\in(\theta^{*},z_{2})$. Because
$\hat{V}_{2}'(\theta^{*})=k$, this implies that $\hat{V}_{2}^{\prime\prime}(\theta^{*+})\ge0$.
Because $\mathscr{A}\hat{V}_{2}(x)+\pi(x)=0$ for $x>\theta^{*}$,
we have 
\[
\mathscr{A}\hat{V}_{2}(\theta^{*+})+\pi(\theta^{*})=\frac{1}{2}\sigma^{2}(\theta^{*})\hat{V}_{2}^{\prime\prime}(\theta^{*+})+\mu(\theta^{*})k-r\hat{V}_{2}(\theta^{*})+\pi(\theta^{*})=0\:,
\]
which implies that $\mu(\theta^{*})k-r\hat{V}_{2}(\theta^{*})+\pi(\theta^{*})\le0$.
For $x<\theta^{*}$, we have $\hat{V}_{2}(x)=k(x-\theta^{*})+\hat{V}_{2}(\theta^{*})$
by the assumed relation (\ref{eq:J-c2}), so the following holds:
\[
\mathscr{A}\hat{V}_{2}(x)+\pi(x)=\rho(x)+rk\theta^{*}-r\hat{V}_{2}(\theta^{*})\:.
\]
Recall that $\rho(x)$ is increasing for $x<\theta^{*}$ by Assumption
\ref{assum:rho_L}(i) because $\theta^{*}<\hat{x}<x^{*}$ by Lemma
\ref{lemm:I0-J0}. Therefore, we find $\mathscr{A}\hat{V}_{2}(x)+\pi(x)<\mathscr{A}\hat{V}_{2}(\theta^{*-})+\pi(\theta^{*})\le0$
for all $x<\theta^{*}$ where the last inequality holds because we
showed $\mu(\theta^{*})k-r\hat{V}_{2}(\theta^{*})+\pi(\theta^{*})\le0$
above. Therefore, $\hat{V}_{2}(\cdot)$ satisfies the differential
properties (\ref{eq:AV_i}) and (\ref{eq:AV_i-ineq}).

Next, we show that $\hat{V}_{2}(x)>(R_{r}\pi)(x)$. By construction,
it is clear that $\hat{V}_{2}(x)>(R_{r}\pi)(x)$ for $x\ge\theta^{*}$
because $\alpha=-I(\theta^{*})>0$ by Lemma \ref{lemm:I0-J0}(ii).
Recall that $\hat{V}_{2}'(x)=k$ for $x<\theta^{*}$ and that $(R_{r}\pi)'(x)>k$
for $x<\theta^{*}$ by Assumption \ref{assum:single-inv} and the
fact that $\theta^{*}<\hat{x}<z^{*}$ by Lemma \ref{lemm:I0-J0}(i).
Thus, $\hat{V}_{2}'(x)<(R_{r}\pi)'(x)$ for all $x<\theta^{*}$, which
implies that $\hat{V}_{2}(x)>(R_{r}\pi)(x)$ for $x<\theta^{*}$ as
well because both $\hat{V}_{2}(\cdot)$ and $(R_{r}\pi)(\cdot)$ are
in $C(\mathscr{I})$.

Now it remains to prove the quasi-variational inequality $\hat{V}_{2}(x)\ge\sup_{\zeta\ge0}[\hat{V}_{2}(x+\zeta)-k\zeta-c_{2}]$.
We first prove that $z_{2}$ maximizes $\hat{V}_{2}(z)-kz$. From
(\ref{eq:V'}), we have $\hat{V}_{2}'(x)=k$ for $x\in(a,\theta^{*})$,
$\hat{V}_{2}'(x)>k$ for $x\in(\theta^{*},z_{2})$, and $\hat{V}_{2}'(x)<k$
for $x\in(z_{2},b)$. It follows that $\hat{V}_{2}(z)-kz$ has a unique
global maximum exactly at $z=z_{2}$. Thus, we obtain the following
expression:
\[
\bar{V}(x):=\sup_{\zeta\ge0}[\hat{V}_{2}(x+\zeta)-k\zeta-c_{2}]=\begin{cases}
\hat{V}_{2}(z_{2})-k(z_{2}-x)-c_{2} & \text{for}\;x<z_{2}\\
\hat{V}_{2}(x)-c_{2} & \text{for}\;x\ge z_{2}
\end{cases}\:.
\]
Then the difference $\hat{V}_{2}(x)-\bar{V}(x)$ is given as follows:
\[
\hat{V}_{2}(x)-\bar{V}(x)=\begin{cases}
0 & x\in(a,\theta^{*})\\
\alpha\phi(x)+(R_{r}\pi)(x)-kx-J(z_{2})+c_{2} & x\in[\theta^{*},z_{2}]\\
c_{2} & x\in(z_{2},b)
\end{cases}\;.
\]
Recall from (\ref{eq:J-c2}) that $J(\theta^{*})-J(z_{2})+c_{2}=\alpha\phi(\theta^{*})+(R_{r}\pi)(\theta^{*})-k\theta^{*}-J(z_{2})+c_{2}=0$
and note that $\hat{V}_{2}'(x)\ge k=\bar{V}'(x)$ for $x\in(\theta^{*},z_{2})$.
Hence, we have $\alpha\phi(x)+(R_{r}\pi)(x)-kx-J(z_{2})+c_{2}\ge0$
for all $x\in(\theta^{*},z_{2})$. We conclude that $\hat{V}_{2}(x)\ge\bar{V}(x)$
for all $x$, so the variational inequality is satisfied.

(ii) We now prove that there exists $F_{1,t}(\cdot)$ which satisfies
all the conditions of Section \ref{subsec:Verification-Theorem}.
We define $\alpha:=-I(\theta^{*})$ and $\beta:=-I(z_{1})$ and construct
$\hat{F}_{1,t}(\cdot)$ as follows: $\hat{F}_{1,t}(x)=\hat{U}_{1}(x)$
if $t\in[T_{m},\tau_{\theta^{*}}^{m})$ and $\hat{F}_{1,t}(x)=\hat{V}_{1}(x)$
if $t\in[\tau_{\theta^{*}}^{m},T_{m+1})$, where $\hat{U}_{1}(\cdot)$
and $\hat{V}_{1}(\cdot)$ are given by
\begin{align*}
\hat{V}_{1}(x) & =\alpha\phi(x)+(R_{r}\pi)(x)\:,\\
\hat{U}_{1}(x) & =\beta\phi(x)+(R_{r}\pi)(x)\:,
\end{align*}
for $x\in[\theta^{*},b)$, and 
\begin{align*}
\hat{V}_{1}(x) & =\hat{g}_{1}(x):=\hat{U}_{1}(z_{1})-k(z_{1}-x)-c_{1}\:,\\
\hat{U}_{1}(x) & =\hat{V}_{1}(x)\:,
\end{align*}
for $x\in(a,\theta^{*})$. We also impose a boundary condition $\hat{U}_{1}(\theta^{*})=q\hat{U}_{1}(z_{2})+(1-q)\hat{V}_{1}(\theta^{*})$
which is derived from the condition that player 2 boosts $X$ to $z_{2}$
with a probability of $q$ when $X$ reaches $\theta^{*}$ from above.
These conditions ensure that $\hat{V}_{1}(\cdot)$ and $\hat{U}_{1}(\cdot)$
satisfy the same differential equations as $V_{1}(\cdot)$ and $U_{1}(\cdot)$
as well as the condition of investment for $x\in(a,\theta^{*})$ up
to $z_{1}$.

We first prove that $\hat{V}_{1}(\cdot)$ is continuous at $\theta^{*}$.
We compute $\hat{V}_{1}(\theta^{*+})-\hat{g}_{1}(\theta^{*-})$ and
see if it equals zero:
\begin{align*}
\hat{V}_{1}(\theta^{*+})-\hat{g}_{1}(\theta^{*-}) & =-I(\theta^{*})\phi(\theta^{*})+(R_{r}\pi)(\theta^{*})+I(z_{1})\phi(z_{1})-(R_{r}\pi)(z_{1})+k(z_{1}-\theta^{*})+c_{1}\\
 & =J(\theta^{*})-J(z_{1})+c_{1}=0\;,
\end{align*}
where the last equality is due to (\ref{eq:J-c1}). Thus, $\hat{V}_{1}(\cdot)$
is continuous everywhere. We then note that $\hat{V}_{1}(\cdot)$
coincides with $\hat{V}_{2}(\cdot)$ constructed in (i); it follows
from the fact that $\alpha=-I(\theta^{*})$, which leads to $\hat{V}_{1}(x)=\hat{V}_{2}(x)$
for $x\ge\theta^{*}$ and that $\hat{V}_{1}'(x)=\hat{V}_{2}'(x)=k$
for $x<\theta^{*}$. Hence, we immediately have $\hat{V}_{1}(\cdot)\in C^{2}(\mathscr{I}\backslash\{\theta^{*}\})\cap C^{1}(\mathscr{I})\cap C(\mathscr{I})$,
which must be satisfied by $V_{1}(\cdot)$. It also follows that $\hat{V}_{1}(x)>(R_{r}\pi)(x)$
for all $x\in\mathscr{I}$.

Next, we show that the boundary conditions of $F_{1,t}(\cdot)$ at
$\theta^{*}$ are satisfied by $\hat{F}_{1,t}(\cdot)$ if (\ref{eq:q2})
holds. It is straightforward to verify that the boundary condition
$\hat{U}_{1}(\theta^{*})=q\hat{U}_{1}(z_{2})+(1-q)\hat{V}_{1}(\theta^{*})$
is equivalent to (\ref{eq:q2}) because
\begin{equation}
q=\frac{\hat{U}_{1}(\theta^{*})-\hat{V}_{1}(\theta^{*})}{\hat{U}_{1}(z_{2})-\hat{V}_{1}(\theta^{*})}\:,\label{eq:q2-1}
\end{equation}
 leads to (\ref{eq:q2}).

Thus far, we have proved that $\hat{F}_{1}(\cdot)$ satisfies the
same differential equations and boundary conditions as $F_{1}(\cdot)$.
It remains to prove that $\hat{U}_{1}(x)>\hat{V}_{1}(x)$ for $x\ge\theta^{*}$
and that $\hat{F}_{1}(\cdot)$ satisfies the quasi-variational inequality
$\hat{V}_{1}(x)\ge\sup_{\zeta\ge0}\{\hat{U}_{1}(x+\zeta)-k\zeta-c_{1}\}$.

We first derive $\hat{U}_{1}(x)>\hat{V}_{1}(x)$ for $x\ge\theta^{*}$
from the equality $J(z_{1})-J(\theta^{*})=c_{1}$. Since $J'(x)<0$
for $x>\hat{x}$ and $J(z_{2})-J(\theta^{*})=c_{2}$, we have $\hat{x}<z_{1}<z_{2}$
because $c_{1}>c_{2}$. Furthermore, $I'(x)>0$ for $x>\hat{x}$,
so that $-I(z_{1})>-I(z_{2})=-I(\theta^{*})$, which means that $\beta>\alpha$.
It follows that $\hat{U}_{1}(x)>\hat{V}_{1}(x)$ for $x\ge\theta^{*}$.

Next, we examine the derivative of $\hat{U}_{1}(z)-kz$ for $z>\theta^{*}$:
\[
\hat{U}_{1}'(z)-k=-I(z_{1})\phi'(z)+(R_{r}\pi)'(z)-k=\phi'(z)[I(z)-I(z_{1})]\:.
\]
For $z\in(\theta^{*},z_{1})$, $I(z)$ starts out higher than $I(z_{1})$
for $z$ sufficiently close to $\theta^{*}$ because (\ref{eq:I-z2-th}),
but $I(z)$ reaches its global minimum at $\hat{x}$ and then increases
in $z$. Then, because $\phi'(\cdot)<0$, $\hat{U}_{1}'(z)-k$ crosses
zero from positive to negative across $z_{1}$. Thus, $\hat{U}_{1}(z)-kz$
has its local interior maximum at $z_{1}$. However, $\hat{U}_{1}(z)-kz$
may have its global maximum at the boundary $\theta^{*}$ because
$\hat{U}_{1}'(z)-k$ is negative for $z$ close to $\theta^{*}$.
Hence, if $\hat{U}_{1}(\theta^{*})-k\theta^{*}<\hat{U}_{1}(z_{1})-kz_{1}$,
then $\hat{U}_{1}(z)-kz$ has its global maximum at $z_{1}$. It can
be verified that $\hat{U}_{1}(\theta^{*})-k\theta^{*}<\hat{U}_{1}(z_{1})-kz_{1}$
holds if and only if (\ref{eq:OptimalU}) is satisfied. (We do not
need to consider $z<\theta^{*}$ because $\hat{U}_{1}'(x)=\hat{V}_{1}'(x)=k$
for $x<\theta^{*}$ so that $\hat{U}_{1}(z)-kz$ cannot have a global
maximum in the interval $(a,\theta^{*})$ unless $\hat{U}_{1}(\theta^{*})-k\theta^{*}\ge\hat{U}_{1}(z_{1})-kz_{1}$).

Lastly, we define $\bar{U}(x)=\sup_{\zeta\ge0}\{\hat{U}_{1}(x+\zeta)-k\zeta-c_{1}\}$
for any $x>\theta^{*}$. Then $\bar{U}(x)=\hat{U}_{1}(z_{1})-k(z_{1}-x)-c_{1}$
for $x\in(\theta^{*},z_{1}]$ and $\bar{U}(x)=\hat{U}_{1}(x)-c_{1}$
for $x>z_{1}$. Therefore,
\[
\hat{V}_{1}(x)-\bar{U}(x)=\begin{cases}
0 & x\le\theta^{*}\\
-I(\theta^{*})\phi(x)+(R_{r}\pi)(x)-kx-J(z_{1})+c_{1} & x\in(\theta^{*},z_{1})\\
-I(\theta^{*})\phi(x)+I(z_{1})\phi(x)+c_{1} & x>z_{1}
\end{cases}
\]
For $x\le\theta^{*}$, $\bar{U}(x)=\hat{g}_{1}(x)$, so it follows
that $\hat{V}_{1}(x)-\bar{U}(x)=0$. Note also that $-I(\theta^{*})\phi(\theta^{*})+(R_{r}\pi)(\theta^{*})-k\theta^{*}-J(z_{1})+c_{1}=J(\theta^{*})-J(z_{1})+c_{1}=0$
from (\ref{eq:J-c1}). We note that the first derivative of $-I(\theta^{*})\phi(x)+(R_{r}\pi)(x)-kx$
is non-negative for $x\leq z_{2}$ because this function coincides
with $\hat{V}_{2}(x)-kx$ which is proved to have non-negative first
derivative for $x\leq z_{2}$ in (i). Then because $z_{1}<z_{2}$,
$-I(\theta^{*})\phi(x)+(R_{r}\pi)(x)-kx-J(z_{1})+c_{1}$ is non-decreasing
in $(\theta^{*},z_{1})$, and hence non-negative in this interval.
Since $\bar{U}(x)$ is continuous at $z_{1}$, it follows that $-I(\theta^{*})\phi(z_{1})+I(z_{1})\phi(z_{1})+c_{1}\ge0$.
We also note that $I(z_{1})<I(\theta^{*})$ because (\ref{eq:I-z2-th}),
$\theta^{*}<\hat{x}<z_{1}<z_{2}$, and Lemma \ref{lemm:I0-J0}, and
that $\phi(\cdot)$ is decreasing. We conclude that $[I(z_{1})-I(\theta^{*})]\phi(x)+c_{1}$
is non-negative for all $x>z_{1}$, and hence, the quasi-variational
inequality is satisfied.

Lastly, note that we already proved above that $\hat{U}_{1}(x)>\hat{V}_{1}(x)$
for $x\ge\theta^{*}$ and we have $\hat{U}_{1}(x)=\hat{V}_{1}(x)$
for $x<\theta^{*}$ by construction. Furthermore, the inequality $\hat{U}_{1}(z_{2})-\hat{V}_{1}(\theta^{*})>0$
must be satisfied because otherwise $q$ in (\ref{eq:q2-1}) cannot
be a positive quantity. Hence, the requisite inequality $\hat{U}_{1}(z_{2})\ge\hat{V}_{1}(\theta^{*})$
is satisfied. We conclude that $\hat{F}_{1,t}(\cdot)$ satisfies all
the conditions of $F_{1,t}(\cdot)$. \hfill{}$\square$

Now we proceed to prove the theorem. It suffices to prove that there
exist $z_{1}$ and $q$ such that $\theta^{*}<\hat{x}<z_{1}<z_{2}$
which satisfy (i) (\ref{eq:J-c1}) and (\ref{eq:q2}), and (ii) (\ref{eq:OptimalU})
if $c_{1}$ is sufficiently close to $c_{2}$.

(i) We prove that given $\theta^{*}$ and $z_{2}$ satisfying (\ref{eq:I-z2-th})
and (\ref{eq:J-c2}), there exist $z_{1}$ and $q$ which satisfy
(\ref{eq:J-c1}) and (\ref{eq:q2}). Recall that $J(\cdot)$ is globally
maximized at $\hat{x}$. If $c_{1}$ is sufficiently close to $c_{2}$
so that $J(\hat{x})-J(\theta^{*})\ge c_{1}$, then there exists a
unique value of $z_{1}\in[\hat{x},z_{2})$ which satisfies $J(z_{1})=J(\theta^{*})+c_{1}$
because $J(\cdot)$ is monotone in the interval $(\hat{x},b)$. (If
$c_{1}>J(\hat{x})-J(\theta^{*})$, then there is no such $z_{1}$
which satisfies $J(z_{1})=J(\theta^{*})+c_{1}$).

Next, from the functional forms $\hat{U}_{1}(\cdot)$ and $\hat{V}_{i}(\cdot)$
defined in the proof of Lemma \ref{lemm:suff-conds}, (\ref{eq:q2})
is equivalent to the relation $q=[\hat{U}_{1}(\theta^{*})-\hat{V}_{1}(\theta^{*})]/[\hat{U}_{1}(z_{2})-\hat{V}_{1}(\theta^{*})]$
from (\ref{eq:q2-1}). We only need to prove that this quantity $q$
is within the interval $(0,1)$. First, because $-I(\theta^{*})<-I(z_{1})$,
we have $\hat{U}_{1}(x)>\hat{V}_{1}(x)$ for all $x\ge\theta^{*}$,
which implies that $\hat{U}_{1}(\theta^{*})-\hat{V}_{1}(\theta^{*})>0$.
Furthermore, from (\ref{eq:J-c1}), for $c_{1}$ sufficiently close
to $c_{2}$, $z_{1}$ can be made sufficiently close to $z_{2}$ in
which case $\hat{U}_{1}(x)=-I(z_{1})\phi(x)+(R_{r}\pi)(x)$ and $\hat{V}_{1}(x)=-I(z_{2})\phi(x)+(R_{r}\pi)(x)$
can be made arbitrarily close to each other in the interval $(\theta^{*},z_{2})$.
More specifically, $\hat{U}_{1}(x)$ uniformly converges to $\hat{V}_{1}(x)$
in the compact interval $[\theta^{*},z_{2}]$ in the limit $c_{1}\downarrow c_{2}$.
Since $\hat{V}_{1}(\cdot)$ is an increasing function in the interval
$(\theta^{*},z_{2})$, we have $\hat{V}_{1}(z_{2})>\hat{V}_{1}(\theta^{*})$,
which implies that $\hat{U}_{1}(z_{2})>\hat{U}_{1}(\theta^{*})$ for
$c_{1}$ sufficiently close to $c_{2}$. Because $\hat{U}_{1}(z_{2})>\hat{U}_{1}(\theta^{*})>\hat{V}_{1}(\theta^{*})$,
we conclude that $q\in(0,1)$ for $c_{1}$ sufficiently close to $c_{2}$.

(ii) Next, we prove that (\ref{eq:OptimalU}) holds if $c_{1}$ is
sufficiently close to $c_{2}$. From the functional form of $\hat{U}_{1}(\cdot)$
in the proof of Lemma \ref{lemm:suff-conds}, the inequality $\hat{U}_{1}(\theta^{*})-k\theta^{*}<\hat{U}_{1}(z_{1})-kz_{1}$
is equivalent to (\ref{eq:OptimalU}). We already established in (i)
above that $z_{1}\rightarrow z_{2}$ in the limit $c_{1}\rightarrow c_{2}$.
Furthermore, since $\hat{U}_{1}(x)$ uniformly converges to $\hat{V}_{1}(x)$
in the compact interval $[\theta^{*},z_{2}]$ in the limit $c_{1}\rightarrow c_{2}$,
we arrive at $\hat{U}_{1}(\theta^{*})-k\theta^{*}<\hat{U}_{1}(z_{1})-kz_{1}$
from the inequality $\hat{V}_{1}(\theta^{*})-k\theta^{*}<\hat{V}_{1}(z_{2})-kz_{2}$.
\eproof

\textbf{Proof of Proposition \ref{prop:NumericalEx}}: From (\ref{eq:I-exp})
and (\ref{eq:J-exp}), we have
\begin{align*}
I(x) & =\frac{1}{\gamma_{-}}(\frac{\alpha x^{\alpha-\gamma_{-}}}{r-\delta(\alpha)}-kx^{1-\gamma_{-}})\\
J(x) & =(1-\frac{\alpha}{\gamma_{-}})\frac{x^{\alpha}}{r-\delta(\alpha)}-kx(1-\frac{1}{\gamma_{-}})\:.
\end{align*}
It follows that we have $I(0)=J(0)$ and $J(z^{*})>0$. In addition,
by the definition of $z^{*}$, we have$I(x)<0$ for $x<z^{*}$ and
$I(x)>0$ for $x>z^{*}$.

By Lemma \ref{lemm:I0-J0}, $I(\cdot)$ and $J(\cdot)$ are strictly
monotone functions in each of the intervals $(0,\hat{x})$ and $(\hat{x},z^{*})$.
Then we can define bijective continuous functions $I_{1}:(0,\hat{x})\mapsto(I(\hat{x}),0)$
and $I_{2}:(\hat{x},z^{*})\mapsto(I(\hat{x}),0)$ such that $I_{1}(x)=I(x)$
for $x\in(0,\hat{x})$ and $I_{2}(x)=I(x)$ for $(\hat{x},z^{*})$.
We now define another continuous function $Z(x):=I_{2}^{-1}(I_{1}(x))$
which maps $(0,\hat{x})$ to $(\hat{x},z^{*})$, where $I_{2}^{-1}$
is the inverse function of $I_{2}$. Thus defined, $x$ and $Z(x)$
satisfies $I(x)=I(Z(x))$ where $x\in(0,\hat{x})$ and $Z(x)\in(\hat{x},z^{*})$.

Next, we define yet another continuous function $C(x):=J(Z(x))-J(x)$.
Since $\hat{x}$ is the point of the global minimum of $I(\cdot)$,
$\lim_{x\rightarrow\hat{x}}Z(x)=\hat{x}$. Hence, $\lim_{x\rightarrow\hat{x}}C(x)=0$.
Furthermore, $\lim_{x\rightarrow0}Z(x)=z^{*}$, so $\lim_{x\rightarrow0}C(x)=J(z^{*})-J(0)=J(z^{*})$.
Since $C(\cdot)$ continuously varies from 0 to $J(z^{*})$ as $x$
varies from $\hat{x}$ to $0$, there is a value $\theta^{*}\in(0,\hat{x})$
that satisfies $C(\theta^{*})=c_{2}$ as long as $c_{2}\in(0,J(z^{*}))$.
Thus, if we set $z_{2}=Z(\theta^{*})$, Assumption \ref{assump:single-player-impulse}
is satisfied. \eproof

\textbf{Proof of Lemma \ref{lemm:1-player-impulse}}: We only need
to verify that the function $V_{s}(x;c)$ satisfies the conditions
of Lemma 2.1 from \citet{Alvarez2008}. Note that $V_{s}(x;c)$ coincides
with the function $\hat{V}_{2}(x)$ defined in the proof of Lemma
\ref{lemm:suff-conds} if we identify $c=c_{2}$ and $z=z_{2}$. In
the proof of Lemma \ref{lemm:suff-conds}, it is established that
$\hat{V}_{2}(\cdot)$ satisfies all the conditions of Lemma 2.1 from
\citet{Alvarez2008}. In particular, we proved that (i) $(\hat{V}_{2}-(R_{r}\pi))(x)>0$
(non-negativity), (ii) $\mathscr{A}(\hat{V}_{2}-(R_{r}\pi))(x)=\mathscr{A}\hat{V}_{2}(x)+\pi(x)\leq0$
for all $x\neq\theta^{*}$ ($r$-superharmonicity for $X$) , and
(iii) $\hat{V}_{2}(x)\ge\sup_{\zeta\ge0}[\hat{V}_{2}(x+\zeta)-k\zeta-c_{2}]$
(quasi-variational inequality).  Finally, from its functional form,
it is evident that $V_{s}(x;c)$ is the payoff from the stated policy
of a single player. \eproof

\textbf{Proof of Proposition \ref{prop:pure-mixed}}: (i) We first
construct a pure strategy MPE. We let $\nu_{1}^{P}$ be player $1$'s
strategy of no investment at all and $\nu_{2}^{P}$ be player $2$'s
strategy of investment at $\tau_{\theta^{*}}=\inf\{t\ge0:X_{t}^{\nu^{P}}\le\theta^{*}\}$
to boost $X^{\nu^{P}}$ up to $z_{2}$. We only need to prove that
the two players' strategies are best responses to each other. Since
the strategy profile $\nu^{P}$ is Markov, we can let $V_{i}^{\nu^{P}}(x)$
denote the payoff function to players $i$ with the current state
variable $x$.

First, we prove that $\nu_{1}^{P}$ is the best response to $\nu_{2}^{P}$.
Note that $\mathscr{A}V_{1}^{\nu^{P}}(x)+\pi(x)=0$ for $x>\theta^{*}$
and $V_{1}^{\nu^{P}}(x)=V_{1}^{\nu^{P}}(z_{2})$ for $x\le\theta^{*}$
with the boundary condition $V_{1}^{\nu^{P}}(\theta^{*})=V_{1}^{\nu^{P}}(z_{2})$
from player 2's strategy of boosting $X^{\nu^{P}}$ up to $z_{2}$
at the hitting time of $\theta^{*}$. Based on the differential equation
and the boundary condition, we immediately obtain the following:
\begin{align*}
V_{1}^{\nu^{P}}(x) & =\beta\phi(x)+(R_{r}\pi)(x)\;\text{for}\;x\ge\theta^{*}\\
 & =V_{1}^{\nu^{P}}(z_{2})\;\text{for}\;x<\theta^{*}.,
\end{align*}
where $\beta$ is given in (\ref{eq:pure-suff}).

Next, we prove the variational inequality $V_{1}^{\nu^{P}}(x)\ge\sup_{z\ge x}\{V_{1}^{\nu^{P}}(z)-k(z-x)-c_{1}\}$.
We note that this inequality is always satisfied if $(V_{1}^{\nu^{P}})'(x)-k<0$
for all $x>\theta^{*}$. From the definition of $I(\cdot)$, we find
that the inequality $(V_{1}^{\nu^{P}})'(x)-k<0$ can be re-expressed
as $\beta>-I(x)$ for all $x>\theta^{*}$. This inequality is always
satisfied if $\beta>-I(\hat{x})$ because $\hat{x}$ is the global
minimizer of $I(\cdot)$. Thus, (\ref{eq:pure-suff}) is the sufficient
condition for the variational inequality.

In sum, $V_{1}^{\nu^{P}}(x)$ satisfies the variational inequality
and the differential equation $\mathscr{A}V_{1}^{\nu^{P}}(x)+\pi(x)=0$
for $x>\theta^{*}$ as well as the inequality $V_{1}^{\nu^{P}}(x)>(R_{r}\pi)(x)$
because $\beta>0$. Then we can employ the same line of argument used
to prove Lemma \ref{lemm:verification} (also similarly to Lemma 2.1
of \citealp{Alvarez2008}) and conclude that $V_{1}^{\nu^{P}}(x)$
dominates player 1's payoff from any arbitrary strategy against $\nu_{2}^{P}$.
Therefore, $\nu_{1}^{P}$ is the best response to $\nu_{2}^{P}$.

Next, we show that $\nu_{2}^{P}$ is the best response to $\nu_{1}^{P}$.
If player 1 never invests, then the investment decision for player
2 reduces to that of a single player. By Lemma \ref{lemm:1-player-impulse},
$\nu_{2}^{P}$ is the optimal policy for a single player with cost
$c_{2}$. Therefore, $\nu_{2}^{P}$ is the best response to $\nu_{1}^{P}$.

(ii) Let $\nu^{*}$ denote the mixed strategy SPE obtained in Theorem
\ref{thm:Existence_SPE}. We show that $\nu^{P}$ is more efficient
than $\nu^{*}$. We first note that $V_{2}^{\nu^{P}}(x)=V_{2}^{\nu^{*}}(x)$
for all $x$, so it remains to show that $V_{1}^{\nu^{P}}(x)>F_{1,t}(x)$
for all $x$ and $t$, where $F_{1}(\cdot)$ is given by (\ref{eq:F1}).

First, we consider $F_{1,t}(x)=V_{1}(x)$ for $t\ge\tau_{\theta^{*}}^{m}$.
Recall that $V_{1}(x)=V_{2}(x)=V_{s}(x;c_{2})=-I(\theta^{*})\phi(x)+(R_{r}\pi)(x)$
for $x\ge\theta^{*}$. In (i), we already established that $V_{1}^{\nu^{P}}(x)>V_{s}(x;c_{2})$
because $\beta>-I(\hat{x})$, so it follows that $V_{1}^{\nu^{P}}(x)>F_{1}(x;t)$
if $t>\tau_{\theta^{*}}$.

Next, we consider $F_{1,t}(x)=U_{1}(x)$ for $t<\tau_{\theta^{*}}^{m}$.
Because $V_{1}(x)=V_{s}(x;c_{2})$ and $U_{1}(\theta^{*})=qU_{1}(z_{2})+(1-q)V_{1}(\theta^{*})$,
we can construct $U_{1}(\cdot)$ as a payoff to a single player with
a policy characterized by $\theta^{*}$ and $z_{2}$ with a dynamic
upfront cost. To be specific, his upfront cost is initially zero,
but it permanently increases to $c_{2}$ with a probability of $1-q$
at each investment event; once it increases to $c_{2}$, the cost
does not change anymore. In comparison, $V_{1}^{\nu^{P}}(\cdot)$
can be obtained as a payoff associated with a policy of $\theta^{*}$
and $z_{2}$ but with zero upfront cost. Therefore, $V_{1}^{\nu^{P}}(x)>U_{1}(x)$.
\eproof
\noindent \begin{flushleft}
\textsf{\textbf{\large{}References}}{\small{}\begin{btSect}[INFORMS2011]{ImpulseGame}
\btPrintCited
\end{btSect}
}{\small\par}
\par\end{flushleft}

\end{btUnit}
\end{document}